%% file: iga.tex
\documentclass[preprint,10pt]{elsarticle}

\usepackage[hyphens]{url}
\usepackage{graphicx}
\usepackage{subfig}
\usepackage{epsfig}
\usepackage{psfrag}
\usepackage{hyperref}
\usepackage{amsmath,amssymb,amsfonts,amsthm,stmaryrd}
\usepackage{dsfont}
\usepackage{eucal}
\usepackage{mathrsfs}
\usepackage{listings}    
\usepackage{color}
\usepackage{float}
\usepackage{natbib}
\usepackage{booktabs}
\usepackage{tabularx}
\usepackage{setspace}
\usepackage{tikz}
\usepackage{enumerate}
\usetikzlibrary{calc}
\usetikzlibrary{shapes}

\input{pre_commands.tex}

\newcommand{\grbf}[1]{\mbox{\boldmath{$#1$}}}

\floatstyle{ruled}
\newfloat{Program}{thp}{lop}
\floatname{Program}{Program}

\floatstyle{ruled}
\newfloat{Fbox}{thp}{lop}
\floatname{Fbox}{Box}

\newcolumntype{C}{>{\centering\arraybackslash}X}

\definecolor{darkgray}{rgb}{0.95,0.95,0.95}
\definecolor{mygreen}{rgb}{0,0.6,0}

\lstdefinestyle{Matlab}
{
 basicstyle=\footnotesize, numbers=none, numberstyle=\tiny,%
 showstringspaces=false, language=Matlab, escapechar=|,frame=tb,%
 commentstyle=\color{mygreen}
}

\lstdefinestyle{Matlab-num}
{
 basicstyle=\footnotesize, numbers=left, numberstyle=\tiny,%
 showstringspaces=false, language=Matlab, escapechar=|,frame=tb,%
commentstyle=\color{mygreen}
}

\newcommand{\tty}[1]{\textnormal{\texttt{#1}}}
\newcommand{\sym}[1]{\textnormal{\textit{#1}}}

\lstnewenvironment{snippet}[1][]
{
 \lstset{style=Matlab, xleftmargin=5mm, gobble=4, #1}
}
{}

\lstnewenvironment{snippet1}[1][]
{
 \lstset{style=Matlab-num, xleftmargin=5mm, gobble=4, #1}
}
{}

\makeatletter
\newdimen\@myBoxHeight%
\newdimen\@myBoxDepth%
\newdimen\@myBoxWidth%
\newdimen\@myBoxSize%
\newcommand{\SquareBox}[2][]{%
    \settoheight{\@myBoxHeight}{#2}
    \settodepth{\@myBoxDepth}{#2}
    \settowidth{\@myBoxWidth}{#2}
    \pgfmathsetlength{\@myBoxSize}{max(\@myBoxWidth,(\@myBoxHeight+\@myBoxDepth))}%
    \tikz \node [shape=rectangle, shape aspect=1,draw=red,inner sep=2\pgflinewidth, minimum size=\@myBoxSize,#1] {#2};%
}%
\makeatother

\begin{document}


\definecolor{MyDarkBlue}{rgb}{1, 0.9, 1}
\lstset{language=Matlab,
       basicstyle=\footnotesize,
       commentstyle=\itshape,
       stringstyle=\ttfamily,
       showstringspaces=false,
       tabsize=2}
\lstdefinestyle{commentstyle}{color=\color{green}}

\theoremstyle{remark}
\newtheorem{thm}{Theorem}[section]
\newtheorem{rmk}[thm]{Remark}


\definecolor{red}{gray}{0}
\definecolor{blue}{gray}{0}

\doublespacing

\begin{frontmatter}

\title{Isogeometric analysis: an overview and computer implementation aspects}

\author[cardiff]{Vinh Phu Nguyen  \corref{cor1}\fnref{fn1}}
\author[cardiff]{St\'{e}phane P.A. Bordas \fnref{fn3}}
\author[weimar]{Timon Rabczuk \fnref{fn4}}

\cortext[cor1]{Corresponding author}

\address[cardiff]{School of Engineering, Institute of Mechanics and Advanced
Materials, Cardiff University, Queen's Buildings, The Parade, Cardiff \\
CF24 3AA}
\address[weimar]{Institute of Structural Mechanics, Bauhaus-Universit\"{a}t
Weimar, Marienstra\ss{}e 15 99423 Weimar}

\fntext[fn1]{\url nguyenpv@cardiff.ac.uk}
\fntext[fn3]{\url stephane.bordas@alum.northwestern.edu}
\fntext[fn4]{\url timon.rabczuk@uni-weimar.de}

\begin{abstract}
	Isogeometric analysis (IGA) represents a recently developed technology in computational mechanics that offers 
	the possibility of integrating methods for analysis and
	 Computer Aided Design (CAD) into a single, unified process. The implications to practical engineering design scenarios are profound, since the time taken from design to analysis is greatly reduced, leading to dramatic gains in efficiency. 
	The tight coupling of CAD and analysis within IGA requires knowledge from both
	fields and it is one of the goals of the present paper to outline much of the commonly used notation. 
	In this manuscript, through
	a clear and simple Matlab\textsuperscript{\textregistered{}} implementation, we present an
	introduction to IGA applied to the Finite Element (FE) method and related computer 
	implementation aspects. 
	Furthermore, implementation of the extended IGA which incorporates
	enrichment functions through the partition of unity method (PUM)  is
	also presented, where several examples for both two-dimensional and three-dimensional fracture are 
illustrated. The
	open source Matlab\textsuperscript{\textregistered{}}  code which accompanies the present
	paper can be applied to one, two and three-dimensional problems for
	linear elasticity, linear elastic fracture mechanics, structural mechanics (beams/plates/shells including
              large displacements and rotations)
        and Poisson problems with or without enrichment. The
        B\'{e}zier extraction concept that        allows FE analysis
        to be performed efficiently on T-spline geometries is also
        incorporated. The article includes a summary of recent trends
        and developments within the field of IGA.
\end{abstract}


\end{frontmatter}


\section{Introduction}

\subsection{Underlying concepts of isogeometric analysis}

The predominant technology that is used by CAD to represent complex
geometries is the Non-Uniform Rational B-spline (NURBS). This allows certain
geometries to be represented exactly that are only approximated by
polynomial functions, including conic and circular sections. There is a
vast array of literature focused on NURBS
(e.g. \cite{piegl_book}, \cite{Rogers2001}) and as a result of several decades
of research, many efficient computer algorithms exist for their fast evaluation and
refinement. The key concept outlined by Hughes et al.
\cite{hughes_isogeometric_2005} was to employ NURBS not only as a geometry discretisation technology, but also as a discretisation tool for analysis, attributing such methods to the field of `Isogeometric Analysis'
(IGA). Since this seminal paper, a monograph
dedicated entirely to IGA has been
published \cite{cottrel_book_2009} and applications can now be found
in several fields including structural mechanics, solid
mechanics, fluid mechanics and contact mechanics. 
We give in this section an overview of some of these recent developments 
while outlining the benefits and present shortcomings of IGA.
It should be emphasized that the idea of using CAD technologies in finite elements dates back
at least to \cite{NME:NME292,Kagan2000539} where B-splines were used as shape functions in FEM. In addition, similar methods which adopt
subdivision surfaces have been used to model shells \cite{Cirak_2000}. We also review
some of the recent attempts at simplifying the CAD-FEA integration by separating boundary and domain discretisations.

\subsection{Applications}

In  contact formulations using conventional geometry discretisations, the presence of faceted surfaces can lead to jumps and
oscillations in traction responses unless very fine meshes are used. The benefits of
using NURBS over such an approach are evident, since smooth contact
surface are obtained, leading to more physically accurate contact stresses. Recent work in this area includes 
\cite{temizer_contact_2011,jia_isogeometric_2011,temizer_three-dimensional_2012,
de_lorenzis_large_2011,Matzen201327}.

IGA has also shown advantages over traditional approaches in the context of optimisation problems 
\cite{wall_isogeometric_2008,manh_isogeometric_2011,qian_isogeometric_2011,xiaoping_full_2010} where the tight coupling with CAD
models offers an extremely attractive approach for industrial
applications.  Another attractive class of methods include those that require only a boundary discretisation, creating a truly direct coupling with CAD. Isogeometric boundary element methods for elastostatic analysis were presented in 
\cite{simpson_two-dimensional_2012,Scott2013197}, demonstrating that mesh generation can be completely circumvented by using CAD discretisations for analysis. 

Shell and plate problems are another field where IGA has demonstrated 
compelling benefits over conventional approaches
\cite{benson_isogeometric_2010,kiendl_isogeometric_2009,benson_large_2011,
beirao_da_veiga_isogeometric_2012,uhm_tspline_2009,Echter2013170,Benson2013133}. 
The smoothness of the NURBS basis functions allows for a straightforward
construction of plate/shell elements. Particularly for thin shells, rotation-free
formulations can be easily constructed \cite{kiendl_isogeometric_2009,kiendl_bending_2010}.
Note that for multi-patch NURBS surfaces, rotation-free IGA elements require special treatment
at patch boundaries where the basis functions are found to be $C^0$ continuous. 
Furthermore, isogeometric plate/shell elements exhibit much less
pronounced shear-locking 
compared to standard FE plate/shell elements. Elements with smooth boundaries such as 
circular and cylindrical elements were successfully constructed using the IGA concept
\cite{jia_circular_2009,lu_cylindrical_2011}.

The smoothness of NURBS basis functions is attractive for analysis of 
fluids \cite{gomez_isogeometric_2010,nielsen_discretizations_2011,Bazilevs:2010:LES:1749635.1750210} and for
fluid-structure interaction problems 
\cite{bazilevs_isogeometric_2008,bazilevs_patient-specific_2009}. 
In addition, due to the ease of constructing high order continuous basis functions, IGA has been 
used with great success in solving PDEs that incorporate fourth order (or
higher)
derivatives of the field variable such as the Hill-Cahnard equation 
\cite{gomez_isogeometric_2008}, explicit gradient damage models \cite{verhoosel_isogeometric_2011-1} and gradient
elasticity \cite{fischer_isogeometric_2010}. 
The high order NURBS basis has also found potential applications in the Kohn-Sham equation for electronic 
structure modeling of semiconducting materials \cite{Masud2012112}.

NURBS provide advantageous properties for structural vibration problems 
\cite{cottrell_isogeometric_2006,Hughes20084104,NME:NME4282,Wang2013}
where $k$-refinement (unique to IGA)
has been shown to provide more robust and accurate frequency spectra than 
typical higher-order FE $p$-methods. Particularly, the optical branches of frequency spectra,
which have been identified as contributors to Gibbs phenomena in wave propagation problems 
(and the cause of rapid degradation of higher modes in the $p$-version of FEM), 
are eliminated. However when lumped mass matrices were used, the accuracy is limited to second order
for any basis order. High order isogeometric lumped mass matrices are not yet available.
The mathematical properties of IGA were studied in detail by Evans et al.\cite{evans_n-widths_2009}.

The isogeometric concept has also spread to the field of meshfree methods such as \cite{Shaw20081541,NME:NME4360} in which spline-based meshfree methods were presented.
Industrial applications of IGA have been presented in  \cite{GroÃŸmann2012519,Scott2013197} along with applications in experimental mechanics \cite{elguedj_isogeometric_2011}
where NURBS-based DIC (Digital Image Correlation) was shown to outperform standard FE DIC.

\subsection{Shortcomings of NURBS and alternative geometry discretisations}

NURBS are ubiquitous in CAD but are known to exhibit major
shortcomings from a computational geometry standpoint. Perhaps the greatest difficulty encountered is the inability of NURBS
to produce watertight geometries, often
complicating mesh generation.  From an analysis perspective, the tensor product structure of NURBS proves to be inefficient, caused by the global nature of refinement operations.  In turn, this leads to inefficient error estimation and adaptivity algorithms.   One solution which has gathered momentum from both the computational geometry and analysis communities is the use of T-splines \cite{Sederberg:2003:TT:882262.882295} which overcome the limitations of NURBS while retaining the familiar
structure of NURBS algorithms. T-splines correct the
deficiencies of NURBS by creating a single patch, watertight geometry which can be locally refined and coarsened.  Buffa et al. \cite{Buffa_TSplines_2010} note that linear independence of the T-spline basis
functions is not guaranteed on generic T-meshes leading to the definition of analysis-suitable T-splines
\cite{Scott_TSplines_2012}, a mildly restricted subset of T-splines
which meet the demands of both design and analysis.
Utilisation of T-splines in an IGA framework has been illustrated in
\cite{bazilevs_isogeometric_2010,doerfel_adaptive_2010}, and by adopting a B\'{e}zier extraction process, Scott et al.\cite{scott_isogeometric_2011} showed that T-splines can be incorporated efficiently into existing FE codes.

Alternatives to T-splines include polycube splines \cite{Wang_polycube}, PHT-splines 
\cite{Deng_PSH} and LR-splines \cite{Dokken_ECCM}. 
PHT-splines (polynomial spline over hierarchical T-meshes) have been extended to 
rational splines and applied in \cite{nguyen-thanh_isogeometric_2011,nguyen-thanh_rotation_2011} to problems 
in elasticity for continua and thin structures.  Adaptive refinement with PHT-splines is 
particularly simple. Although T-splines allow for local adaptive
refinement, the complexity of knot insertion under adaptive refinement
is complex, particularly in 3D. However, we note that research is currently being pursued on hierarchical T-spline refinement algorithms that address this issue.

Another direction of IGA research includes hierarchical B-splines 
\cite{vuong_hierarchical_2011,Schillinger2012,Bornemann2013584} and unstructured 
Powell-Sabin splines \cite{Speleers2012132}. 
The hierarchical B-splines finite cell method \cite{Schillinger2012} furnishes a
seamless CAD-FEA integration for very complex geometries.
We refer also to \cite{Kleiss2012168} for IGA combined with finite element based local refinement capabilities.
Different subdivision surface techniques (Catmull-Clark, Loop) have
also been utilized
for solid and shell modeling \cite{CGF:CGF1766,Polthier-subdivision}.

In computer aided geometric design, patching multiple NURBS
parameterizations to form
complex topologies is far from trivial if certain continuity
requirements are to be maintained. Trimming techniques provide a promising alternative for
representing complex NURBS domains. 
In \cite{kim_isogeometric_2010}, a trimmed surface based analysis framework has been proposed 
where NURBS-enhanced FEM \cite{sevilla_3d_2011,sevilla_nurbsenhanced_2008} was applied to define a suitable 
integration domain within parameter space. In a recent contribution \cite{Schmidt201293},  the authors
presented an alternative method to handle trimmed NURBS geometries.

\subsection{Discontinuities and fracture}

IGA has been applied to cohesive fracture \cite{verhoosel_isogeometric_2011}, outlining a framework for
modeling debonding along material interfaces using NURBS and propagating cohesive
cracks using T-splines. The method relies upon the ability to specify the continuity of NURBS and T-splines
through a process known as knot insertion.
As a variation of the eXtended Finite Element Method (XFEM)
\cite{mos_finite_1999}, IGA was applied to Linear Elastic Fracture Mechanics (LEFM) using the partition
of unity method (PUM) to capture two dimensional strong discontinuities and
crack tip singularities efficiently \cite{de_luycker_xfem_2011,ghorashi_extended_2012}. The method is usually
referred to as XIGA (eXtended IGA).
In \cite{Tambat20121} an explicit isogeometric enrichment technique was proposed for modeling
material interfaces and cracks exactly. Note that this method
is contrary to PUM-based enrichment methods which define cracks implicitly.

A phase field model for dynamic fracture was presented in \cite{Borden201277} using adaptive T-spline refinement to provide an effective method for simulating fracture in three dimensions.
In \cite{Nguyen2013} high order B-splines were adopted to efficiently 
model delamination of composite specimens and in \cite{nguyen_cohesive_2013}, an isogeometric framework
for two and three dimensional delamination analysis of composite laminates was presented where
the authors showed that using IGA can significantly reduce the usually time consuming pre-processing
step in generating FE meshes (solid elements and cohesive interface elements) for delamination computations.
A continuum description of fracture using 
explicit gradient damage models was also studied using NURBS \cite{verhoosel_isogeometric_2011-1}.

\subsection{Alternatives to IGA}

Other techniques which integrate CAD and analysis include the use of subdivision surfaces to model shells \cite{Cirak_2000},
NURBS-enhanced finite elements
\cite{sevilla_3d_2011,sevilla_nurbsenhanced_2008} and NURBS for BEM shape
optimisation \cite{Cervera2005}. Immersed boundary methods 
\cite{Sanches20111432} (and references therein),
the finite cell method \cite{Schillinger2012,Rank2012104} and the structured XFEM
\cite{Moumnassi2011774,legrain-xfem-nurbs} are yet other alternatives
which aim to combine analysis and design technologies. In general, it
is found that for CAD and analysis technologies to work seamlessly
together, the underlying discretisation must either be directly
compatibile or easily converted between the two. 

IGA has offered
significant advances towards the goal of a unified design and analysis
framework, but much research is still needed before this goal is
realised. There are several indications of the future promise of IGA for industrial design but ultimately, the litmus test of success for IGA will be whether the
approach is widely adopted by industry.

\subsection{Computational aspects}\label{sec:comp_aspects}

Some major computational aspects of IGA which have been studied so far include (i) locking issues, 
(ii) sensitivity to mesh distortion, (iii) impact of high continuity of NURBS on direct solvers,
(iv) collocation methods, (v) competing demands of analysis and 
computational geometry discretisations, (vi) construction of
trivariate solids from given bivariate surface representations and (vii) optimal quadrature rules.

\begin{enumerate}[(i)]
\item Although the smoothness of NURBS basis functions reduces to some extent the locking phenomena for
constrained problems such as incompressible media, thin-walled structures, NURBS-based FEs are not 
locking free \cite{auricchio_fully_2007,elguedj_and_2008,Echter2010374,iga-Fbar,NME:NME4328}. 
Existing locking removal
techniques in standard FEMs were successfully adapted to IGA such as the Discrete Strain Gap method \cite{Echter2010374}, the F/B-bar method \cite{elguedj_and_2008,iga-Fbar}, the enhanced assumed strain method \cite{NME:NME4328} 
and mixed formulations  \cite{NME:NME3048}.

\item The effect of mesh distortion on the performance of IGA for solid mechanics was discussed in 
\cite{Lipton2010357} in which it was found that higher-order NURBS functions are able to 
somewhat alleviate the impact of the distortions.

\item The high order continuity offered by NURBS has a negative impact on the performance of direct solvers
as pointed out in \cite{Collier2012353}. The authors
found that for a fixed number of unknowns and basis degree, a higher degree of continuity
drastically increases the CPU time and RAM needed to solve the problem when using a direct solver.

\item In an attempt to compete with low order FEs with one-point quadrature that are extensively used in
industrial applications, isogeometric collocation methods were developed 
\cite{auricchio-collocation1,Auricchio20122}.

\item Due to the fact that meshes in an isogeometric framework  are defined by the parametrisation of 
the object of interest, the quality of the geometry parametrisation
plays an important role in ensuring mesh quality. This issue has, however, been addressed by only a few researchers
\cite{takacs_existence_2011,xu_parameterization_2011,cohen_analysis-aware_2010,schmidt_iga_cad,
   Zhou:2005:NGM:1060244.1060253}.
In particular, in \cite{cohen_analysis-aware_2010}, the authors proposed the concept of ``analysis-aware
geometry modeling".

\item In CAD, solids are defined as boundary surfaces in which the
  interior is not explicitly modeled. In FEA, a solid representation
  is necessary and therefore, the transition from CAD solids to FEA
  solids demands a step in which the CAD representations are converted
  to solid FEA representations. Initial developments have been reported 
in \cite{zhang_patient-specific_2007,Zhang2012185,zhang-solid-tsplines,
   Martin:2009:VPT:1542572.1543076,aigner_swept,escobar_new_2011}.
Note that in this regard, the isogeometric boundary element method (IGABEM) can be considered
a truly isogeometric method 
\cite{simpson_two-dimensional_2012,Scott2013197,Peake201393,Simpson20132}
since BEM analysis requires only the definition of a boundary
discretisation, completely defined by CAD.

\item Gaussian quadrature is not optimal for IGA.
Research is currently focussed on optimal integration techniques such
as that in \cite{hughes_efficient_2010,Auricchio201215} in which 
(nearly) optimal quadrature rules have been presented.
\end{enumerate}

\subsection{Available implementations}

Some implementation aspects of IGA were reported in 
\cite{cottrel_book_2009} and more recently, an open source IGA Matlab\textsuperscript{\textregistered{}}  
code was
described in \cite{vuong_isogat:_2010} with a restriction to 2D scalar PDEs.
An excellent open source IGA code written in Matlab\textsuperscript{\textregistered{}} is given in 
\cite{falcao_geopde_2011}.
Incorporating IGA within an object-oriented C++ FE code was discussed in \cite{rypl_iga_2012}. Implementation
details for enriched formulations within an IGA framework are reported in
\cite{benson_gfem_2010} using commercial FE software. An IGA BEM code written in 
Matlab\textsuperscript{\textregistered{}} was presented in \cite{Simpson20132}.
Isogeometric analysis was also incorporated into  FEAP \cite{feap,NME:NME3048}.
A high performance IGA code was given in \cite{petiga} which is based on PETSc, 
the Portable, Extensible Toolkit for Scientific Computation. The Python programming language has
also been adopted to implement IGA \eg \cite{python-iga,pigasus}.

\subsection{Contributions and outline}

In this paper, we present in detail an isogeometric (Bubnov-)Galerkin finite
element method applied to two- and three-dimensional elasto-static solid/structural mechanics
problems and traction-free crack problems. 
The discussion is confined to NURBS for the sake of simplicity. 
Although no new fundamental findings are presented the contribution of the paper are  

\begin{itemize}
\item an overview of IGA, applications and recent developments are presented;
\item implementation details for 1D/2D/3D solid mechanics and rotation-free
      plate elements are provided;
\item implementation for 2D and 3D XIGA for traction-free cracks;      
\item visualization techniques for (X)IGA are discussed;
\item techniques to impose Dirichlet boundary conditions including the penalty method, the Lagrange multipliers
      method, the least square projection method are implemented;
\item different techniques to model discontinuities in the context of IGA are discussed;      
\end{itemize}
where some implementations are not available in the literature.

The paper is structured as follows: Section \ref{sec-nurbs} outlines B-spline and NURBS technology used to construct curves,
surfaces and solids; the use of NURBS for discretisation within a finite element framework is treated in Section
\ref{nurbs-fem}; implementation of an isogeometric finite element method for two-dimensional elasticity is described in Section~\ref{sec:elast-two-dimens} and 
an extended isogeometric formulation based on the concept of the PUM is the subject of Section \ref{sec:xiga}.  
A detailed description of our IGA Matlab code is given in Section \ref{sec:migfem} and
Application of IGA to 
structural mechanics problems is presented in Section \ref{sec:structural-mechanics}.
Numerical examples including 2D/3D fracture mechanics and large deformation shell problems
are given in Section \ref{examples}.

\subsection{Notation}
\label{sec:notation}

We use lowercase indices to indicate a local index and uppercase indices to indicate a global index. We denote $d_p$ and $d_s$ as the number of parametric directions and spatial directions respectively. \textbf{Boldfont} is used to indicate matrices and vectors where the number of components is implied.

\section{A brief introduction to B-splines/NURBS}\label{sec-nurbs}

As a precursor to Section~\ref{nurbs-fem}, an understanding of the discretisation technology which underpins IGA is essential, since the beneficial properties which are found through analysis emanate directly from the underlying CAD basis functions. We give a brief outline of parametric functions, then state the terminology and basis function definitions which allow curves, surfaces and solids to be represented by B-splines and NURBS. The algorithms which define B-splines and NURBS in the present Matlab\textsuperscript{\textregistered{}} code are based on those given in \cite{piegl_book} in which explicit implementations are given for all algorithms.

\subsection{Parametric representation of geometry}

What is fundamental to all the discretisation technology used in the present paper (and the majority of technology in the CAGD community), is the representation of geometry through \textit{parametric functions}. These define a mapping from a given parameter to the desired geometry. We can imagine that as we move through the parameter domain, the parametric function `sweeps' out the desired shape. For example, if we consider the case of a circle of radius 1, the implicit form of this equation is given by
\begin{equation}
x^2 + y^2 = 1
\label{eq:circle}
\end{equation}
or alternatively, if we define a mapping $f: [0,2\pi] \to \mathbb{R}^2$, the parametric form of the same circle is given by
\begin{align}
f(t) =
\left(
\cos t,
\sin t \right).
\label{eq:circle-parametric}
\end{align}
This is found to be much more conducive for graphical implementation. To see this, consider the task of plotting the circle through both Eq.~\eqref{eq:circle} and Eq.~\eqref{eq:circle-parametric}. In the case of the latter,  by simply determining $t$ at a discrete set of points in the interval $[0,2\pi]$, the desired result is obtained. In addition, many algorithms which perform geometrical transformations become much simpler when parametric functions are used. Both B-splines and NURBS are based on parametric functions.

\subsection{B-splines}

The demands of Computer Aided Geometrical Design place certain
requirements on the types of discretisation that can be used to
represent geometrical objects where, for example,  in the case of car design, surfaces of $G^2$ continuity\footnote{\ref{sec:disc-cont} outlines the differences between $C^k$ and $G^k$ continuity} are required to avoid unnatural reflections. Additional requirements include local control of geometrical
features, the ability to apply refinement algorithms and numerical
stability of high order curves. It is found that
representations of objects using B\'{e}zier curves or Lagrange
polynomials do not meet many of these requirements and alternatives
are sought. The CAGD community have settled on the use of
B-spline-based technology since it is
found to provide many of the desired properties for interactive
geometrical design, realised through the properties of the underlying
B-spline basis functions. 

B-splines can be considered as a mapping  from \textit{parametric space} $\hat{\Omega} \subset \mathbb{R}^{d_p}$
to 
\textit{physical space} $\Omega \subset \mathbb{R}^{d_s}$ \footnote{Further discussion of these spaces in given in Section~\ref{sec:relevant-spaces}.}.  In this sense, a B-spline can be considered to `sweep' out a curve, surface or volume as we move through the range of parameter values. In the present work we define coordinates in parameter space as $\boldsymbol{\xi} = (\xi,\eta,\zeta) = (\xi^1, \xi^2, \xi^3)$ and coordinates in physical space as $\mathbf{x} = (x,y,z) = (x^1, x^2, x^3)$. These are simplified accordingly in the case of one- and two-dimensional problems.  What remains is to determine the particular form of the mapping $\mathbf{x}: \hat{\Omega} \to \Omega$. 

To construct a B-spline, a \textit{knot vector} must be specified and is defined as an ordered set of increasing parameter values $\Xi=\{\xi_1,\xi_2,\ldots,\xi_{n+p+1}\}$, $\xi_i \leq \xi_{i+1}$ where $\xi_i$ is the \textit{i}th knot, $n$ is the number of basis functions and $p$ is 
the polynomial order. The knot vector divides the parametric space into intervals usually referred to
as \textit{knot spans} with the number of 
coincident knots for a particular knot value referred to as a knot with a certain \textit{multiplicity} $k$. That is, a knot has a multiplicity $k$ if it is repeated $k$ times in the knot vector. Most commonly, 
\textit{open} knot vectors are used where the first and last knots have a
multiplicity $k=p+1$. \ref{sec:knot-vect-conv} outlines alternative knot vector notations which remove redundant information, but it can be assumed that the above knot vector notation is used throughout the present work.

An important property of B-splines formed from open knot vectors is that they are interpolatory at their start and end points, greatly facilitating the imposition of boundary conditions
for analysis.  In geometrical terms, this corresponds to a B-spline that coincides with its start and end control points. A \textit{uniform knot vector} is associated to evenly
distributed knots. Otherwise it is classified as a \textit{non-uniform knot
vector}. 

Knot vectors are not commonly used by CAD designers, and in most CAD software the ability to modify knot vectors is not provided. It is much more common to tailor the geometry through modification of the polynomial order, control points and weightings.

\subsubsection{B-spline basis functions}

Given a knot vector $\Xi$, the associated set of B-spline basis functions $\{N_{i,p}\}_{i=1}^n$ are
defined recursively by the Cox-de-Boor formula, starting with the zeroth order basis
function ($p=0$)
\begin{equation}
  N_{i,0}(\xi) = \begin{cases}
    1 & \textrm{if $ \xi_i \le \xi < \xi_{i+1}$},\\
    0 & \textrm{otherwise},
  \end{cases}
  \label{eq:basis-p0}
\end{equation}
and for a polynomial order $p \ge 1$
\begin{equation}
  N_{i,p}(\xi) = \frac{\xi-\xi_i}{\xi_{i+p}-\xi_i} N_{i,p-1}(\xi)
               + \frac{\xi_{i+p+1}-\xi}{\xi_{i+p+1}-\xi_{i+1}}
	       N_{i+1,p-1}(\xi).
  \label{eq:basis-p}
\end{equation}
in which fractions of the form $0/0$ are
defined as zero.

Some salient properties of B-spline basis functions include: (1) they constitute 
a partition of unity, (2) each basis function is non-negative over the entire parametric domain,
(3) they are linearly independent, (4) the support of a B-spline basis function $N_{i,p}$  is given by $p+1$ knot spans denoted by $[\xi_i,\xi_{i+p+1}]$, (5) they exhibit $C^{p-m_i}$ continuity across knot $\xi_i$ where $m_i$ is the multiplicity of knot $\xi_i$ and (6) 
B-spline basis functions in general are not interpolatory, equivalent to saying that the Kronecker-delta property is not guaranteed at control points. This last point requires careful treatment when imposing non-homogeneous Dirichlet boundary conditions, and is discussed later
in Section \ref{sec:BCs}.

\begin{figure}[htbp]
  \centering 
  \includegraphics[width=0.6\textwidth]{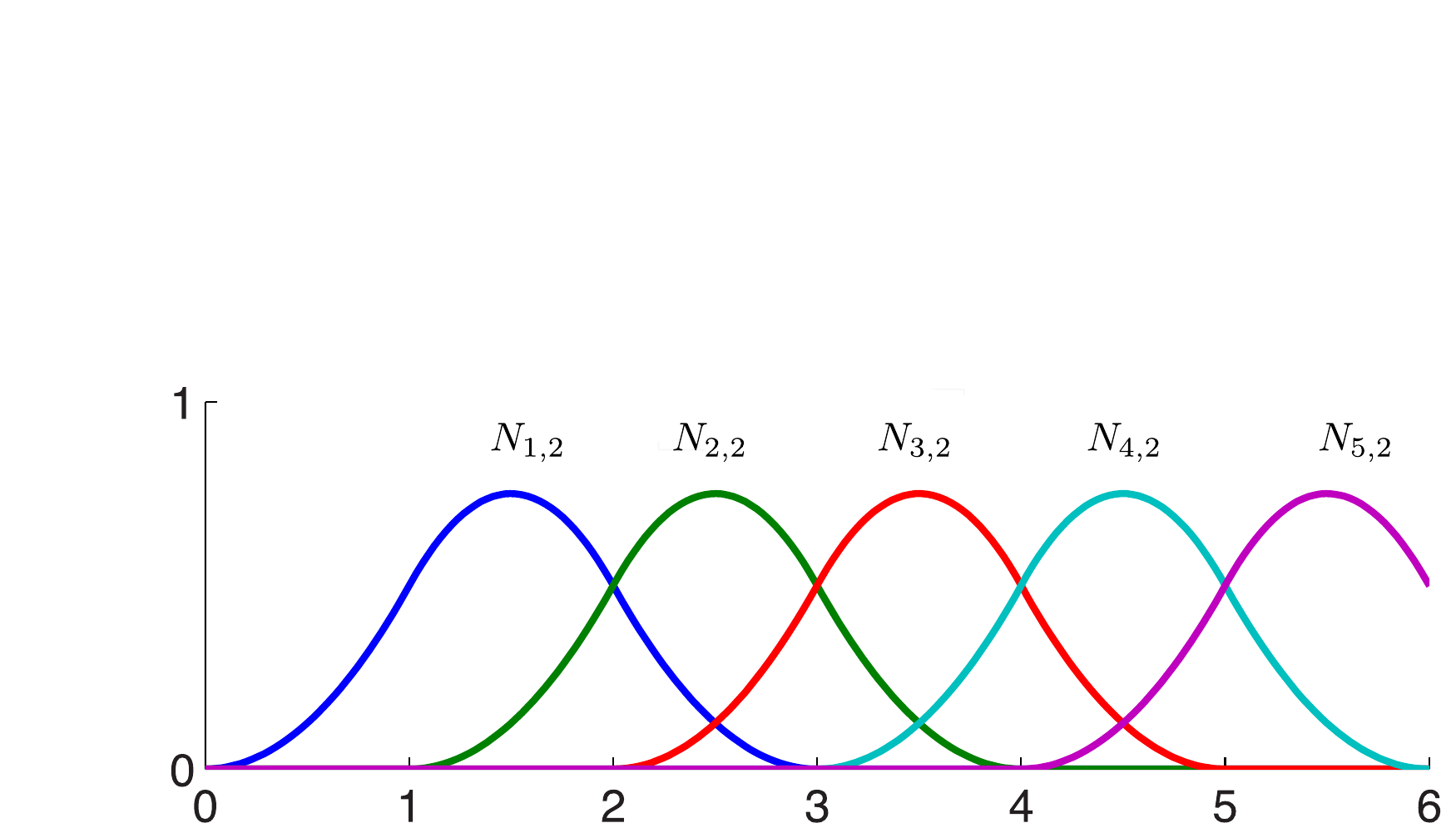}
  \caption{Quadratic B-spline basis functions defined for the uniform knot vector
  $\Xi=\{0,1,2,3,4,5,6\}$.} 
  \label{fig:bspline-quad-uniform} 
\end{figure}

\begin{figure}[h!]
  \centering 
  \includegraphics[width=0.7\textwidth]{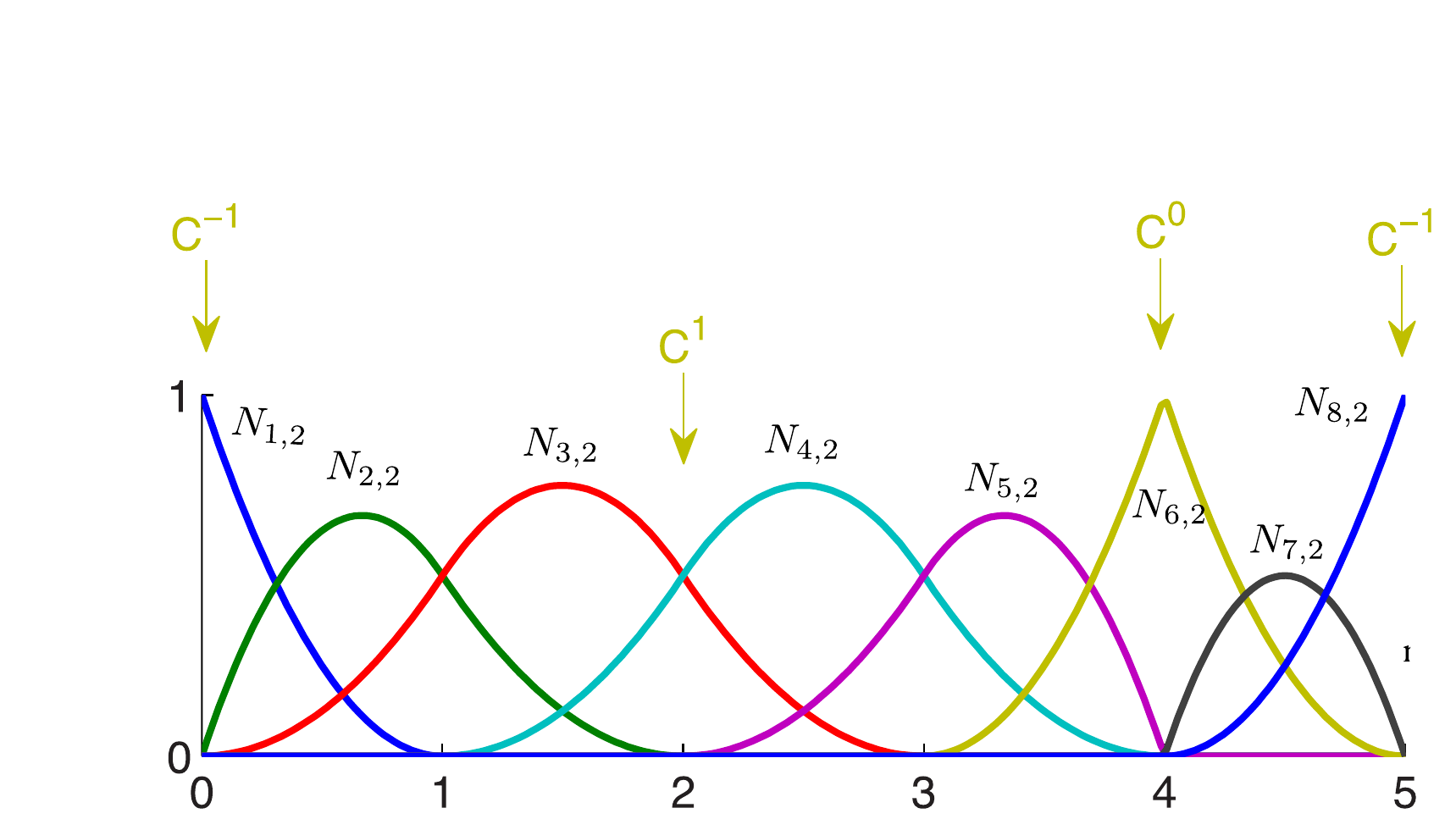}
  \caption{Quadratic B-spline basis functions defined for the open, non-uniform knot vector
  $\Xi=\{0,0,0,1,2,3,4,4,5,5,5\}$. Note the flexibility in the construction of
  basis functions with varying degrees of regularity.} 
  \label{fig:bspline-quad-open} 
\end{figure}

To demonstrate the effect of alternative knot vectors, Fig. \ref{fig:bspline-quad-uniform} illustrates a set of quadratic ($p=2$) B-spline
basis functions for a uniform knot vector. Fig. \ref{fig:bspline-quad-open} illustrates a corresponding set of basis functions for an open, non-uniform knot vector. Of particular note is the interpolatory nature of the basis function at each
end of the interval created through an open knot vector, and the reduced continuity at $\xi = 4$ due to the presence of the location of a repeated knot where $C^0$ continuity is attained. Elsewhere, the functions are $C^1$ continuous ($C^{p-1}$).
The ability to control continuity 
by means of knot insertion is particularly useful for modeling discontinuities
such as cracks or material interfaces \cite{verhoosel_isogeometric_2011,nguyen_cohesive_2013}. 

During implementation, the derivatives of B-splines are required in  both CAD and analysis to compute quantities such as tangent and normal vectors and field variable derivatives. The first derivative of a B-spline basis function is computed recursively from lower order basis functions as
\begin{equation}
  \frac{d}{d\xi}N_{i,p}(\xi) = \frac{p}{\xi_{i+p}   - \xi_i} N_{i,p-1}(\xi) -
                               \frac{p}{\xi_{i+p+1} - \xi_{i+1}}
			       N_{i+1,p-1}(\xi).
  \label{eq:basis-derivatives}
\end{equation}
Derivatives of higher order can be found in \cite{piegl_book}.

\subsection{Refinement algorithms}
\label{sec:refinement}

A feature which is essential for both computational geometry and analysis is the ability to successively refine discretisations to allow intricate geometries to be modelled or to capture rapid variations in the field solution. A significant advantage of adopting B-spline-based functions as a discretisation tool is the  ability to apply a variety of refinement algorithms in a simple manner. We restrict ourselves to a brief discussion of the available algorithms, with a more detailed treatment given in \cite{cottrel_book_2009}.

B-spline refinement algorithms in the context of isogeometric analysis can be assigned to one of three types:
\begin{enumerate}
\item \textbf{Knot insertion} is the process whereby additional knots are inserted into the knot vectors thereby creating additional knot intervals or elements, in the context of analysis. This is directly analogous to h-refinement seen in the conventional FEM.
\item \textbf{Degree elevation} is the process of raising the order of the underlying basis, directly analogous to p-refinement in the  conventional FEM.
\item \textbf{k-refinement} is unique to IGA and consists of a combined process of degree elevation and knot insertion. These processes are not commutative and therefore the order in which these refinements are applied will change the final basis. k-refinement first applies degree elevation proceeded by knot insertion, offering a reduction in degrees of freedom over its counterpart \cite{hughes_isogeometric_2005}.
\end{enumerate}

\subsection{B-spline curves}

\begin{figure}[htbp]
  \centering 
  \includegraphics[width=0.4\textwidth]{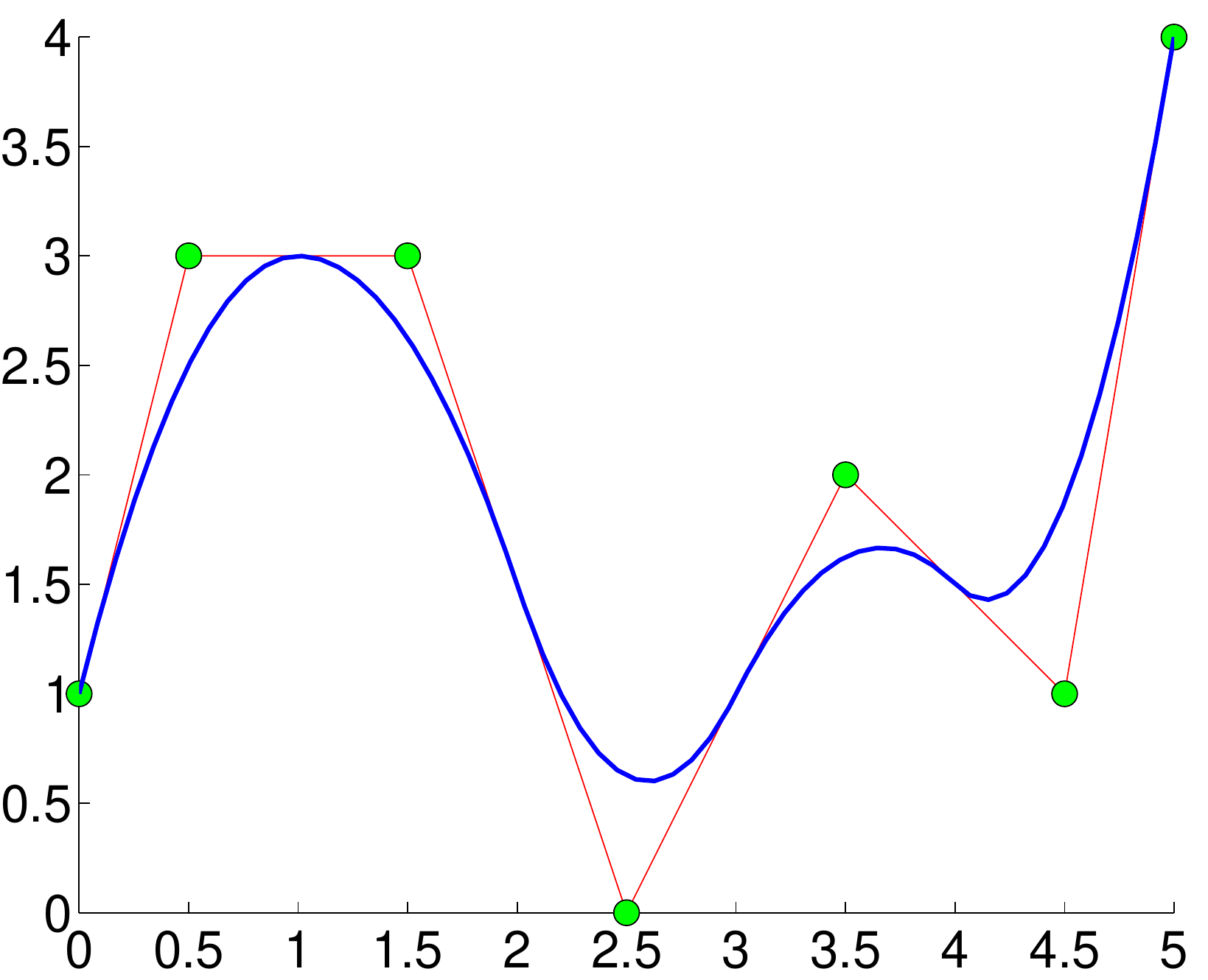}
  \caption{A quadratic ($p=2$) B-spline curve with a uniform open knot vector 
  $\Xi=\{0,0,0,1,2,3,4,5,5,5\}$. Control points are denoted by filled circles
     and the control polygon is denoted by the red line.
  Note that the curve is at least $C^1$ continuous everywhere except at the
  extremities where, due to the presence of repeated knots, $C^{-1}
$ continuity is obtained.} 
  \label{fig:splinecurve1} 
\end{figure}

Given a set of B-spline basis functions $\{N_{A,p}\}_{A=1}^n$ and a set of \textit{control points} $\{\mathbf{P}_A\}_{A=1}^n$ where $\vm{P}_A \in \mathbb{R}^{d_s}$, we can construct a piecewise-polynomial B-spline curve as
\begin{equation}
	\vm{C}(\xi) = \sum_{A=1}^{n} \vm{P}_A N_{A,p}(\xi).
  \label{eq:Bspline-curve}
\end{equation}
An example quadratic B-spline curve is shown in Fig. \ref{fig:splinecurve1} constructed from a uniform open knot vector. Its associated \textit{control polygon}, constructed by piecewise linear interpolation of the control points is also shown.  A B-spline curve inherits all of the continuity properties
of its underlying basis and, as mentioned previously, the use of open knots ensures that the first and
last points are interpolatory. Note that the tangent to the curve at
the first control point is defined by the first two control
points. This property is exploited in rotation-free element formulations 
for thin-walled structures (section \ref{sec:structural-mechanics}).

\subsection{B-spline surfaces and volumes}

\begin{figure}[htbp]
  \centering 
  \psfrag{n}{$ N_{3,3}(\xi)$}
  \psfrag{m}{$ M_{3,3}(\eta)$}
  \psfrag{r}{$N_{3,3}^{3,3}(\xi,\eta)$}
  \includegraphics[width=0.75\textwidth]{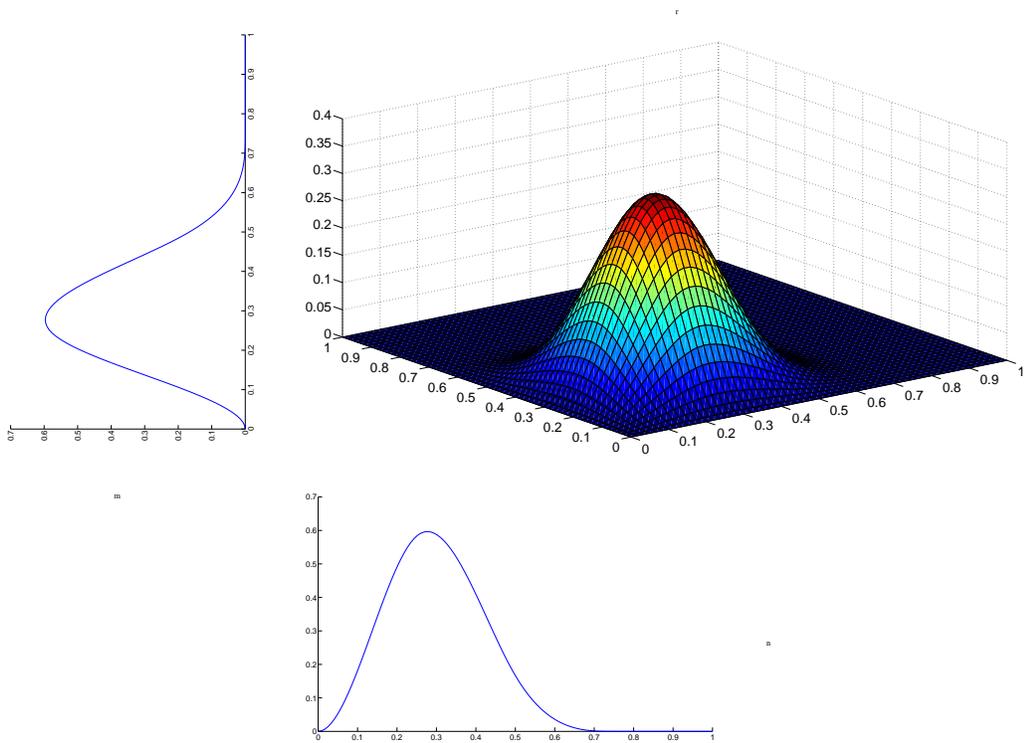}
  \caption{A bivariate cubic B-spline basis function for knot vectors 
  $\Xi^1=\Xi^2=\{0,0,0,0,0.25,0.5,0.75,1,1,1,1\}$.}
  \label{fig:bivariate-bspline} 
\end{figure}

Given two knot vectors for each parametric direction $\Xi^1=\{\xi_1,\xi_2,\ldots,\xi_{n+p+1}\}$
and $\Xi^2=\{\eta_1,\eta_2,\ldots,\eta_{m+q+1}\}$, and a control net $\vm{P}_{i,j} \in \mathds{R}^{d_s}$,
a tensor-product B-spline surface is defined as
\begin{equation}
	\vm{S}(\xi,\eta) =
	\sum_{i=1}^{n}\sum_{j=1}^{m}N_{i,p}(\xi)M_{j,q}(\eta) \vm{P}_{i,j},
  \label{eq:Bspline-surface}
\end{equation}
where $N_{i,p}(\xi)$ and $M_{j,q}(\eta)$ are the univariate B-spline
basis functions of order $p$ and $q$ corresponding to knot vectors $\Xi^1$ and
$\Xi^2$, respectively. Defining a global index as
\begin{equation}
  \label{eq:bspline_surface_mapping}
  A = n( j - 1 ) + i, 
\end{equation}
Eq.~\eqref{eq:Bspline-surface} can be rewritten in a more compact form as
\begin{equation}
	\vm{S}(\boldsymbol{\xi}) = \sum_{A=1}^{n \times m} \vm{P}_A  N_{A}^{p,q}(\boldsymbol{\xi} ),
  \label{eq:Bspline-surface1}
\end{equation}
in which $N_{A}^{p,q}$ is a bivariate B-spline basis function defined as $N_{A}^{p,q}(\boldsymbol{\xi}) = N_{i,p}(\xi)M_{j,q}(\eta)$. 
Figs. \ref{fig:bivariate-bspline} illustrates an example bicubic B-spline basis function.

The extension to B-spline volumes is straightforward, where a trivariate basis is formed through a tensor product of B-spline basis functions as
\begin{equation}
  \label{eq:bspline_volume}
  	\vm{V}(\xi,\eta,\zeta) =
	\sum_{i=1}^{n}\sum_{j=1}^{m}\sum_{k=1}^l N_{i,p}(\xi)M_{j,q}(\eta) L_{k,r}(\zeta)\vm{P}_{i,j,k},
\end{equation}
or, by defining a global index $A$ through
\begin{equation}
  \label{eq:bspline_volume_mapping}
  A = (n \times m) ( k - 1)  + n( j - 1 ) + i, 
\end{equation}
a simplified form of Eq.~\eqref{eq:bspline_volume} can be written as
\begin{equation}
  \label{eq:bspline_vol_simple}
  	\vm{V}(\boldsymbol{\xi}) = \sum_{A=1}^{n \times m \times l} \vm{P}_A  N_{A}^{p,q,r}(\boldsymbol{\xi} ).
\end{equation}

To compute derivatives of B-spline surfaces and volumes, the chain rule is applied to \eqref{eq:basis-derivatives} with a detailed treatment of this topic provided in \cite{piegl_book}.

\subsection{NURBS}\label{nurbs}

B-splines are convenient for free-form modelling, but they lack the ability to exactly 
represent some simple shapes such as circles and ellipsoids. 
This is why today, the \textit{de facto} standard technology in CAD is
a generalisation of B-splines referred to as NURBS (Non-Uniform Rational B-Splines). NURBS are formed through rational functions of B-splines, forming a superset of B-splines. They inherit all the  favourable properties of B-splines and are favoured over their counterpart due to their ability to form exact representations of conic
sections such as spheres, ellipsoids, paraboloids and hyperboloids. In addition, there exist
efficient algorithms for their evaluation and refinement.

\subsubsection{NURBS basis functions}

NURBS basis functions are defined as
\begin{equation}
	R_{i,p}(\xi) = \frac{N_{i,p}(\xi)w_i}{W(\xi)} =
	\frac{N_{i,p}(\xi)w_i}{\sum_{\hat{i}=1}^{n}N_{\hat{i},p}(\xi)w_{\hat{i}}}
  \label{eq:rational-basis}
\end{equation}
where $\{N_{i,p}\}_{i=1}^n$ is the set of B-spline basis functions of
order $p$ and $\{w_i\}_{i=1}^n, w_i >0$ is the set of NURBS weights. 
Selecting appropriate weights permits the description of many
different types of curves including polynomials and circular arcs.
For the special case in which all weights are equal, the
NURBS basis reduces to the B-spline basis. NURBS weights for certain simple geometries are given in \cite{piegl_book}, but in general, weights are user-defined through CAD packages such as Rhino \cite{rhino}.

The first derivative of a NURBS basis function $R_{i,p}$ is computed using the
quotient rule as
\begin{equation}
	\frac{d}{d\xi} R_{i,p}(\xi)=w_i
	\frac{N_{i,p}'(\xi)W(\xi) - N_{i,p}(\xi)W'(\xi)}{W(\xi)^2},
  \label{eq:nurbs-derivatives-1}
\end{equation}
where $N_{i,p}'(\xi) \equiv \frac{d}{d\xi}N_{i,p}(\xi)$ and
\begin{equation}
	W'(\xi)=\sum_{\hat{i}=1}^nN_{\hat{i},p}'(\xi)w_{\hat{i}}.
	\label{dddd}
\end{equation}

\subsubsection{NURBS curves, surfaces and volumes}

In a similar fashion to B-spline curves, the NURBS curve associated with a set of control points and weights  $\{\vm{P}_A, w_A\}_{A=1}^n$
and basis functions $\{R_{A,p}\}_{A=1}^n$ is defined as
\begin{equation}
	\vm{C}(\xi) = \sum_{A=1}^{n} \vm{P}_A R_{A,p}(\xi).
  \label{eq:NURBS-curve}
\end{equation}
NURBS surfaces are constructed from a linear combination of bivariate NURBS basis functions, control points $\vm{P}_{i,j}\in \mathds{R}^{d_s} $ and weights $w_{i,j} > 0$ as 
\begin{equation}
	\vm{S}(\xi,\eta) =
	\sum_{i=1}^{n}\sum_{j=1}^{m} \vm{P}_{i,j}  R_{i,j}^{p,q}(\xi,\eta),
  \label{eq:NURBS-surface1}
\end{equation}
where the bivariate NURBS basis functions are defined as
\begin{equation}
	R_{i,j}^{p,q}(\xi,\eta) = \frac{N_{i}(\xi) M_j(\eta) w_{i,j}}{
  \sum_{\hat{i}=1}^{n} \sum_{\hat{j}=1}^{m} N_{\hat{i}}(\xi) M_{\hat{j}}(\eta)
  w_{\hat{i},\hat{j}}}.
  \label{eq:rational-basis2}
\end{equation}
Alternatively, using the mapping defined by Eq.~\eqref{eq:bspline_surface_mapping}, Eq.~\eqref{eq:NURBS-surface1} can be written more succinctly as
\begin{equation}
  \label{eq:nurbs_surface_simplified}
  \vm{S}(\boldsymbol{\xi}) = \sum_{A=1}^{n \times m} \vm{P}_A R_A^{p,q}(\boldsymbol{\xi}).
\end{equation}
NURBS volumes are constructed from control points $\vm{P}_{i,j,k}\in \mathbb{R}^{d_s}$, weights $w_{i,j,k} > 0$ as
\begin{equation}
	\vm{V}(\xi,\eta,\zeta) = \sum_{i=1}^{n}\sum_{j=1}^{m}\sum_{k=1}^l
	\vm{P}_{i,j,k}  R_{i,j,k}^{p,q,r}(\xi,\eta,\zeta)
  \label{eq:NURBS-solid1}
\end{equation}
where the trivariate NURBS basis functions $R_{i,j,k}^{p,q,r}$ are given by
\begin{equation}
  R_{i,j,k}^{p,q,r}(\xi,\eta,\zeta) = \frac{N_{i}(\xi) M_j(\eta) P_k(\zeta) w_{i,j,k}}{
  \sum_{\hat{i}=1}^{n} \sum_{\hat{j}=1}^{m}  \sum_{\hat{k}=1}^{l}N_{\hat{i}}(\xi) M_{\hat{j}}(\eta)
  P_{\hat{k}}(\zeta) w_{\hat{i},\hat{j},\hat{k}}}.
  \label{eq:rational-basis3}
\end{equation}
Using the mapping given by Eq.~\eqref{eq:bspline_volume_mapping}, Eq.~\eqref{eq:NURBS-solid1} can also be written as
\begin{equation}
  \label{eq:nurbs_solid_2}
  \vm{V}(\boldsymbol{\xi}) = \sum_{A=1}^{n \times m \times l} \vm{P}_A R_A^{p,q,r}(\boldsymbol{\xi})
\end{equation}

\begin{rmk}
As mentioned in Section~\ref{sec:comp_aspects}, CAD representations
are usually composed of surface models or boundary-representations. Trivariate
discretisations defined by \eqref{eq:nurbs_solid_2} are not normally
explicitly given, and therefore some preprocessing is required 
before domain based numerical methods such as the FEM can be
applied. Except for the case of simple cases such as extruded-surface
models and swept models, this task is far from trivial. This is currently an open research topic in IGA.
\end{rmk}

\section{NURBS as a basis for analysis: isogeometric finite element formulation}\label{nurbs-fem}

Our
attention now focusses on the use of B-splines and NURBS as a discretisation tool for analysis, outlining the core
concepts of isogeometric analysis. In this section the important
spaces and mappings are defined, followed by
the isogeometric FEM formulation in which we use NURBS as a basis
for analysis. The discussion is made more concrete through a
one-dimensional example.

\subsection{Relevant spaces}
\label{sec:relevant-spaces}

Familiarity must be gained with the spaces that are commonplace
in isogeometric analysis and the relationships that exist between each. Those that are considered presently in the context of B-splines and NURBS include: index, parameter, physical and parent space.

\subsubsection{Index space}
\label{sec:index-space}

\begin{figure}
  \centering
  \includegraphics{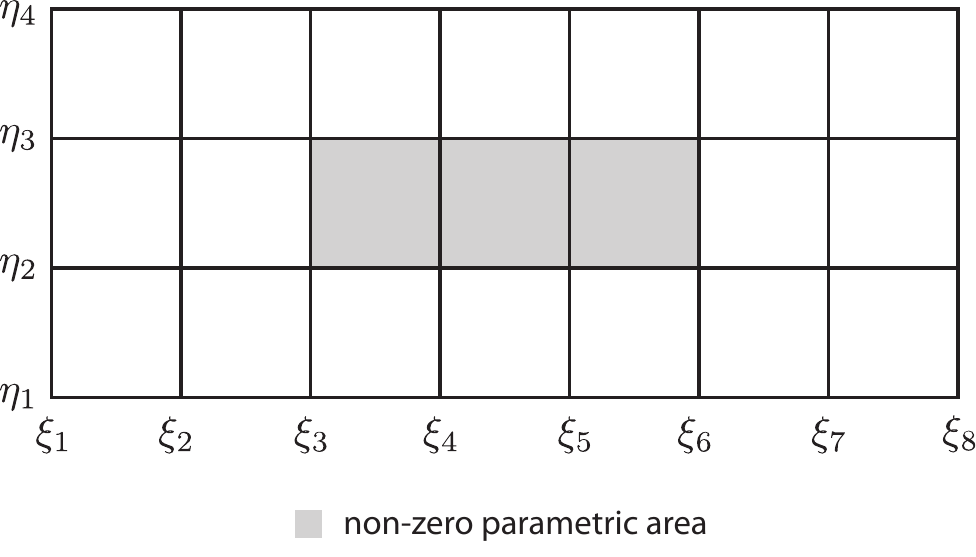}
  \caption{Creation of index space from knot vectors $\Xi^1 =
    \{0,0,0,1,2,3,3,3\}$, $\Xi^2 = \{0,0,1,1\}$ with non-zero parametric area highlighted.}
  \label{fig:index_space}
\end{figure}

Index space is formed through the specified knot vectors by giving each knot value a distinct coordinate, regardless of whether the knot is repeated or not. As an example, consider a NURBS patch defined through bivariate NURBS
basis functions with knot vectors $\Xi^1 = \{0,0,0,1,2,3,3,3\}$,
$\Xi^2 = \{0,0,1,1\}$ in each of the parametric directions $\xi,\eta$
respectively. This will form the index space as illustrated in
Fig.~\ref{fig:index_space} where the presence of repeated knots
leads to several regions of zero parametric area. Index space is often
used during implementation\footnote{T-splines being one notable
  example, with the local basis function mesh directly analogous to
  index space}, discarding elements of non-zero parametric
area, but in the present work we choose to only consider those
elements that have a non-zero parametric area, obviating the need for
index space.

\subsubsection{Parametric space}
\label{sec:parametric-space}

\begin{figure}
  \centering
  \includegraphics{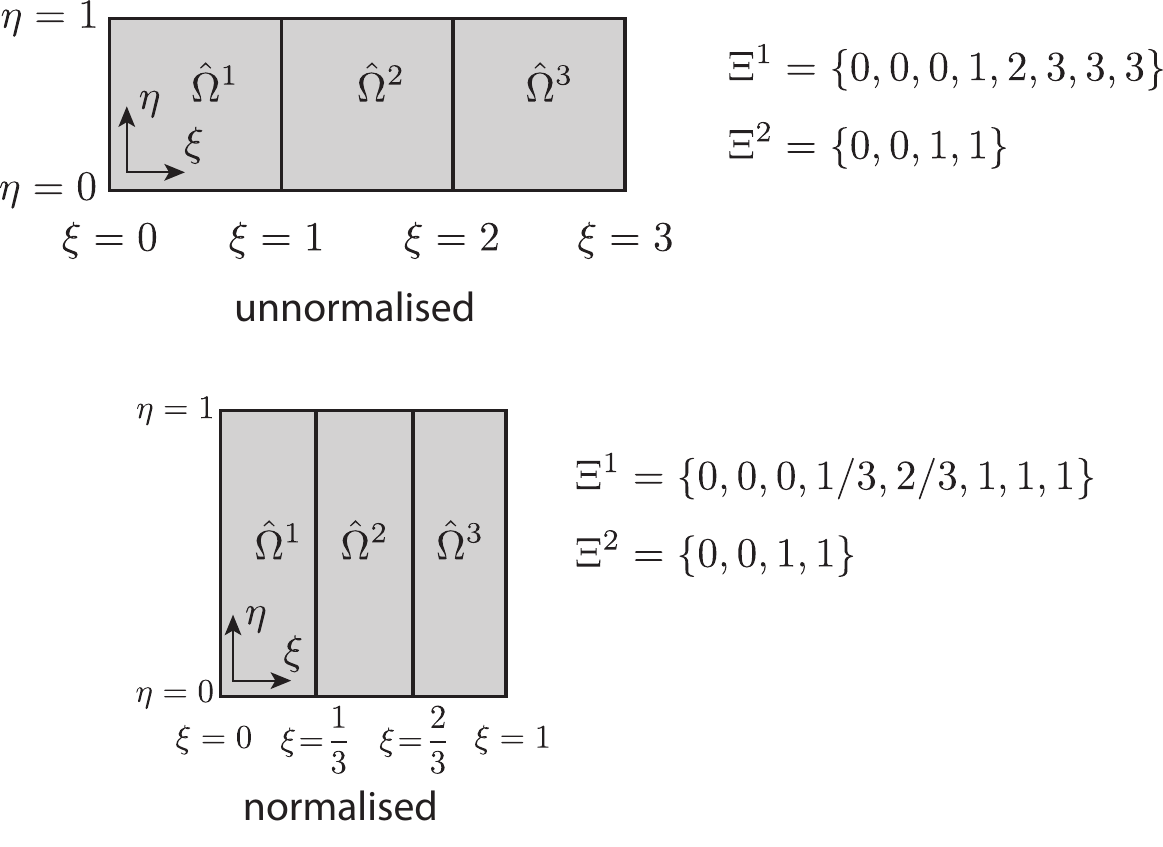}
  \caption{Parametric space defined by non-zero knot intervals. All parametric spaces defined for a B-spline or NURBS patch can be normalised to a unit interval, unit square, unit cube in 1D, 2D, 3D respectively. Knot lines provide a natural definition of element boundaries.}
  \label{fig:parameter_space}
\end{figure}

Parametric space (sometimes referred to as the `pre-image' of the
NURBS mapping) is formed by considering only the non-zero intervals
between knot values. For the knot vectors considered previously, the
parametric space is illustrated in Fig.~\ref{fig:parameter_space}
which can subsequently be reduced to a unit square through  appropriate
normalisation. All parametric spaces can be reduced to a unit
interval ($d_p=1$), square ($d_p=2$) or cube ($d_p=3$) in this manner.  We define the 
parametric space as $\hat{\Omega} \subset \mathbb{R}^{d_p}$ with a associated set of
parametric coordinates $\boldsymbol{\xi} =
(\xi,\eta, \zeta) = (\xi^1, \xi^2, \xi^3)\in \hat{\Omega}$ ($d_p=3$). If normalisation is performed,
 $\hat{\Omega} = [0,1]^{d_p}$.

Fig.~\ref{fig:parameter_space} also reveals that regions bounded by knot lines with non-zero parametric area lead to a natural definition of element domains. More formally, a set $\mathcal{S} \subset \Xi$ of unique knot values can be defined as
\begin{equation}
  \label{eq:unique_knot_vector_defn}
  \mathcal{S}= \{\xi_1,\xi_2,\cdots,\xi_{n_{s}} \} \quad \xi_i \neq \xi_{i+1} \,\, \textrm{for}  \,\, 1\leq i \leq n_s - 1
\end{equation}
where $n_s$ is the number of unique knot values. This is generalised
to $\mathcal{S}^i \subset \Xi^i$ which represents the unique knot
values for each parametric direction $i=1,2,\dots d_p$. Elements can now be defined in the general multivariate case as
\begin{align}
  \label{eq:parametric_element_defn}
 \hat{\Omega}^e = [\xi_i, \xi_{i+1}] \otimes [\eta_j, \eta_{j+1}] \otimes [\zeta_k, \zeta_{k+1}]  \quad &1\leq i \leq n_{s}^1 - 1,\\
&1\leq j \leq n_{s}^2 - 1,\notag\\
&1\leq k \leq n_{s}^3 - 1,\notag\\
&\xi_i\in\mathcal{S}^1, \,\eta_j\in\mathcal{S}^2, \,\zeta_k\in\mathcal{S}^3 \notag
\end{align}
where $n_{s}^1,n_{s}^2$ and $n_{s}^3$  represent the number unique knots in the $\xi$,$\eta$ and $\zeta$ parametric directions respectively. This leads to a natural numbering scheme for elements over a patch as 
\begin{equation}
  \label{eq:element_numbering}
  e = k( n_{s}^2 - 1 )( n_{s}^1 - 1) +  j( n_{s}^1 - 1) + i. 
\end{equation}
Eq.~\eqref{eq:parametric_element_defn} and
\eqref{eq:element_numbering} can be simplied accordingly for $d_p=1,2$.

\subsubsection{Physical space}
\label{sec:physical-space}

\begin{figure}
  \centering
  \includegraphics[width=0.6\textwidth]{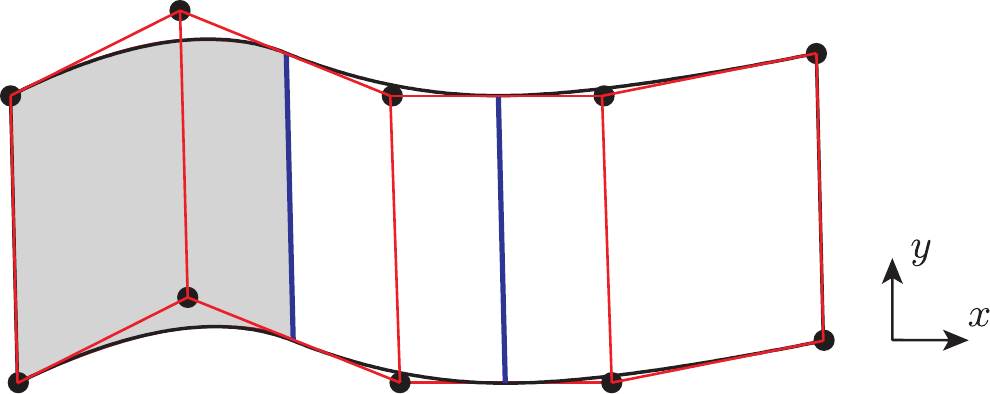}
  \caption{A 2D NURBS surface defined for knot vectors $\Xi^1 =
    \{0,0,0,1,2,3,3,3\}$, $\Xi^2 = \{0,0,1,1\}$. The control mesh is
    shown in red with control points denoted by black circles. Knot
    lines shown in blue indicate element boundaries.}
  \label{fig:physical_space}
\end{figure}

The B-spline and NURBS mappings of Eqs.~\eqref{eq:Bspline-surface1} and
\eqref{eq:NURBS-surface1} transform coordinates in parameter space to
physical space $\Omega \subset \mathbb{R}^{d_s}$. For
three-dimensional domains, we associate a coordinate system $\mathbf{x} = (x,y,z)
=(x^1, x^2, x^3)$ for physical space, which appropriate modifications for one-
and two-dimensional problems. Fig.~\ref{fig:physical_space} illustrates a NURBS mapping
for the parametric space shown in Fig.~\ref{fig:parameter_space} for
an arbitrary set of control points and weights.  The control grid
(which defines the connectivity between control points) is also
shown. The non-interpolatory nature of control points in the interior
of the domain is evident, and represents a notable difference over
conventional Lagrangian meshes. 

\subsubsection{Parent space}
\label{sec:parent-space}

The previous three spaces are inherent to B-splines and NURBS, but for analysis to be performed we
require the definition of an additional space, commonly referred to as
parent space $\tilde{\Omega} = [-1,1]^{d_p}$. This is required for the
use of numerical integration routines which are often defined over the interval
$[-1,1]$. Parent space coordinates are denoted as $\tilde{\boldsymbol{\xi}} =
(\tilde{\xi}, \tilde{\eta}, \tilde{\zeta}) = (\tilde{\xi}^1,
\tilde{\xi}^2, \tilde{\xi}^3)$ with corresponding simplifications for $d_p=1,2$. 

\subsection{Isogeometric formulation}
\label{sec:isog-form}

Before outlining the details of isogeometric analysis in a FE context, it is instructive to consider the similarities and differences of IGA over conventional discretisation technology. Lagrangian basis functions are most commonly used to discretise both the geometry and unknown fields in an \textit{isoparametric} fashion.  In this way, exactly the same basis functions are used for both. For the majority of cases, the geometry is always approximated incurring a geometrical error which may lead to erroneous results, particularly for highly oscillatory problems.  In addition, once a discretisation is generated from a given CAD model, the geometry information that has been lost can never be retrieved, forming a one-way process from discretisation to analysis. This has serious repercussions for efficiency, attenuated by the iterative nature of design.

Isogeometric analysis also makes use of an isoparametric formulation, but a key difference over its Lagrangian counterpart is the use of basis functions generated by CAD to discretise both the geometry and unknown fields. Not only is task of discretisation (meshing) greatly reduced or eliminated entirely, but a direct link is made with CAD, forming a bi-directional process. In addition, the use of CAD discretisations ensures that the \textit{exact} geometry is used at all stages of analysis, incurring no geometrical error. But the key compelling feature of IGA is the unified nature of design and analysis.

We can summarise these points as:
\begin{itemize}
  \item \textbf{Conventional finite element analysis}: the basis which is chosen to approximate the unknown field
	  is also used to approximate the known geometry. This most
          commonly takes the form of (low order) Lagrangian basis functions. In most cases the
          geometry is only approximated. CAD and analysis are disparate.
  \item \textbf{Isogeometric analysis}: the basis is generated by CAD  which captures the
	  geometry exactly. This basis is also used to approximate the unknown
          field. Refinement may be required for the unknown fields,
          but the exact geometry is maintained at all stages of analysis. CAD and analysis are combined to form a unified process.
\end{itemize}

\subsection{Isogeometric discretisation}

The B-spline and NURBS discretisations outlined in Section~\ref{sec-nurbs} are written in terms of parametric coordinates, but to use such discretisations for analysis, we must provide a mapping that allows us to operate at the parent element level. We will outline the appropriate mappings in Section~\ref{sec:mappings}, but for now let us assume that B-spline and NURBS basis functions can be written in terms of parent coordinates. This allows us to state the isoparametric discretisation used to approximate both the geometry and fields in IGA. For a given element $e$, the geometry is expressed as
\begin{equation}
  \label{eq:iga_geometry_discretisation}
  \mathbf{x}^e(\tilde{\boldsymbol{\xi}}) = \sum_{a=1}^{n_{en}} \vm{P}_a^e R_a^e(\tilde{\boldsymbol{\xi}})
\end{equation}
where $a$ is a local basis function index, $n_{en} = (p+1)^{d_p}$ is the number of non-zero basis functions over element $e$ and $\vm{P}_a^e$,$R_a^e$ are the control point and NURBS basis function  associated with index $a$ respectively. We employ the commonly used notation of an element connectivity mapping \cite{hughes-fem-book} which translates a local basis function index to a global index through
\begin{equation}
  \label{eq:element_connectivity_array}
  A = \textrm{IEN}( a, e ).
\end{equation}
Global and local control points are therefore related through $\vm{P}_A \equiv \vm{P}_{\textrm{IEN}(a,e)} \equiv \vm{P}_a^e$ with similar expressions for $R_a^e$.  A field $\vm{u}(\mathbf{x})$ which governs our relevant PDE can also be discretised in a similar manner to \eqref{eq:iga_geometry_discretisation} as
\begin{equation}
  \label{eq:iga_field_discretisation}
  \vm{u}^e(\tilde{\boldsymbol{\xi}}) =  \sum_{a=1}^{n_{en}} \vm{d}_a^e R_a^e(\tilde{\boldsymbol{\xi}})
\end{equation}
where $\vm{d}^e_a$ represents a control (nodal) variable. In contrast to conventional discretisations, these coefficients are not in general interpolatory at nodes. This is similar to the case of meshless
methods built on non-interpolatory shape functions such as the moving least squares
(MLS) \cite{efg-nayroles,NME:NME1620370205,nguyen_meshless_2008}.

\subsection{Mappings (change of variables)}
\label{sec:mappings}

\begin{figure}
  \centering
  \includegraphics{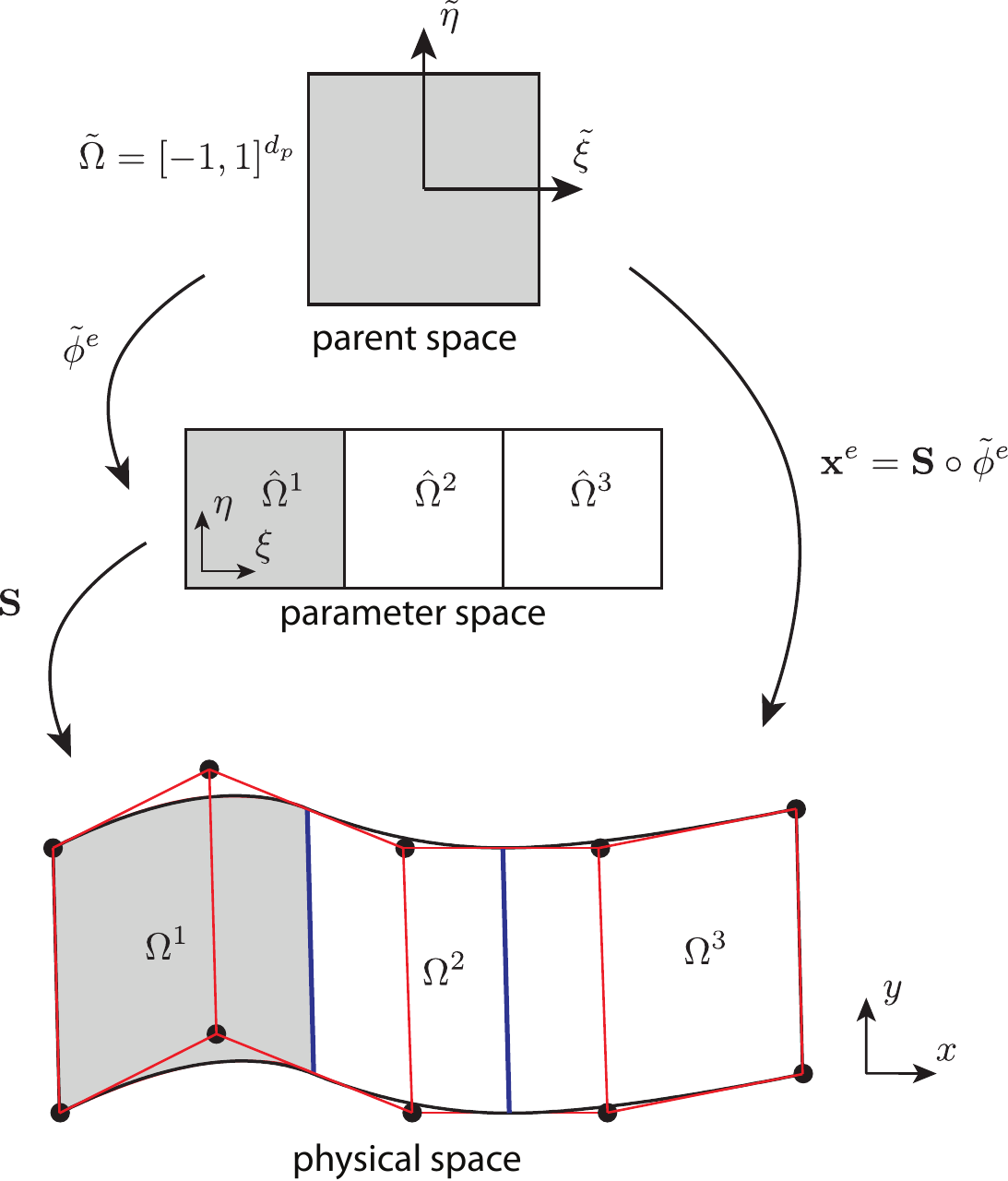}
  \caption{Diagrammatic interpretation of mappings from parent space through parametric space to physical space. }
  \label{fig:iga_mappings}
\end{figure}

The use of NURBS basis functions for discretisation introduces the concept of parametric space which is absent in conventional FE implementations. The consequence of this additional space is that an additional mapping must be performed to operate in parent element coordinates. As shown in Fig.~\ref{fig:iga_mappings}, two mappings are considered for IGA with NURBS: a mapping $\tilde{\phi}^e: \tilde{\Omega} \to \hat{\Omega}^e$ and $\vm{S}: \hat{\Omega} \to \Omega$. The mapping $\vm{x}^e: \tilde{\Omega} \to \Omega^e$ is given by the composition $\vm{S}\circ \tilde{\phi}^e$. 

Taking the case $d_p = d_s = 2$,  an element defined by $\hat{\Omega}^e = [\xi_i, \xi_{i+1}]\otimes [\eta_i, \eta_{i+1}]$ is mapped from parent space to parametric space through 
\begin{align}
  \label{eq:phi_mapping}
  \tilde{\phi}^e(\tilde{\boldsymbol{\xi}})  =
\left\{ 
\begin{matrix}
\frac{1}{2}[(\xi_{i+1}-\xi_i)\tilde{\xi} + (\xi_{i+1}+\xi_i)]\\
 \frac{1}{2}[(\eta_{j+1}-\eta_j)\tilde{\eta} + (\eta_{j+1}+\eta_j)]
\end{matrix}
\right\}
\end{align}
with an associated Jacobian determinant given by
\begin{equation}
  \label{eq:jacob_determ_parent}
    |J_{\bsym{\tilde{\xi}}}| = \frac{1}{4}(\xi_{i+1}-\xi_i)(\eta_{j+1}-\eta_j).
\end{equation}

Similarly, the mapping from parametric space to physical space is given by the NURBS expressions of Eqs.~\eqref{eq:NURBS-curve}, \eqref{eq:nurbs_surface_simplified} and \eqref{eq:nurbs_solid_2}. In the case $d_p=d_s=2$, the Jacobian of transformation for this mapping is represented by the matrix
\begin{equation}
  \vm{J}_{\bsym{\xi}} = \left[\begin{array}{cc}
      \dfrac{\partial x}{\partial \xi} & \dfrac{\partial x}{\partial \eta}\\[2ex]
      \dfrac{\partial y}{\partial \xi} & \dfrac{\partial y}{\partial \eta}
  \end{array}\right]
  \label{eq:jacobian-map}
\end{equation}
in which the components are calculated as
\begin{equation}
  \label{eq:nurbs_spatial_derivs}
  \frac{\partial \vm{x}}{\partial \boldsymbol{\xi}} = \sum_{a=1}^{n_{en}} \vm{P}_a^e \frac{\partial R_a^e(\boldsymbol{\xi})}{\partial \boldsymbol{\xi}}.
\end{equation}
The associated Jacobian determinant is denoted by $|J_{\bsym{\xi}}|$.

The mapping $x^e: \tilde{\Omega} \to \Omega^e$, formed through the composition of the previous mappings, can be written as
\begin{align}
  \label{eq:final_mapping}
  \vm{x}^e(\tilde{\boldsymbol{\xi}})  &= \left. \left( \sum_{A=1}^{N} \vm{P}_A R_A\left(\tilde{\phi}^e(\tilde{\boldsymbol{\xi}})\right) \right)\right|_e\\
&= \left. \left( \sum_{a=1}^{n_{en}} \vm{P}_a R_a(\tilde{\phi}^e(\tilde{\boldsymbol{\xi}})) \right)\right|_e\\
&= \sum_{a=1}^{n_{en}} \vm{P}_a^e R_a^e(\tilde{\boldsymbol{\xi}}) \\
\end{align}
where $N$ is the number of global basis functions and $(\cdot)|_e$ denotes that the expression $(\cdot)$ is restricted to element $e$. 
The Jacobian determinant for this mapping is given by 
\begin{equation}
  \label{eq:jacob_det_final_mapping}
  |J| = |J_{\bsym{\xi}}| | J_{\bsym{\tilde{\xi}}}|.
\end{equation}

With this final mapping and Jacobian determinant, it is possible to integrate a function $f: \Omega \to \mathbb{R}$ over the physical domain as
\begin{align}
  \label{eq:domain_integration}
  \int_\Omega f(\vm{x})\, \mathrm{d}\Omega &= \sum_{e=1}^{n_{el}} \int_{\Omega^e} f(\vm{x})\, \mathrm{d}\Omega\\
&= \sum_{e=1}^{n_{el}} \int_{\hat{\Omega}^e} f(\vm{x}(\boldsymbol{\xi}))  \, |J_{\bsym{\xi}}|  \, \mathrm{d}\hat{\Omega}\\
&= \sum_{e=1}^{n_{el}} \int_{\tilde{\Omega}} f(\vm{x}(\tilde{\phi}^e(\tilde{\boldsymbol{\xi}}))) \,  |J_{\bsym{\xi}}| \, | J_{\bsym{\tilde{\xi}}}| \, \mathrm{d}\tilde{\Omega}\\
&=  \sum_{e=1}^{n_{el}} \int_{\tilde{\Omega}} f(\tilde{\boldsymbol{\xi}}) \,|J|\mathrm{d}\tilde{\Omega}
\end{align}
with the final integral in a suitable form for application
of standard Gauss-Legendre
quadrature (hereafter named Gaussian quadrature).
As detailed in \cite{hughes-fem-book} for Lagrangian basis functions,
a rule of $(p+1)\times(q+1)$ Gaussian quadrature can be applied for
two-dimensional elements in which $p$ and $q$ denote the orders of
the chosen basis functions in the $\xi$ and $\eta$ direction. The same
procedure is also used for NURBS basis functions in the present work,
although it should be emphasised that Gaussian quadrature is not optimal for IGA.
Research is currently focussed on optimal integration techniques such
as that in \cite{hughes_efficient_2010,Auricchio201215} in which an
optimal quadrature rule, known as the half-point rule, has been applied.

Spatial derivatives of basis functions are also required for element assembly algorithms, and are calculated as
\begin{equation}
  \left[ \begin{array}{c} \dfrac{\partial R_a^e}{\partial \vm{x}}\end{array}
    \right] = \left[ \begin{array}{c} \dfrac{\partial R_a^e}{\partial \vm{\boldsymbol{\xi}}}\end{array} \right] 
\left[ \begin{array}{cc} \dfrac{\partial \xi}{\partial x} & \dfrac{\partial \xi}{\partial y} \\[2ex]
	\dfrac{\partial \eta}{\partial x} & \dfrac{\partial \eta}{\partial y}\end{array}
    \right] = \left[ \begin{array}{c} \dfrac{\partial R_a^e}{\partial \vm{\boldsymbol{\xi}}}\end{array} \right]    \vm{J}_{\bsym{\xi}}^{-1},
  \label{eq:sd}
\end{equation}
with the derivatives $\partial R_a^e/\partial \vm{\boldsymbol{\xi}}$ obtained through Eq.~\eqref{eq:nurbs-derivatives-1}.

\subsection{One-dimensional IGA formulation}\label{1dAssembly}

To illustrate the use of the isogeometric concept in a finite element
context,  a one-dimensional IGA formulation is developed through an example. Efforts are made to maintain the notation adopted in \cite{hughes-fem-book} and \cite{cottrel_book_2009}. In the following, we define the domain $\Omega \subset \mathbb{R}$ with boundary $\Gamma \equiv \partial \Omega$. The boundary is partitioned as $\Gamma = \overline{\Gamma_D \cup \Gamma_N}$, $\Gamma_N \cap \Gamma_D = \emptyset$ where $\Gamma_D$ and $\Gamma_N$ denote the Dirichlet and Neumann boundaries respectively with an overline representing a closed set. Robin boundary conditions are not considered in the present work. 

Letting $\Omega = (0,1)$ and $\Gamma_D = \Gamma = \{0,1\}$ with $\Gamma_N=\emptyset$, we seek the solution $u : \overline{\Omega} \to \mathbb{R}$ such that
\begin{equation}
\frac{d^2 u (x)}{dx^2} +  b(x) = 0, \quad x \in \overline{\Omega},
	\label{eq:1d-iga-pde}
\end{equation}
with the Dirichlet boundary conditions specified through $g: \Gamma_D \to \mathbb{R}$ as
\begin{equation}
  \label{eq:1d-iga-bcs}
\quad g(0) = 0,  \quad g(1) = 0.
\end{equation}
Choosing $b(x)=x$, it can be shown that the exact solution to Eq.~\eqref{eq:1d-iga-pde} subject to Eq.~\eqref{eq:1d-iga-bcs} is given by
\begin{equation}
	u(x) = -\frac{1}{6}x^3 + \frac{1}{6}x.
\end{equation}

To construct the Galerkin formulation of the preceding problem, the infinite dimensional spaces $\mathscr{U}$ and $\mathscr{V}$ which represent the usual trial and test spaces respectively\footnote{Specifically, $\mathscr{U} = \{ u : u \in H^1(\Omega), u=g \text{ on } \Gamma_D \}$, $\mathscr{V} = \{w : w \in H^1(\Omega), w=0 \text{ on } \Gamma_D\}$ where $H^1(\Omega)$ represents a Sobolev space of order $1$} are required. By multiplying Eq. \eqref{eq:1d-iga-pde} by a test function $w \in \mathscr{V}$ and applying integration by parts, the weak form of the problem reads: find $u \in \mathscr{U}$, such that 
\begin{equation}
  \label{eq:1d_poisson_weak_form}
\int_\Omega\, \frac{d w}{dx}\frac{du}{dx} \, \mathrm{d}\Omega = \int_\Omega wb\, \mathrm{d}\Omega \quad \text{ for all } w \in \mathscr{V}.
\end{equation}
This is reduced to a finite-dimensional problem by creating finite-dimensional subspaces $\mathscr{U}^h \subset \mathscr{U}$ and $\mathscr{V}^h \subset \mathscr{V}$ both formed through a NURBS basis.  Letting $g^h \in \mathscr{U^h}$ denote the finite-dimensional representation of $g$, the solution we wish to seek can be written as
\begin{equation}
  \label{eq:1d_iga_discrete_soln}
  u^h = v^h + g^h
\end{equation}
which is valid for all $v^h \in \mathscr{V}^h$. Expression~\eqref{eq:1d_iga_discrete_soln} underpins a key concept of the Galerkin formulation in that the solution $u^h$ and test (weighting) function $w^h$ are constructed from the same space of functions. We also see that the boundary condition data is `built in' to the solution by the inclusion of the term $g^h$. 

The approximation to our solution can be written in terms of NURBS basis functions as
\begin{equation}
  \label{eq:1d_iga_nurbs_approx_soln}
  u^h = \sum_{B=1}^{n_{eq}} d_B R_B(x) + \sum_{B=n_{eq}+1}^{n_{np}} g_B R_B(x)
\end{equation}
where $n_{eq}$ is the number of equations (number of unknowns) and $n_{np}$ is the total number of nodal points (total degrees of freedom) in the system. $\{d_B\}_{B=1}^{n_{eq}}$ represents the set of all unknown control variables and $\{g_B\}_{B=n_{eq}+1}^{n_{np}}$ is the set of known Dirichlet control variables. Likewise, the test function is discretised as
\begin{equation}
  \label{eq:1d_iga_discrete_weighting}
  w^h = \sum_{A=1}^{n_{eq}} c_A R_A(x)
\end{equation}
Substitution of Eqs.~\eqref{eq:1d_iga_nurbs_approx_soln} and \eqref{eq:1d_iga_discrete_weighting} into 
Eq.~\eqref{eq:1d_poisson_weak_form} and noting that the values of the set $\{c_A\}$ are arbitrary, the discrete form is now written as
\begin{equation}
  \label{eq:1d_iga_discrete_form}
 \left( \sum_{A=1}^{n_{eq}}  \sum_{B=1}^{n_{eq}} \int_\Omega \frac{dR_A}{dx}  \frac{dR_B}{dx}\, \mathrm{d}\Omega \right) \, d_B = \sum_{A=1}^{n_{eq}} \int_\Omega R_A b \, \mathrm{d}\Omega - \left( \sum_{A=1}^{n_{eq}}  \sum_{B=n_{eq}+1}^{n_{np}} \int_\Omega \frac{dR_A}{dx}  \frac{dR_B}{dx}\, \mathrm{d}\Omega \right) \, g_B
\end{equation}
or, by defining the element stiffness matrix and element force vector as
\begin{align}
  \label{eq:1d_elemental_iga_matrix}
  \mathbf{K}_{AB} &=  \int_\Omega \frac{dR_A}{dx}  \frac{dR_B}{dx}\, \mathrm{d}\Omega\\
\label{eq:1d_elemental_iga_vector}
\mathbf{F}_A &=   \int_\Omega R_A b \, \mathrm{d}\Omega -  \sum_{B=n_{eq}+1}^{n_{np}} \int_\Omega \frac{dR_A}{dx}  \frac{dR_B}{dx}\, \mathrm{d}\Omega,
\end{align}
\eqref{eq:1d_iga_discrete_form} can expressed in matrix notation as
\begin{equation}
  \label{eq:1d_iga_global_matrix}
  \mathbf{K} \mathbf{d} = \mathbf{F}
\end{equation}
where 
\begin{align}
  \label{eq:1d_iga_matrix_notation}
  \mathbf{K} &= [\mathbf{K}_{AB}]\\
  \mathbf{d} &= \{ d_B\}\\
  \mathbf{F} &= \{ \mathbf{F}_A \} \quad A,B = 1,2,\ldots n_{eq}.
\end{align}
Eq.~\eqref{eq:1d_iga_global_matrix} is in a form amenable for computation. 

\begin{figure}[htbp]
  \centering 
  \psfrag{xi}{$\xi$}
  \psfrag{n1}{$R_1$} \psfrag{n2}{$R_2$}
  \psfrag{n3}{$R_3$} \psfrag{n4}{$R_4$}
  \psfrag{xi1}{$\xi_1=\xi_2=\xi_3=0$}
  \psfrag{xi3}{$\xi_5=\xi_6=\xi_7=1$}
  \psfrag{xi2}{$\xi_4$}
  \psfrag{b1}{$\vm{P}_1$} \psfrag{b2}{$\vm{P}_2$}
  \psfrag{b3}{$\vm{P}_3$} \psfrag{b4}{$\vm{P}_4$}
  \psfrag{e1}{1} \psfrag{e2}{2}
  \psfrag{curve}[c]{A quadratic B-spline curve with 4 control points}
  \psfrag{para}[c]{Parametric space}
  \psfrag{basis}[c]{4 quadratic basis functions}
  \includegraphics[width=0.65\textwidth]{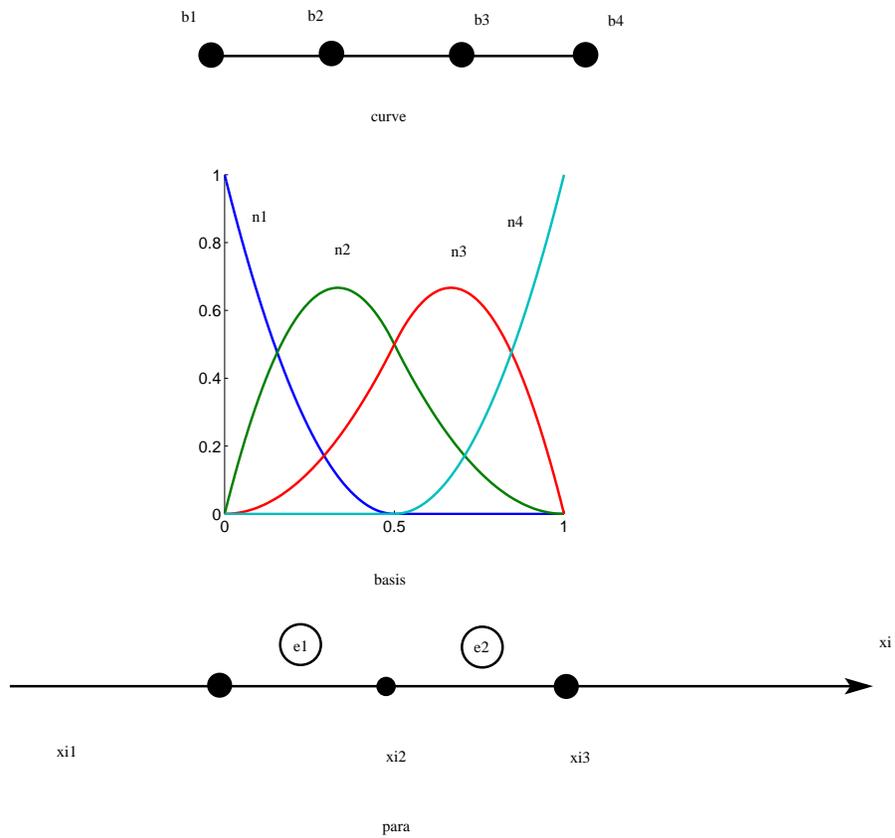}
  \caption{One dimensional isogeometric analysis example: exact geometry, quadratic
           basis functions and mesh in the parametric space. IEN(:,1)=$[1\;2\;3]$ and IEN(:,2)=$[2\;3\;4]$. } 
  \label{fig:assembly1d} 
\end{figure}

\subsubsection{Assembly of system of equations}
\label{sec:assembly-syst-equat}

The global stiffness matrix $\mathbf{K}$ and global force vector $\mathbf{F}$ are assembled by looping over each element (non-zero knot interval) and inserting the appropriate element entries into their relevant rows and columns in the global system of equations. To allow integration over each element, a parent coordinate system is adopted making use of the mappings and Jacobian matrices outlined in Section~\ref{sec:mappings}.  For example, for a particular element $e$ with local basis function indices $a,b=1,2,\dots,p+1$, we can compute a local element stiffness matrix $\mathbf{k}^e_{ab}$ using parent coordinates as
\begin{equation}
  \label{eq:1d_iga_element_stiffness_parent}
  \mathbf{k}^e_{ab} = \int_{\tilde{\Omega}} \frac{dR_a^e}{dx} \frac{dR_b^e}{dx} |J| \, \mathrm{d}\tilde{\Omega}
\end{equation}
where use has been made of expressions \eqref{eq:sd} and \eqref{eq:jacob_det_final_mapping}. This can be inserted into the global stiffness matrix through the element mapping given by Eq.~\eqref{eq:element_connectivity_array}.

To make this discussion more concrete, the element assembly process for a simple quadratic ($p=2$) NURBS discretisation with $\Xi=\{0,0,0,0.5,1,1,1\}$ is outlined. We apply the discretisation to the problem defined by 
Eqs.~\eqref{eq:1d-iga-pde} and \eqref{eq:1d-iga-bcs}. Fig.~\ref{fig:assembly1d} illustrates the relevant NURBS discretisation including control points, basis functions and element definitions. Note that there are $p+1$ non-zero basis functions over an element $e$ given by $\{R_a^e\}_{a=1}^{p+1}$. We can therefore summarise the element data for this discretisation as
\begin{equation}
  \begin{array}{ccccc}
	  \textrm{element} & \textrm{knot interval} & \textrm{non-zero basis functions} &
	  \textrm{control variables} & \textrm{control points}\\
	  1       & [\xi_3,\xi_4] & R_1,R_2,R_3 & d_1, d_2, d_3 &
	  \vm{P}_1,\vm{P}_2, \vm{P}_3\\
	  2       & [\xi_4,\xi_5] & R_2,R_3,R_4 & d_2, d_3, d_4 &
	  \vm{P}_2,\vm{P}_3,\vm{P}_4\\
  \end{array}
  \label{eq:dep}
\end{equation}
The implementation required to compute the element stiffness matrix  is shown in Box \ref{box-k-example}. The code
follows the same standard FEM routines presented in \cite{chessa-fem} where the notation $\norm{\cdot}$ is used to signify the $L^2$ norm\footnote{That is, for a vector $\mathbf{x} = (x_1, x_2,\ldots,x_n)$, $||\mathbf{x}|| = (x_1^2 + x_2^2 +\ldots+x_n^2)^{1/2}$}.

\begin{Fbox}
	\caption{Element stiffness matrix evaluation for element $\hat{\Omega}^e=[\xi_i,\xi_{i+1}]$}
  \begin{enumerate}
	  \item $\vm{P}=[\vm{P}_1;\vm{P}_2;\vm{P}_3]$
	  \item $sctr=[\text{IEN}(1,e)\; \text{IEN}(2,e)\; \text{IEN}(3,e)]$
	  \item Set $\vm{k}^e=0$
	  \item Loop over Gauss points (GPs) $\{\tilde{\xi}_j, \tilde{w}_j\} \quad j=1,2,\ldots,n_{gp}$$^*$
  \begin{enumerate}
	  \item Compute parametric coordinate
		  $\xi=\tilde{\phi}^e(\tilde{\xi}_j)$ (see Eq.~\eqref{eq:phi_mapping})
	  \item Compute derivatives $R_{a,\xi}^e$ ($a=1,2,3$) at point $\xi$
	  \item Define vector $\vm{R}_\xi=[R^e_{1,\xi}\; R^e_{2,\xi}\;
		  R^e_{3,\xi}]$
	  \item Compute $|J_\xi|=||\vm{R}_{\xi}\vm{P}||$
	  \item Compute $|J_{\tilde{\xi}}|=0.5(\xi_{i+1}-\xi_i)$
	  \item Compute shape function derivatives
		  $\vm{R}_{x}=J_\xi^{-1}\vm{R}_\xi\trans$ 
	  \item $\vm{k}^e = \vm{k}^e +   \vm{R}_x \vm{R}_x\trans \,|J_\xi| \, |J_{\tilde{\xi}}|\, \tilde{w}_j$
  \end{enumerate}
  \item End loop over GPs
  \item Assemble $\vm{k}^e$ into the global matrix as
	  $\vm{K}(sctr,sctr)=\vm{K}(sctr,sctr)+\vm{k}^e$
  \end{enumerate}
  $^*$ $\tilde{w}_j$ denotes a Gauss point weight and $n_{gp}$ is the total number of Gauss points.
  \label{box-k-example}
\end{Fbox}

\begin{figure}[htbp]
	\centering
	\subfloat[Two element NURBS mesh]{\includegraphics[width=0.31\textwidth]{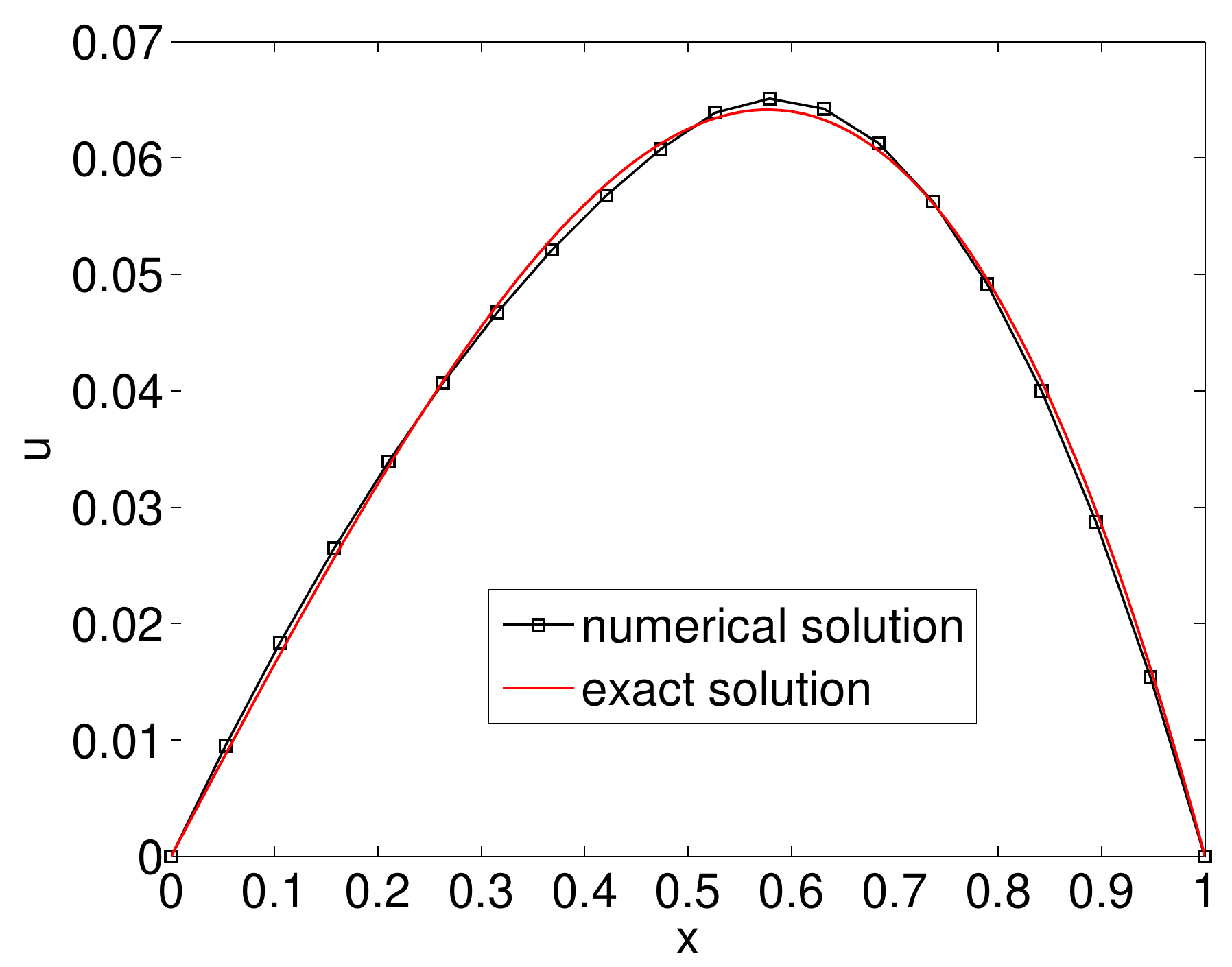}}\;\;\;\;\;
	\subfloat[Four element NURBS mesh.]{\includegraphics[width=0.31\textwidth]{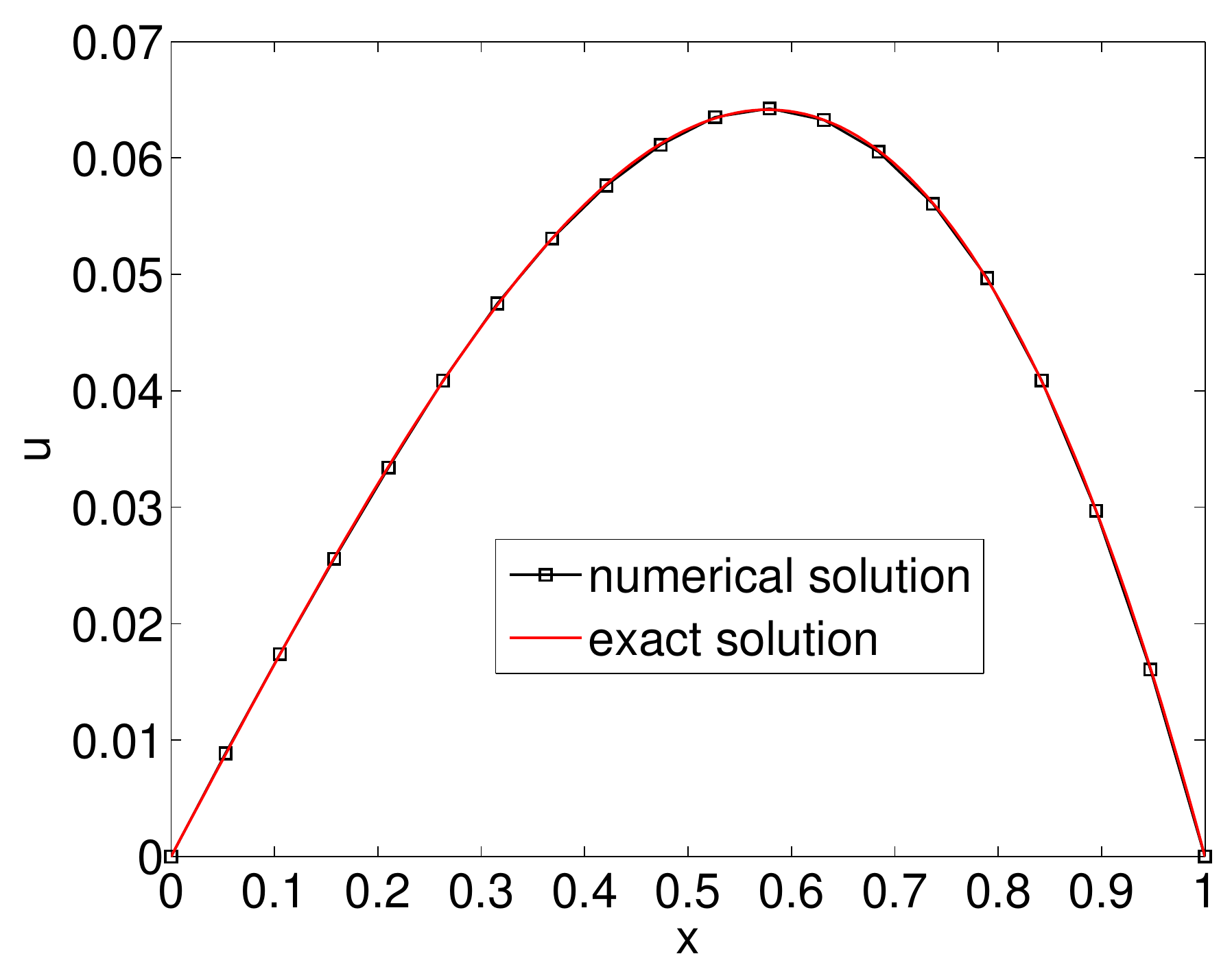}}\;\;\;\;\;
	\caption{Quadratic NURBS solution for one-dimensional Poisson equation.}
	\label{fig:iga-1d-result}
\end{figure}

Fig.~\ref{fig:iga-1d-result} illustrates the solution obtained for a NURBS discretisation for both two and four elements.  Convergence to the exact solution is observed, verifying the present IGA implementation.


In order to assess how high order B-spline elements behave for problems with localized
gradients, let us consider the following problem \cite{nguyen_meshless_2008}

\begin{equation}
u_{,xx}(x) +  b(x) = 0 \quad x \in [0,1]; \quad u(0) = 0,  \quad u(1) = 1,
	\label{eq:1d-pde}
\end{equation}
\noindent with

    \begin{equation}
         b(x) = \left\{\begin{array}{ll}
               \biggl \{ 2\alpha^2 - 4[\alpha^2(x-0.5)]^2 \biggr \} \exp
               \biggl \{ - [\alpha(x-0.5)]^2\biggr \} &  x \in [0.42,0.58] \\
               0 & \textrm{otherwise}
             \end{array}
      \right.
    \end{equation}

\noindent The exact solution of this problem is

\begin{equation}
u(x) = x + \exp \biggl \{ - [\alpha(x-0.5)]^2 \biggr \} x\in[0,1]
\end{equation}
\noindent We use a value of 50 for $\alpha$ and the exact solution exhibits a sharp peek
at location $x=0.5$. We are going to solve this problem using elements of order ranging from
one (linear elements) to five (quintic elements). To remove error in the numerical integration 
of the body force term, Eq. (\ref{eq:k-example}), 10 GPs were used for each element. We use the
$k$-refinement (to be discussed in detail in Section \ref{sec:hpk-refinement}) in building meshes
of different basis orders. The initial mesh consists of one single linear element with knot vector
$\Xi=\{0,0,1,1\}$. The parametrization is thus linear and after performing $k$-refinement, the
parametrization is still linear. Therefore, a point with $x=0.3$ corresponds to $\xi=0.3$ in the
parameter space. The Matlab file for this problem is \textbf{iga1DStrongGradient.m} in the 
\textbf{iga} folder of our source code.

\begin{figure}[htbp]
\centering
\subfloat[16 elements]{\includegraphics[width=0.46\textwidth]{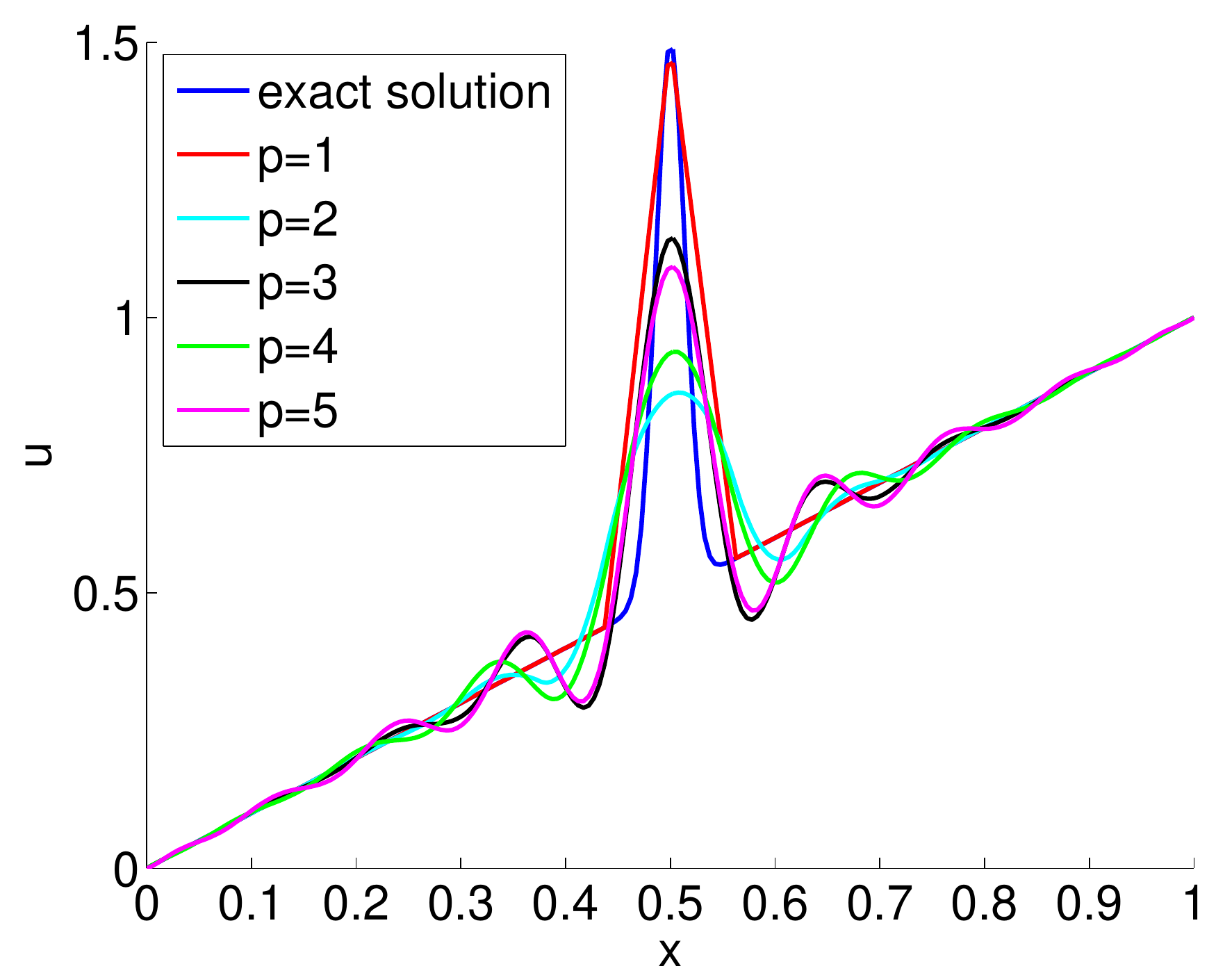}}\;
\subfloat[32 elements]{\includegraphics[width=0.46\textwidth]{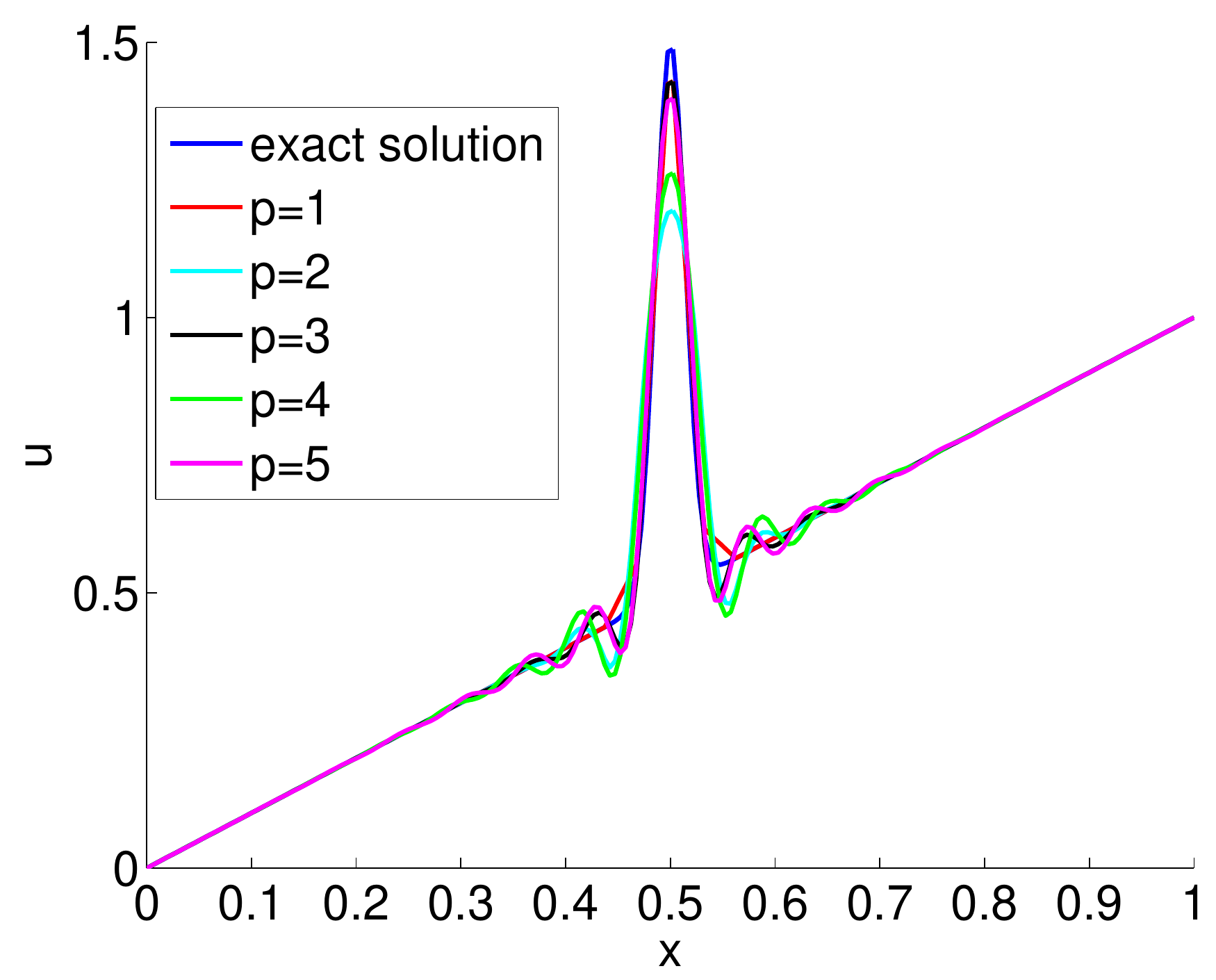}}\;
\subfloat[32 elements with $C^0$ continuity at $x=0.5$]{\includegraphics[width=0.46\textwidth]{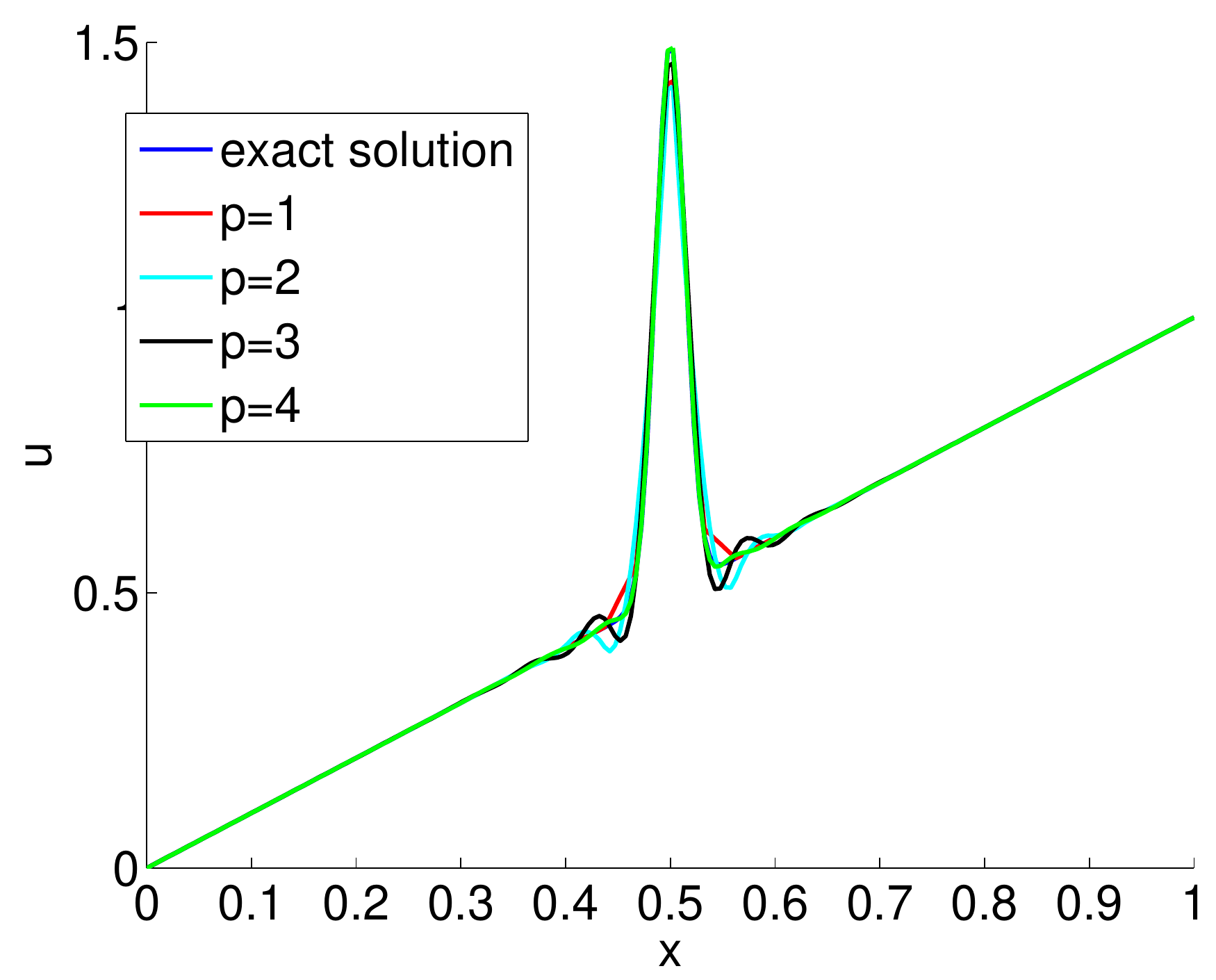}}\;
\subfloat[32 elements with $C^0$ continuity at $x=0.5$ (close up)]{\includegraphics[width=0.46\textwidth]{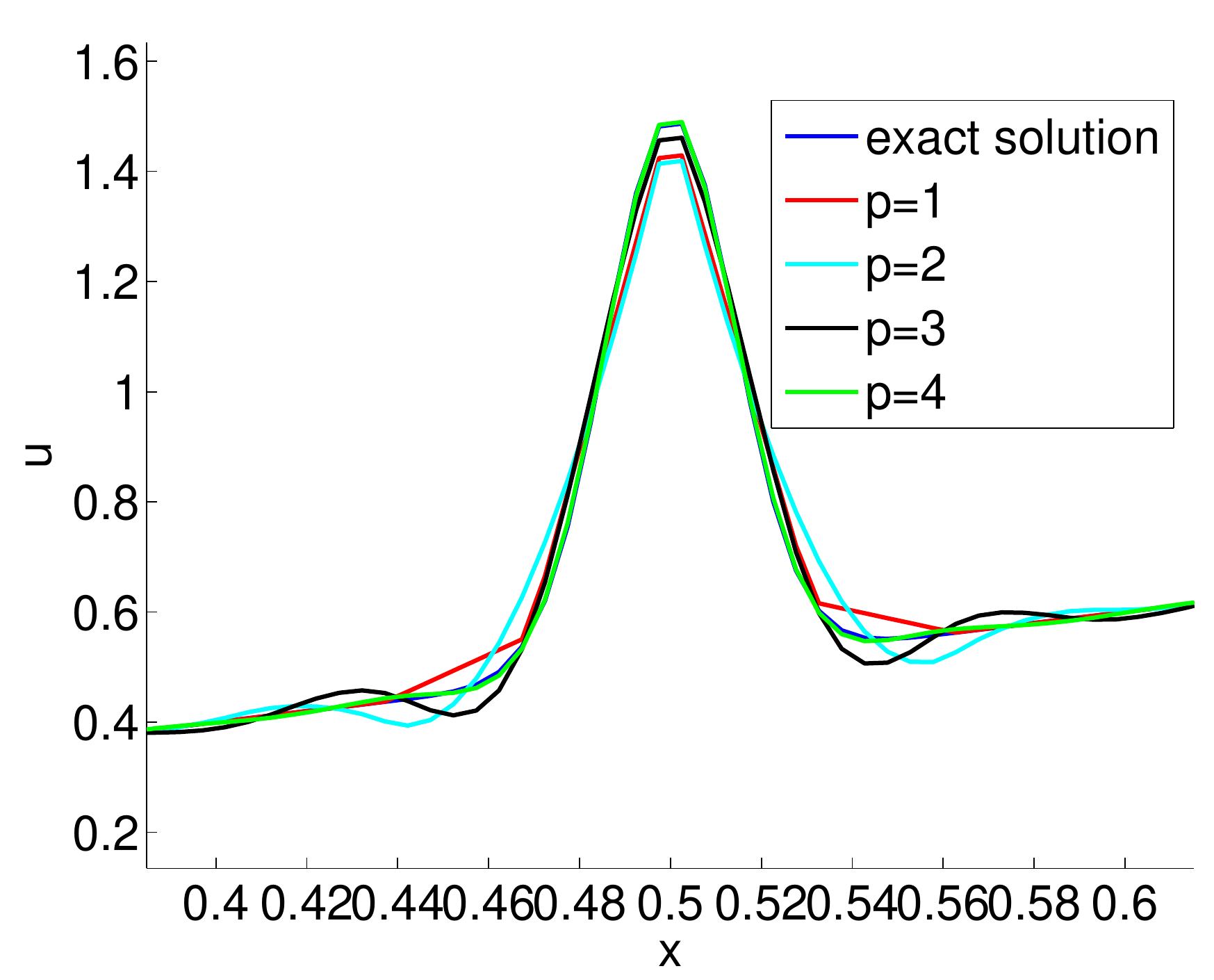}}\;
\caption{Comparison of the IGA result against the exact solution for the one dimensional PDE given 
   in Eq. (\ref{eq:1d-pde}).}
\label{fig:iga-1d-gradient}
\end{figure}

Figs. \ref{fig:iga-1d-gradient}a,b shows the results obtained with
meshes consisting of 16 
and 32 elements of $C^{p-1}$ continuity where $p$ denotes the B-spline basis order.
It is obvious that smooth basis is not suitable to problems with sharp gradients. Linear elements
which have a $C^0$ continuity at location $x=0.5$ (precisely at every knots)
give better results than high order B-spline elements.

\begin{figure}[htbp]
	\centering
	\subfloat[16 elements with $C^0$ continuity at $x=0.42,0.5,0.58$]
{\includegraphics[width=0.46\textwidth]{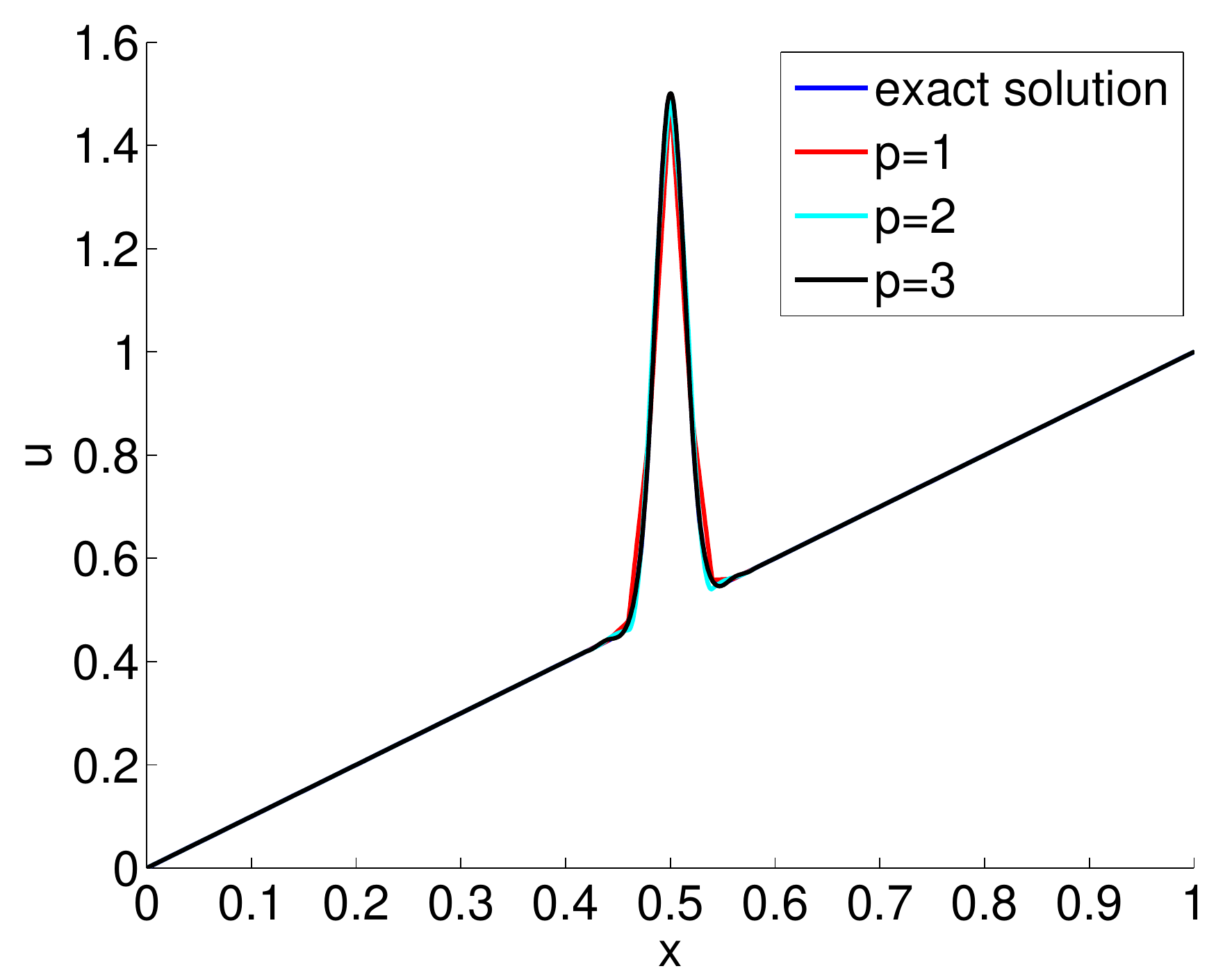}}\;
	\subfloat[16 elements with $C^0$ continuity at $x=0.42,0.5,0.58$ (close up)]
{\includegraphics[width=0.46\textwidth]{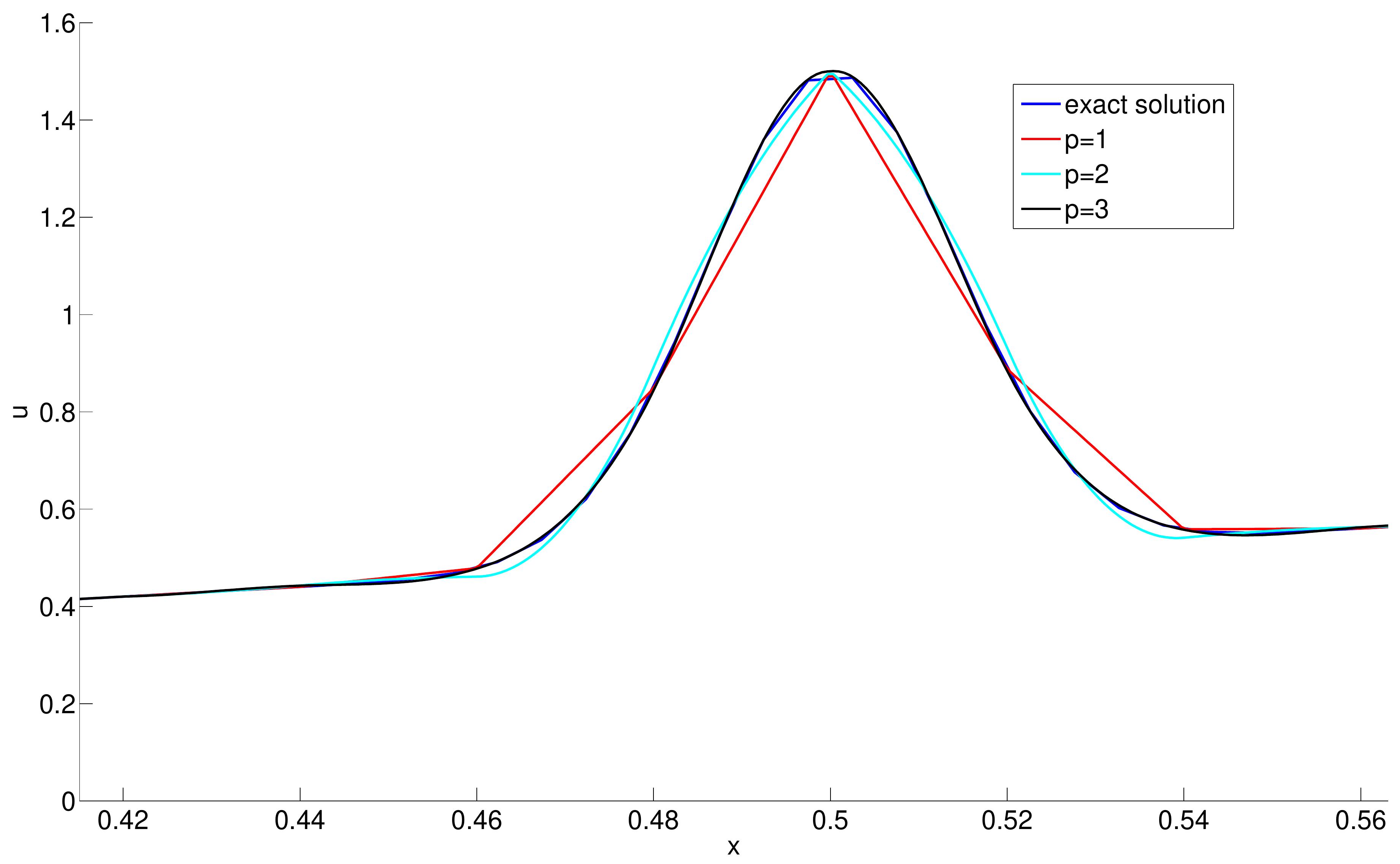}}
	\caption{Comparison of the IGA result against the exact solution
	for the one dimensional PDE given in Eq. (\ref{eq:1d-pde}).}
	\label{fig:iga-1d-gradient1}
\end{figure}

Knot insertion was made to insert the value $0.5$ $p$ times to the initial knot vector so that 
the basis is $C^0$ continuous at $x=0.5$. For example for $p=2$, a vector $\{0.5,0.5\}$ was inserted to
the knot vector. The results are given in Fig. \ref{fig:iga-1d-gradient}c,d where the peak was captured
better.
Finally $0.42,0.5,0.58$ were added to the knots $p$ times so that the basis are $C^0$ continuous
at locations $x=0.42,0.5,0.58$. The corresponding results are given in Fig. \ref{fig:iga-1d-gradient1}.
With only 16 cubic elements (25 CPs) the exact solution was well captured. This example showed the flexibility
of B-splines/NURBS--high order functions and any level of continuity can be easily achieved. 


\section{Elasticity: two-dimensional implementation}
\label{sec:elast-two-dimens}

The application of the FEM to elasticity is common and represents a
familiar language to many researchers.  We therefore outline a
two-dimensional implementation of IGA for linear elasticity,
highlighting the differences over conventional discretisations.

\subsection{Assembly process for two dimensional elastostatic analysis}
\label{sec:2DAssembly}

Consider a domain $\Omega $, bounded by $\Gamma $. The boundary is partitioned
into two sets: $\Gamma _u$ and $\Gamma _t $ with displacements prescribed on $\Gamma _u $
and tractions $\overline{\mathbf t}$ prescribed on $\Gamma _t$:
$\Gamma = \overline{\Gamma_t \cup \Gamma_u}$, $\Gamma_t \cap \Gamma_u = \emptyset$. 
The weak form of a linear elastostatics problem
is to find $\mathbf u$ in the trial space \footnote{contains $C^0$ functions}, such that for
all test functions $\delta \mathbf u$ in the test space \footnote{contains $C^0$ functions
vanishing on $\Gamma _u$},

\begin{equation}
\int_{\Omega} \boldsymbol{\varepsilon}(\mathbf u):
\mathbf D :\boldsymbol{\varepsilon}(\delta \mathbf u) \mathrm{d}
\Omega  = \int_{\Gamma_t} \overline{\mathbf t} \cdot \delta \mathbf u
\mathrm{d} \Gamma + \int_{\Omega} \mathbf b \cdot \delta \mathbf u
\mathrm{d} \Omega,
\label{eq:hang}
\end{equation}
where the elasticity matrix is denoted by $\vm{D}$ and $\vm{b}$
refers to a body force. Using the Galerkin method where the
same shape functions $R_a(\tilde{\bsym{\xi}})$ are used for both $\vm{u}$ and $\delta\vm{u}$,  we can write

\begin{equation}
  \vm{u}(\vm{x}) = \sum_{A=1}^{n_{np}}R_{A}(\tilde{\bsym{\xi}}) \vm{u}_{A}, \quad
  \delta \vm{u}(\vm{x}) = \sum_{A=1}^{n_{np}}R_{A}(\tilde{\bsym{\xi}}) \delta \vm{u}_{A},
\end{equation}
where $\vm{u}_A,\delta \vm{u}_{A}$ denote the nodal displacement and its variations, respectively
and $n_{np}$ is the total number of control points. In 2D, each control point has two
unknowns--the $x$ and $y$ displacements, hence one writes $\vm{u}_A=\{u_{xA},u_{yA}\}$.
Proper modification can be made for 3D problems.

Substitution of these approximations
into Eq. (\ref{eq:hang}) and using the arbitrariness of the nodal variations 
gives the standard discrete set of equations ${\bf K} \ {\bf u} = {\bf f}$ with

\begin{equation}
\mathbf K_{AB} = \int_{\Omega} \mathbf B_A^\mathrm{T} \mathbf D
\mathbf B_B \mathrm{d} \Omega, \quad
\mathbf f_{A} = \int_{\Gamma_t} R_A \mathbf{\overline{t}}
\mathrm{d} \Gamma + \int_{\Omega} R_A \mathbf b \mathrm{d} \Omega, \quad A,B=1,2,\ldots,n_{np}.
\end{equation}

\noindent In two dimensions, the strain-displacement matrix $\mathbf B_A$ is given by

\begin{equation}
   \mathbf B_A = \begin{bmatrix}
     R_{A,x} & 0 \\
     0 & R_{A,y} \\
     R_{A,y} & R_{A,x} \\
   \end{bmatrix},
     \label{eq:strain-disp}
\end{equation}

\noindent where the shape function derivatives are computed according to Eq. (\ref{eq:sd})
and $R_{A,x}\equiv dR_A/dx$.



\begin{figure}[htbp]
  \centering 
  \psfrag{N11}{$N_{1}(\xi)$}
  \psfrag{N21}{$N_{2}(\xi)$}
  \psfrag{N12}{$M_{1}(\eta)$} \psfrag{N22}{$M_{2}(\eta)$}
  \psfrag{N32}{$M_{3}(\eta)$} \psfrag{N42}{$M_{4}(\eta)$}
  \psfrag{xi1}{$\xi_{1,2}=0$}
  \psfrag{xi2}{$\xi_{3,4}=1$}
  \psfrag{eta1}{$\eta_{1,2,3}=0$}
  \psfrag{eta2}{$\eta_4=0.5$}
  \psfrag{eta3}{$\eta_{5,6,7}=1$}
  \psfrag{11}[c]{1} \psfrag{12}[c]{2}
  \psfrag{b11}{$\vm{B}_{11}(1)$}
  \psfrag{b12}{$\vm{B}_{12}(3)$}
  \psfrag{b13}{$\vm{B}_{13}(5)$}
  \psfrag{b14}{$\vm{B}_{14}(7)$}
  \psfrag{b21}{$\vm{B}_{21}(2)$}
  \psfrag{b22}{$\vm{B}_{22}(4)$}
  \psfrag{b23}{$\vm{B}_{23}(6)$}
  \psfrag{b24}{$\vm{B}_{24}(8)$}
  \psfrag{1}{1} \psfrag{2}{2} \psfrag{3}{3} \psfrag{4}{4}
  \psfrag{5}{5} \psfrag{6}{6} \psfrag{7}{7} \psfrag{8}{8}
  \psfrag{0}{0} \psfrag{0}{0}\psfrag{a}[c]{(a) Exact
  geometry}\psfrag{b}[c]{(b) Mesh in
  physical space}\psfrag{c}[c]{(c) Mesh in parametric space and basis functions}
  \includegraphics[width=0.8\textwidth]{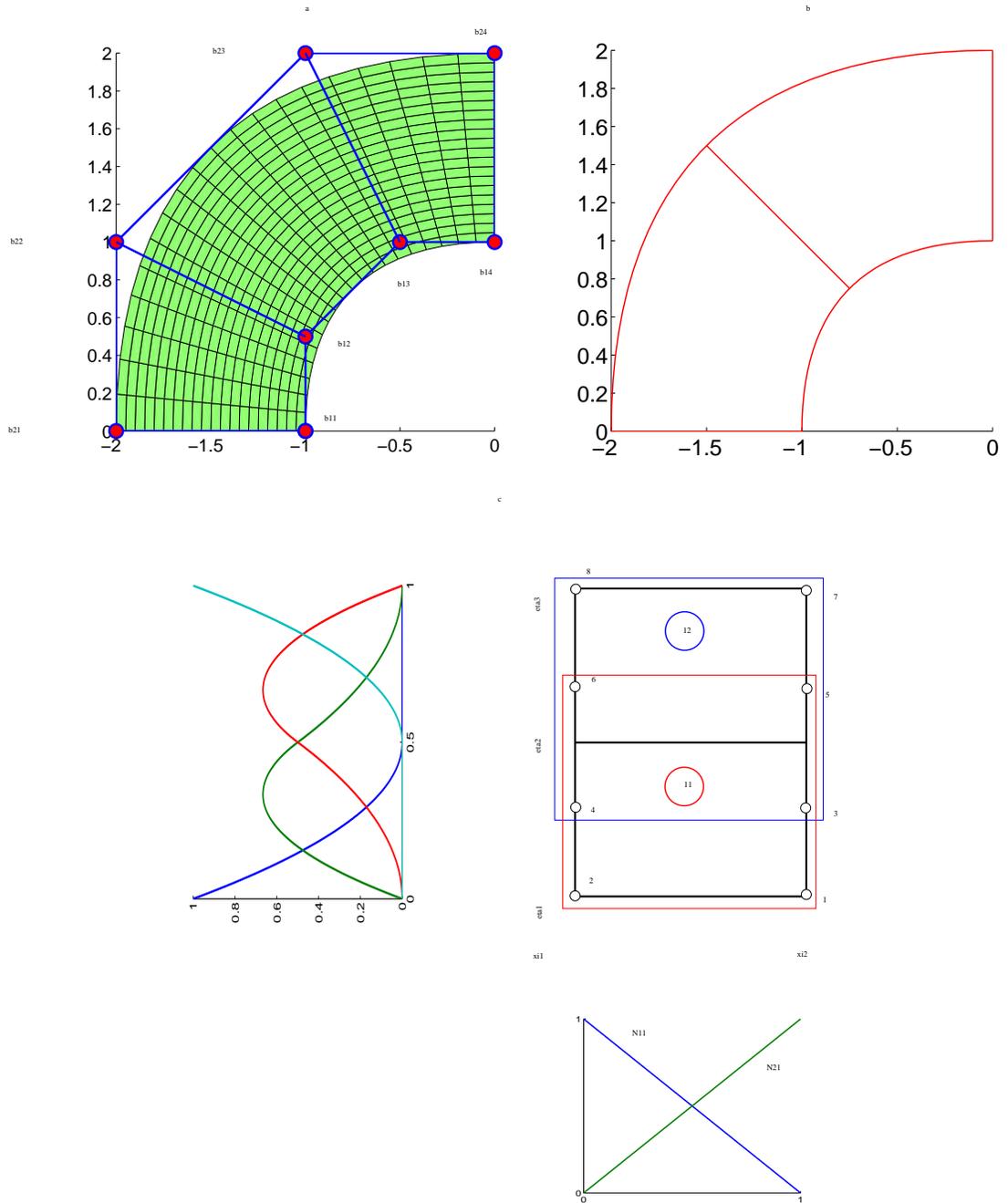}
   \caption{Two dimensional isogeometric analysis example: (a) exact geometry,
   (b) mesh in physical space and (c) mesh in the parametric space. 
   The two rectangles are used to illustrate the control points belonging to each
   element.
   } 
  \label{fig:assembly2d} 
\end{figure}

We now consider a concrete two dimensional problem 
shown in Fig.~\ref{fig:assembly2d}. In this case, the knot vectors are
$\Xi^1=[0, 0, 1, 1]$ and $\Xi^2=[0, 0, 0, 0.5, 1, 1, 1]$, respectively. 
The orders of the basis functions are $p=1$ and $q=2$. There two control points ($n=2$) along the $\xi$ 
direction and four control points ($m=4$) along the $\eta$ direction that results in 8 control points
($n_{np}=n\times m=8$).

The domain and non-zero basis functions for Element 1 are given by 

\begin{equation}
  \begin{array}{ccc}
	  \text{direction}       & \text{knot interval} & \text{non zero basis} \\	  
     \textrm{$\xi$ }     & [\xi_2,\xi_3] & N_1,N_2\\
     \textrm{$\eta$ }    & [\eta_3,\eta_4] & M_1,M_2,M_3\\
  \end{array}
  \label{eq:dep1}
\end{equation}

\noindent Hence there are six non-zero basis functions on element $e=1$
which can be assembled into a vector $\vm{R}$ as follows

\begin{equation}
  \vm{R}_1 = [N_1M_1, N_2M_1, N_1M_2, N_2M_2, N_1M_3, N_2M_3].
  \label{eq:d}
\end{equation}

\noindent These six basis functions are associated with six global basis indices given by
(in Matlab\textsuperscript{\textregistered{}} notation)

\begin{equation}
	\text{IEN}(:,1) = [1,2,3,4,5,6].
  \label{eq:hong}
\end{equation}
Similarly, for element 2, the shape function vector is given by

\begin{equation}
  \vm{R}_2 = [N_1M_2, N_2M_2, N_1M_3, N_2M_3, N_1M_4, N_2M_4],
  \label{eq:du}
\end{equation}

\noindent with the associated global indices

\begin{equation}
  \text{IEN}(:,2) = [3,4,5,6,7,8].
  \label{eq:hong1}
\end{equation}
Control points are stored in a two dimensional matrix of dimensions $n_{np}\times2$.
The connectivity data is stored in a two dimensional matrix of dimensions
$(p+1)*(q+1) \times n_{el}$ where $n_{el}$ denotes the number of elements. 
For the example under consideration, these two matrices are given by

\begin{equation}
	\text{controlPts}=
	\begin{bmatrix}
	\vm{B}_{11}&\vm{B}_{21}&\vm{B}_{12}&\vm{B}_{22}&
	                   \vm{B}_{13}&\vm{B}_{23}&\vm{B}_{14}&\vm{B}_{24}
		   \end{bmatrix}\trans, \quad
        \text{IEN}=\begin{bmatrix}
	1 &2& 3& 4& 5& 6 \\
	3 &4& 5& 6& 7& 8
\end{bmatrix}\trans. 
\end{equation}

\noindent The knot intervals along the $\xi$ and $\eta$ directions are stored in
the following matrices 

\begin{equation}
\text{elRangeU}=\begin{bmatrix}0&1\end{bmatrix},\quad		   
\text{elRangeV}=\begin{bmatrix}0&0.5\\0.5&1\end{bmatrix},		   
\end{equation}

\noindent where the number of rows is equal to the number of elements in each direction
($n_s^1-1,n_s^2-1$ respectively).

With this vector of basis functions $\vm{R}$, the control points associated to this element, 
we can compute the derivatives of the basis functions. The derivatives of the basis
functions with respect to $x$ are stored in the following vector

\begin{equation}
	\vm{R}_{,x} = \begin{bmatrix}
		R_{1,x} & R_{2,x} & R_{3,x} & R_{4,x} & R_{5,x} & R_{6,x}
	\end{bmatrix}\trans,
	\label{eq:derivatives-vector}
\end{equation}

\noindent Similarly, we have $\vm{R}_{,y}$ for the derivatives of the basis
functions with respect to $y$.
Having these basis function derivatives, we  are now ready to define the 
$\vm{B}$ matrix for any element $e$

\begin{equation}
  \vm{B}_e = \left[\begin{array}{ccccc}
	  \vm{R}[1]_{,x} & 0 & \vm{R}[2]_{,x} & 0 & \cdots \\
    0 & \vm{R}[1]_{,y} & 0 & \vm{R}[2]_{,y} & \cdots \\
   \vm{R}[1]_{,y} & \vm{R}[1]_{,x} & \vm{R}[2]_{,y} & \vm{R}[2]_{,x} & \cdots
  \end{array}\right],
  \label{eq:Bmatrix}
\end{equation}

\noindent where the control point displacement vector is stored 
in the following order
$\vm{u}=[u_{x1},u_{y1},u_{x2},u_{y2},\ldots,u_{xn_{np}},u_{yn_{np}}]\trans$. 
The element stiffness matrix is then given by

\begin{equation}
  \vm{K}_e = \int_{\Omega_e} \vm{B}_e \trans \vm{D} \vm{B}_e \mathrm{d} \Omega_e,
  \label{eq:stiffness}
\end{equation}

\noindent which is then assembled to the global stiffness matrix using the
element connectivity matrix and the fact that a control point $I$ corresponds
to positions $2*I-1$ and $2*I$ in the global displacement vector $\vm{u}$.

\begin{figure}[htbp]
  \centering 
  \psfrag{1}{1} \psfrag{2}{2} 
  \psfrag{4}{4} \psfrag{8}{8}\psfrag{12}{12} 
  \psfrag{3}{3} \psfrag{4}{4}
  \psfrag{5}{5} \psfrag{6}{6} \psfrag{7}{7} \psfrag{8}{8}
  \psfrag{9}{9} \psfrag{10}{10} \psfrag{11}{11} \psfrag{12}{12}
  \includegraphics[width=0.35\textwidth]{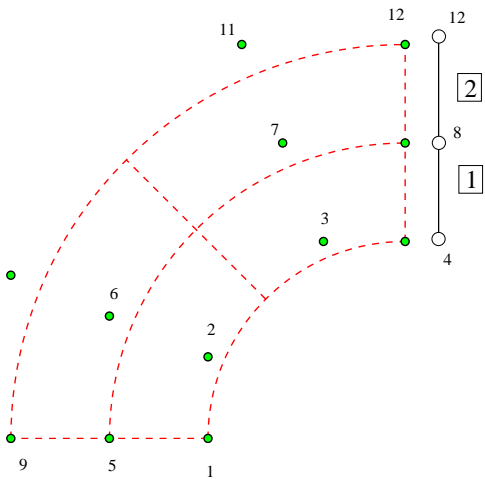}
  \caption{Boundary mesh for external force computation. A linear basis is
  used in the $\eta$ direction. Assume a traction is applied on the edge
  containing  nodes 4, 8 and 12. The boundary mesh is composed of  two linear isogeometric elements 
  \SquareBox{1} and \SquareBox{2}.}
   \label{fig:fext}
\end{figure}

For implementation convenience,  Box \ref{box:elasticity2d}
gives a procedure of an isogeometric analysis of 2D elasticity problems. Note that this
aims for a Matlab implementation and the whole elementary stiffness matrix is computed
in one go.
In order to compute the external force vector, it is convenient to define a
boundary mesh as shown in Fig. \ref{fig:fext}. 
The computation of the external force vector $\int_\Gamma
\vm{R}\trans \bar{\vm{t}}\di\Gamma$ then follows the one dimensional assembly procedure
given in Section \ref{1dAssembly}.

\begin{Fbox}
  \caption{Procedure for isogeometric analysis of 2D elasticity problems.}
  \begin{enumerate}
    \item Loop over elements, $e=1,\ldots,nel$
      \begin{enumerate}
        \item Determine NURBS coordinates $[\xi_i,\xi_{i+1}]\times[\eta_j,\eta_{j+1}]$ using elRangeU and elRangeV
	\item Get connectivity array, $sctr=\text{IEN}(:,e)$
	\item Define $sctrB(1,1:2:2*nn) = 2*sctr-1;sctrB(1,2:2:2*nn) = 2*sctr$ for assembly$^*$
	\item $\vm{K}_e=\vm{0}$
	\item Loop over Gauss points, $\{\tilde{\bsym{\xi}}_j, \tilde{w}_j\} \quad j=1,2,\ldots,n_{gp}$$^{**}$
           \begin{enumerate}
	     \item Compute $\bsym{\xi}$ corresponding to $\tilde{\bsym{\xi}}_j$
		     (Eq. \ref{eq:phi_mapping})
	     \item Compute $|{J}_{\tilde{\xi}}|$ (Eq. \ref{eq:jacob_determ_parent})
	     \item Compute derivatives of shape functions  $\vm{R}_{,\xi}$ and $\vm{R}_{,\eta}$ at  $\bsym{\xi}$ 
	     \item Compute $\vm{J}_\xi$ using $\text{controlPts}(sctr,:)$,
		     $\vm{R}_{,\xi}$ and $\vm{R}_{,\eta}$ (Eq. \ref{eq:jacobian-map})
	     \item Compute Jacobian inverse $\vm{J}^{-1}_\xi$ and determinant
		     $|{J}_\xi|$
	     \item Compute derivatives of shape functions
		     $\vm{R}_{,\vm{x}}=[\vm{R}_{,\xi}\;
		     \vm{R}_{,\eta}]\vm{J}^{-1}_\xi$ (Eq. \ref{eq:sd})
	     \item Use $\vm{R}_{,\vm{x}}$ to build the strain-displacement matrix $\vm{B}$ (Eq. \ref{eq:Bmatrix}
                   )
	     \item Compute $\vm{K}_e = \vm{K}_e + \tilde{w}_j |{J}_{\tilde{\xi}}|
		     |{J}_{{\xi}}| \vm{B}\trans\vm{D}\vm{B}$
           \end{enumerate}
	 \item End loop over Gauss points
	 \item Assemble $\vm{K}_e$: $\vm{K}(sctrB,sctrB) = \vm{K}(sctrB,sctrB) + \vm{K}_e$
      \end{enumerate}
    \item End loop over elements
  \end{enumerate}
  $^*$ $nn$ denotes the number of control points per element \ie nn=length(sctr).\\
  $^{**}$ $\tilde{w}_j$ denotes a Gauss point weight and $n_{gp}$ is the total number of Gauss points.
  \label{box:elasticity2d}
\end{Fbox}

%

\subsection{Boundary condition enforcement}\label{sec:BCs}

\begin{figure}[htbp]
	\centering
  \psfrag{N11}{$N_{1}(\xi)$}
  \psfrag{N21}{$N_{2}(\xi)$}
  \psfrag{N12}{$M_{1}(\eta)$} \psfrag{N22}{$M_{2}(\eta)$}
  \psfrag{N32}{$M_{3}(\eta)$} \psfrag{N42}{$M_{4}(\eta)$}
  \psfrag{xi1}{$\xi_{1,2}=0$}
  \psfrag{xi2}{$\xi_{3,4}=1$}
  \psfrag{11}[c]{1} \psfrag{12}[c]{2}
  \psfrag{1}{1} \psfrag{2}{2} \psfrag{3}{3} \psfrag{4}{4}
  \psfrag{5}{5} \psfrag{6}{6} \psfrag{7}{7} \psfrag{8}{8}
  \psfrag{a}{A}\psfrag{b}{B}\psfrag{c}{C}\psfrag{d}{D}
	\includegraphics[width=0.4\textwidth]{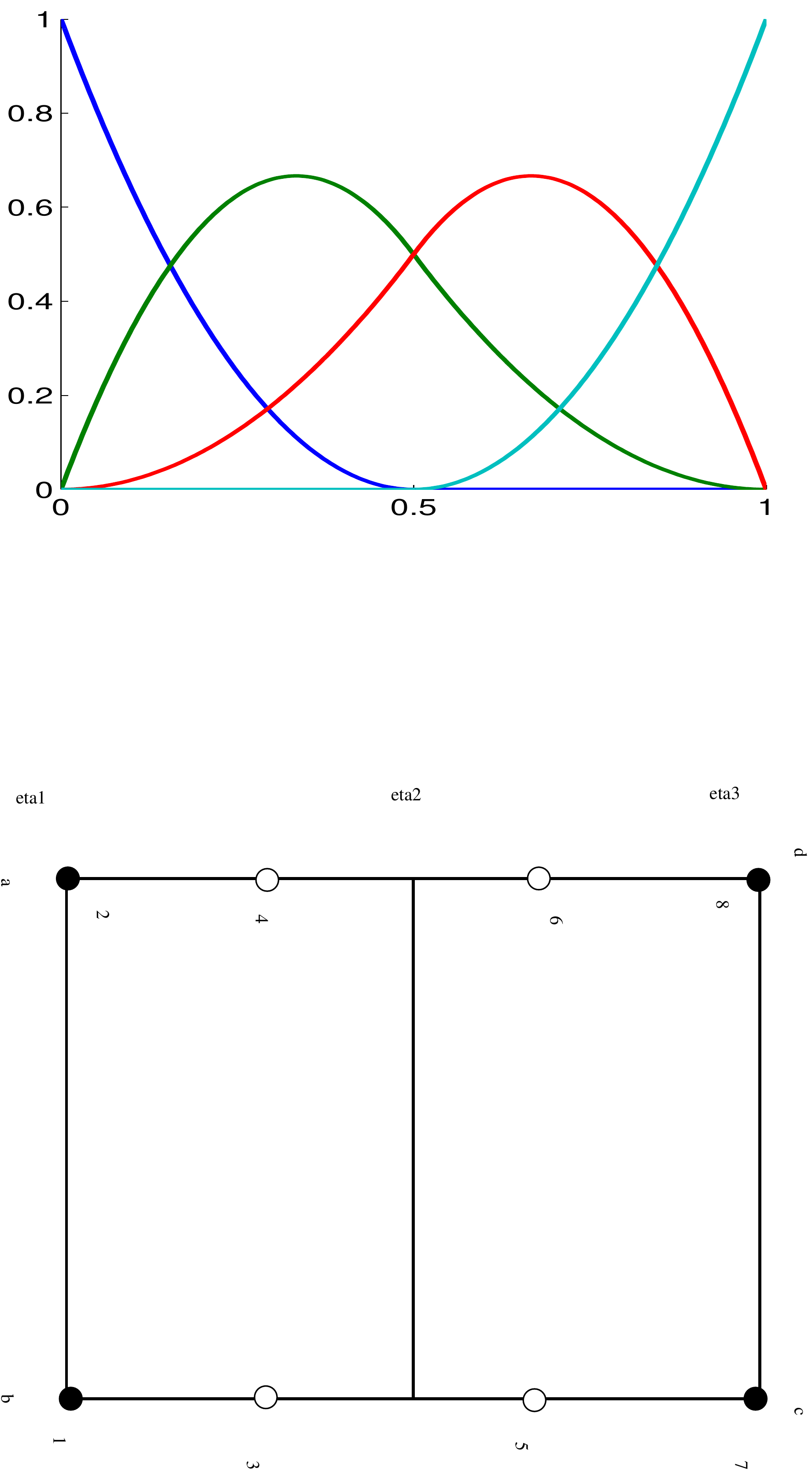}
	\caption{Imposing Dirichlet BCs:
	black points denote corner control points where the NURBS basis
	satisfies the Kronecker delta property.}
	\label{fig:boundary-condition}
\end{figure}

Fig. \ref{fig:boundary-condition} illustrates two kinds of Dirichlet
boundary conditions: on edge $AD$, $u_x=0$ and on edge $BC$, $u_y=\bar{u}$.
The former is a homogeneous Dirichlet boundary condition (BC) while
the latter is referred to as uniform inhomogeneous Dirichlet BCs. 
Homogeneous Dirichlet BCs can be enforced by setting the corresponding 
control variables as zeros (in this example, setting $u_{xI}=0,I=2,4,6,8$).
For edge $BC$, the inhomogeneous Dirichlet BCs
can also be satisfied by setting $u_{yI}=\bar{u},I=1,3,5,7$. This is due to
the partition of unity property of the NURBS basis \ie $u_y^{BC}=M_1(\eta)u_{y1}+
M_2(\eta)u_{y3}+M_3(\eta)u_{y5}+M_4(\eta)u_{y7}=(M_1(\eta)+M_2(\eta)+M_3(\eta)+M_4(\eta))\bar{u}=\bar{u}$.
Note that a boundary of a NURBS surface is a NURBS curve due to the use of open knots.
Inhomogeneous Dirichlet BCs applied on corner control points (black
points in Fig. \ref{fig:boundary-condition}) are enforced by simply setting
the corner control variables equal to the prescribed values since
the NURBS shape functions at these points satisfy the Kronecker delta property (assuming the use of open
knot vectors). This is called direct imposition 
of Dirichlet BCs.

For cases other than the ones previously discussed such as a prescribed
displacement imposed at interior control point 3 or a non-uniform Dirichlet
BC--Dirichlet BCs that vary from point to point
applied on edge $BC$, special treatment of
Dirichlet BCs have to be employed as is the case for meshless methods.
Techniques available include the Lagrange multiplier method, the penalty method, the 
augmented Lagrangian method and we refer to \cite{nguyen_meshless_2008} for an
overview of these techniques in the context of meshless methods. In 
\cite{wang_improved_2010} a transformation method was proposed to impose
inhomogeneous Dirichlet BCs in IGAFEM. However this method requires modifications to the
stiffness matrix which breaks the usual structure of a FE code. The authors in \cite{de_luycker_xfem_2011}
presented a weak enforcement of general inhomogeneous Dirichlet BCs using a
least squares minimization and the implementation is described in what follows.
The same procedure was used in imposing BCs in meshless methods 
when coupling a fluid to a solid domain through a master-slave concept \cite{Rabczuk_Immersed}.
Imposing Dirichlet boundary conditions with Nitsche's method was presented in 
\cite{embar_imposing_2010} for spline-based finite elements.

The basic idea of the least squares method is to find the parameters of
the boundary control points that minimize the following quantity

\begin{equation}
\begin{split}
	J &= \frac{1}{2}\sum_{C} \norm{ \vm{u}(\vm{x}_C) -
	\bar{\vm{u}}(\vm{x}_C)}^2\\
	  &= \frac{1}{2}\sum_{C} \norm{ \sum_A R_A(\vm{x}_C) \vm{q}_A -
	  \bar{\vm{u}}(\vm{x}_C)}^2
\end{split},
	\label{eq:delta-uh-uB}
\end{equation}

\noindent where $\vm{x}_C$ denotes a set of collocation points distributed
on the essential boundary $\Gamma_u$, $\vm{q}_A$ are the parameters of the control points
defining $\Gamma_u$, $\bar{\vm{u}}(\vm{x})$ represents the prescribed displacements
and $R_A$ represents the NURBS basis functions that are
non-zero at $\vm{x}_C$ which are univariate NURBS functions $R_A(\xi)$ for the boundary of a NURBS surface
is a NURBS curve. For the sake of clarity, let us consider the case
where there is only one collocation point and a quadratic basis (thus there
are 3 non-zero $R_A$ at $\vm{x}_C$). So, we have

\begin{equation}
	J = \frac{1}{2}\norm{ R_1(\vm{x}_C)\vm{q}_1 + R_2(\vm{x}_C)\vm{q}_2 +
	R_3(\vm{x}_C)\vm{q}_3 - \bar{\vm{u}}(\vm{x}_C) }^2.
\end{equation}

\noindent The partial derivatives of $J$ with respect to $\vm{q}_i$ are given
by

\begin{equation}
\begin{split}
	\pderiv{J}{\vm{q}_1} &= [R_1(\vm{x}_C)\vm{q}_1 + R_2(\vm{x}_C)\vm{q}_2 +
	R_3(\vm{x}_C)\vm{q}_3 - \bar{\vm{u}}(\vm{x}_C)]R_1(\vm{x}_C)\\
	\pderiv{J}{\vm{q}_2} &= [R_1(\vm{x}_C)\vm{q}_1 + R_2(\vm{x}_C)\vm{q}_2 +
	R_3(\vm{x}_C)\vm{q}_3 - \bar{\vm{u}}(\vm{x}_C)]R_2(\vm{x}_C)\\
	\pderiv{J}{\vm{q}_3} &= [R_1(\vm{x}_C)\vm{q}_1 + R_2(\vm{x}_C)\vm{q}_2 +
	R_3(\vm{x}_C)\vm{q}_3 - \bar{\vm{u}}(\vm{x}_C)]R_3(\vm{x}_C).
\end{split}
\end{equation}

\noindent The condition $\pderiv{J}{\vm{q}}=0$ thus gives the following
linear system 	

\begin{equation}
	\begin{bmatrix}
           R_1R_1 & R_2R_1& R_3R_1\\
           R_1R_2 & R_2R_2& R_3R_2\\
           R_1R_3 & R_2R_3& R_3R_3\\
   \end{bmatrix}_{\vm{x}_C}
	\begin{bmatrix}
		q_1^x & q_1^y \\
		q_2^x & q_2^y \\
		q_3^x & q_3^y 
	\end{bmatrix}=
	\begin{bmatrix}
		\bar{u}_x (\vm{x}_C) R_1 (\vm{x}_C) & \bar{u}_y(\vm{x}_C)
		R_1(\vm{x}_C) \\
		\bar{u}_x (\vm{x}_C) R_2 (\vm{x}_C) & \bar{u}_y(\vm{x}_C)
		R_2(\vm{x}_C) \\ 
		\bar{u}_x (\vm{x}_C) R_3 (\vm{x}_C) & \bar{u}_y(\vm{x}_C)
		R_3(\vm{x}_C) 
	\end{bmatrix}.
	\label{eq:least-square-1point}
\end{equation}

\noindent By collecting all the NURBS basis at point $\vm{x}_C$ in a column
vector $\vm{N}(\vm{x}_C)$, the control points displacements in the $x$ and $y$
directions in $\vm{q}_x$ and $\vm{q}_y$, respectively, Eq.
(\ref{eq:least-square-1point}) can be written in a compact form as

\begin{equation}
\begin{split}
	[ \vm{N}(\vm{x}_C) \vm{N}\trans(\vm{x}_C)] \vm{q}_x &=
	\bar{u}_x(\vm{x}_C) \vm{N}(\vm{x}_C) \\
	[ \vm{N}(\vm{x}_C) \vm{N}\trans(\vm{x}_C) ] \vm{q}_y &=
	\bar{u}_y(\vm{x}_C) \vm{N}(\vm{x}_C).
\end{split}
\end{equation}

\noindent Repeating the same analysis for other collocation points $\vm{x}_C$
on the Dirichlet boundary, one obtains the linear system $\vm{A}\vm{q}=\vm{b}$ 
with two different $\vm{b}$ (one for the $x$ component and the other for the $y$
component). The dimension of 
$\vm{A}$ is $n_D \times n_D$ where $n_D$ denotes the number of control points
defining the Dirichlet boundary. Having these boundary control point
displacements, the enforcement of Dirichlet BCs (when solving
$\vm{K}\vm{u}=\vm{f}$) are then treated as in
standard FEM.

We note that this procedure involves only control points that define
the essential boundary. This is in sharp contrast to meshless shape functions
such as the MLS used in the Element Free Galerkin method \cite{bel:efgm1} 
where the displacements at a point on the
essential boundary depend not only on  the nodes on that boundary but also the
neighbouring interior nodes. It is said that NURBS therefore satisfy the so-called weak Kronecker delta
property as local Maxent interpolants approximations \cite{NME:NME1534}.

\begin{figure}[htbp]
  \centering 
  \psfrag{x}{$x$}\psfrag{y}{$y$}
  \psfrag{1}{1}\psfrag{2}{2}\psfrag{3}{3}\psfrag{4}{4}\psfrag{5}{5}\psfrag{6}{6}\psfrag{7}{7}
  \psfrag{control}{control points}\psfrag{col}{collocation points $\vm{x}_C$}
  \includegraphics[width=0.35\textwidth]{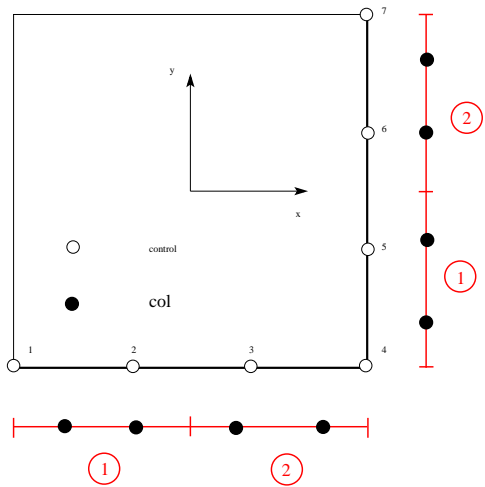}
  \caption{Illustration of the implementation of the least squares method: 
  a quadratic NURBS surface with knot vectors
  $\Xi^1=\{0,0,0,0.5,1,1,1\}$ and $\Xi^2=\{0,0,0,0.5,1,1,1\}$. 
  Essential BCs are imposed on the bottom and right edges. A one dimensional
  mesh for this boundary is created. The control points defining the
  essential boundary are numbered from one to the total number of boundary control
  points (seven in this example). 
  Note that this numbering is required only for assembling the matrix $\vm{A}$.  } 
  \label{fig:least-square}
\end{figure}

Listing \ref{least-square-list} gives the Matlab\textsuperscript{\textregistered{}} implementation of the
least-squares method. Note that in our implementation, the collocation points are uniformly distributed in
the parameter space. We refer to \cite{de_luycker_xfem_2011} for a discussion on the influence of
the collocation points on the accuracy. We refer to Fig. \ref{fig:least-square} for an
illustration of this method.

\begin{snippet}[caption={Matlab\textsuperscript{\textregistered{}} implementation of the least-squares method},
                   label={least-square-list},framerule=1pt]
    A  = zeros(noDispNodes,noDispNodes);
    bx = zeros(noDispNodes,1);
    by = zeros(noDispNodes,1);
    noxC   = 4; % number of collocation pts/element
    % loop over bottom edge
    for ie=1:noElemsU
        sctr   = bottomEdgeMeshIGA(ie,:); % standard connectivity array (global numbering)
        pts    = controlPts(sctr,:);      % control pts of element ie
        sctrA  = bndElement(ie,:);        % connectivity array for A matrix
        xiE    = elRangeU(ie,:);          % parametric coords of ie
        xiArr  = linspace(xiE(1),xiE(2),noxC); % collocation pts x_C
        for ic=1:noxC
            xi             = xiArr(ic);
            [N dNdxi]      = NURBS1DBasisDers(xi,p,uKnot,weights);
            A(sctrA,sctrA) = A(sctrA,sctrA) + N'*N;
            x              = N *pts; 
            % exact displacements
            [ux,uy]   = exact_Griffith(x,E0,nu0,sigmato,xTip,seg,cracklength);
            bx(sctrA) = bx(sctrA) + ux*N';
            by(sctrA) = by(sctrA) + uy*N';
        end
    end	
    % loop over other edges if neccessary
    % solve the system Aq_x=bx and Aq_y=by
    [LL UU] = lu(A);
    qxTemp  = LL\bx;
    qyTemp  = LL\by;
    qx      = UU\qxTemp; qy      = UU\qyTemp;
    % later, before solving Ku=f, qx and qy will be used  to enforce BCs as in conventional FEM.
\end{snippet}

In our Matlab\textsuperscript{\textregistered{}} code implementation are provided for the 
penalty method, the Lagrange
multiplier method and the least squares method. The implementation of
the two former methods are considered standard and the reader is referred to
\cite{nguyen_meshless_2008} for details.

\section{Extended isogeometric finite element method}\label{sec:xiga}

There are basically two ways in which discontinuities can be modeled in the context of IGA:
PUM based enrichment and knot insertion. For the former, there are works of  \cite{benson_gfem_2010,
de_luycker_xfem_2011,ghorashi_extended_2012,Tambat20121}. For the latter, we refer to 
\cite{verhoosel_isogeometric_2011,nguyen_cohesive_2013}. PUM based methods are general for
they can be applied to any kinds of discontinuity such as weak and strong discontinuities whereas
knot insertion used to produce the $C^{-1}$ continuity is only suitable for cracks. Knot insertion
found great applications for delamination analyses \cite{nguyen_cohesive_2013} where the crack path
is known \textit{a priori}. Note that in \cite{verhoosel_isogeometric_2011}, knot insertion was used
in combination with T-splines to model cracks of which trajectory is not known in advance. However the
implementation is tedious.

In the category of PUM based IGA methods which are sometimes called XIGA (eXtended IGA), the method
in \cite{Tambat20121} is different from existing XIGA formulations  \cite{benson_gfem_2010,
de_luycker_xfem_2011,ghorashi_extended_2012} in a way that discontinuities (holes/inclusions/cracks)
are exactly isogeometrically represented by NURBS. Note that in most of XIGA techniques, discontinuities
are defined implicitly using level set method. This section briefly presents the XIGA formulation of
\cite{de_luycker_xfem_2011,ghorashi_extended_2012} and implementation aspects are deferred to
Section \ref{sec:xfem-implementation}.

\begin{figure}[htbp]
  \centering 
  \psfrag{con}{control points}\psfrag{hea}[c]{Heaviside enriched nodes}
  \psfrag{bra}{crack tip enriched nodes}\psfrag{ver}{vertices of physical mesh}
  \includegraphics[width=0.4\textwidth]{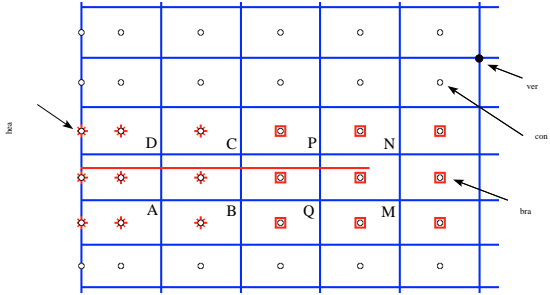}
  \caption{Illustration of enriched node sets $\mathcal{S}^f$ and
  $\mathcal{S}^c$ for a quadratic NURBS mesh. The thick red line denotes the crack.} 
  \label{fig:enriched-nodes}
\end{figure}

The extended finite element method (XFEM) (see \eg \cite{mos_finite_1999},
\cite{fries_extended/generalized_2010} for a recent review and open XFEM library
\cite{bordas_extended_2007}) is a local PUM enrichment method in which internal boundaries
such as holes, inclusions, cracks are modeled independently of the FE discretisation which
allows for crack growth modeling without remeshing. Two dimensional extended
isogeometric finite element formulation (XIGA)  was presented in 
\cite{de_luycker_xfem_2011,ghorashi_extended_2012} in which
the displacement field is enriched for traction-free crack modelling using the following approximation

\begin{equation}
{\rm {\bf u}}^h({\rm {\bf x}})=\sum\limits_{I\in \mathcal{S}}
{R_I} ({\rm {\bf x}}){\rm {\bf u}}_I + \sum\limits_{J\in
\mathcal{S}^{c}} {R_J ({\rm {\bf x}})H ({\rm {\bf x}}){\rm {\bf
a}}_J } +\sum\limits_{K\in \mathcal{S}^{f}} {R_K ({\rm {\bf
x}})\left(\sum\limits_{\alpha =1}^4 B_\alpha{{\rm {\bf b}}_{K}^{\alpha}
}\right)
},
\label{xfemApproximation}
\end{equation}

\noindent where $R_{I,J,K}$ are the NURBS basis functions. In addition to the
standard degrees of freedom (dofs) $\vm{u}_I$, additional dofs $\vm{a}_J$ and
$\vm{b}_K^\alpha$ are introduced. 
The set $\mathcal{S}^{c}$ includes
the control points/nodes whose support is cut by the
crack and the set $\mathcal{S}^{f}$ are control points whose support contains 
the crack tip $\mathbf x_{tip}$, see Fig. \ref{fig:enriched-nodes}.
Note that we use a topological tip enrichment and in the literature
another tip enrichment scheme called geometrical enrichment with a fixed area
(to ensure that the role of enrichment in the approximation space does not
vanish as the mesh is refined) is present 
see \eg \cite{laborde-fix-xfem}.
The Heaviside function $H$ is given by

\begin{equation}
H(\rm {\bf x}) = \left\{\begin{array}{ll}
               +1 & \textrm{if $(\rm{\bf x}-\rm{\bf x}^*)\cdot\rm{\bf n}\geq0 $} \\
               -1 & \textrm{otherwise}
             \end{array}
      \right.,
\end{equation}

\noindent where $\rm{\bf x}^*$ is the projection of point $\rm{\bf
x}$ on the crack and $\vm{n}$ denotes the outward normal vector to the crack. 
And the branch functions, which span the crack tip displacement field, are given by

\begin{equation}
\left[ {B_1 ,B_2 ,B_3 ,B_4 } \right](r,\theta)=\left[
{\sqrt r \sin \frac{\theta }{2},\sqrt r \cos \frac{\theta }{2},\sqrt
r \sin \frac{\theta }{2}\cos{\theta },\sqrt r \cos \frac{\theta
}{2}\cos\theta} \right], \label{branch}
\end{equation}

\noindent where $r$ and $\theta$ are polar coordinates in the local
crack front coordinate system (see \eg \cite{mos_finite_1999} for details). 
It is noted that $B_1$ is discontinuous along the crack face.

Using the standard Galerkin procedure as outlined in Section~\ref{sec:2DAssembly} the
discrete system of equations $\vm{K}\vm{u}=\vm{f}$ are formed by an enlarged $\mathbf B$ matrix given by

\begin{equation}
{\rm{\bf B}} = \left[{
 \begin{array}{c|c}
 {\rm{\bf B}}^{\mathrm{std}} &  {\rm{\bf B}}^{\mathrm{enr}}  \\
 \end{array}}
 \right],
\end{equation}

\noindent where ${\rm{\bf B}}^{\mathrm{std}}$ is the standard strain-displacement matrix \textbf{B} (cf. Eq.~\eqref{eq:strain-disp}) and ${\rm{\bf B}}^\text{enr}$ is the enriched  $\textbf{B}$
matrix of which components are given by

\begin{equation}
 {\rm{\bf B}}_I^{\mathrm{enr}} = \left[{
 \begin{array}{cc}
(R_{I})_{,x}\Psi_{I}+ R_{I}(\Psi_{I})_{,x} & 0 \\
0 & (R_{I})_{,y}\Psi_{I}+ R_{I}(\Psi_{I})_{,y} \\
(R_{I})_{,y}\Psi_{I}+ R_{I}(\Psi_{I})_{,y} &
(R_{I})_{,x}\Psi_{I}+ R_{I}(\Psi_{I})_{,x}
 \end{array}}
 \right],
\end{equation}
$\Psi_{I}$ may represent either the Heaviside
function $H$ or the branch functions $B_{\alpha}$ depending whether control point $I$ is Heaviside or near tip  enriched. The unknown vector $\vm{u}$ contains both 
displacements and enriched dofs. This extended IGAFEM can be implemented within an available
IGAFEM code with little modification following the ideas given in the meshless review
and computer implementation aspects paper \cite{nguyen_meshless_2008}. Some implementation
aspects will be given in Section \ref{sec:xfem-implementation}. In our code enrichment for
holes and material interfaces are also given. They are, however, not covered here because
the implementation follows the ideas given here.

\section{MIGFEM- A Matlab IGA (X)FEM code}\label{sec:migfem}

In this section we describe shortly the open source IGA Matlab (X)FEM program which can be downloaded from
\url{https://sourceforge.net/projects/cmcodes/}.
The code supports one, two and three dimensional linear elasticity problems.
Extended IGA for crack and material interface modelling is also implemented. 
Geometrically nonlinear solid and structural mechanics models are available.
The features of the code include:

\begin{itemize}
\item Global $h$,$p$ and $k$-refinement is provided for one, two and three dimensional meshes. 

\item Extended IGA for 2D/3D stationary traction-free cracks and material interfaces.

\item Visualization of displacements and stresses in Paraview.

\item Inhomogeneous Dirichlet boundary conditions are treated with the penalty
method, the Lagrange multiplier method and the least squares method. 

\item Compatible multi-patch isogeometric formulation for two dimensional problems.

\item Support for T-splines via the B\'{e}zier extraction operators.

\item Structural elements including beams, plates and thin shells. 

\item Implicit Newmark scheme and explicit central difference scheme for time discretisation.
\item Python scripts to extract Rhino3d NURBS surfaces to be used for IGA.
\end{itemize}

\subsection{Data structure}

MIGFEM follows the Matlab FEM code described in \cite{chessa-fem}. 
The main data structures include (1) \textit{element} (store the element connectivity),
(2) \textit{controlPts} (store control point coordinates), (3) \textit{weights} (store the weights) 
(4) \textit{K} (stiffness matrix) and (5)
\textit{f} (external force vector). Contrary to FEM in which the element connectivity and nodal
coordinates are inputs which have been created by a meshing program, in IGA, the input 
consists of CAD data including knots, control points, order of basis functions. Therefore,
one has to construct the \textit{element} matrix based on the knots and the basis orders.\\

\begin{figure}[h!]
  \centering 
  \psfrag{n12}{$N_{1}(\xi)$}
  \psfrag{n22}{$N_{2}(\xi)$}
  \psfrag{n32}{$N_{3}(\xi)$}
  \psfrag{n42}{$N_{4}(\xi)$}
  \psfrag{N12}{$M_{1}(\eta)$} \psfrag{N22}{$M_{2}(\eta)$}
  \psfrag{N32}{$M_{3}(\eta)$} \psfrag{N42}{$M_{4}(\eta)$}
  \psfrag{xi1}{$\xi_{1,2}=0$}
  \psfrag{xi2}{$\xi_{3,4}=1$}
  \psfrag{eta1}{$\eta_{1,2,3}=0$}
  \psfrag{eta2}{$\eta_4=0.5$}
  \psfrag{eta3}{$\eta_{5,6,7}=1$}
  \psfrag{11}[c]{1} \psfrag{12}[c]{2}
  \psfrag{b11}{$\vm{B}_{11}(1)$}
  \psfrag{b12}{$\vm{B}_{12}(3)$}
  \psfrag{b13}{$\vm{B}_{13}(5)$}
  \psfrag{b14}{$\vm{B}_{14}(7)$}
  \psfrag{b21}{$\vm{B}_{21}(2)$}
  \psfrag{b22}{$\vm{B}_{22}(4)$}
  \psfrag{b23}{$\vm{B}_{23}(6)$}
  \psfrag{b24}{$\vm{B}_{24}(8)$}
  \psfrag{1}{1} \psfrag{2}{2} \psfrag{3}{3} \psfrag{4}{4}
  \psfrag{5}{5} \psfrag{6}{6} \psfrag{7}{7} \psfrag{8}{8}\psfrag{9}{9}
  \psfrag{10}{10}\psfrag{11}{11}\psfrag{12}{12}\psfrag{13}{13}\psfrag{14}{14}
  \psfrag{15}{15}\psfrag{16}{16}
  \psfrag{0}{0} \psfrag{0}{0}\psfrag{basis}{Basis functions}
  \psfrag{node}[c]{Node pattern}
  \psfrag{connU}{connU}
  \psfrag{connV}{connV}
  \psfrag{xi}{$\Xi^1=\{0,0,0,0.5,1,1,1\}$}
  \psfrag{et}{$\Xi^2=\{0,0,0,0.5,1,1,1\}$}
  \psfrag{1}{1} \psfrag{2}{2} \psfrag{3}{3} \psfrag{4}{4}
  \includegraphics[width=0.95\textwidth]{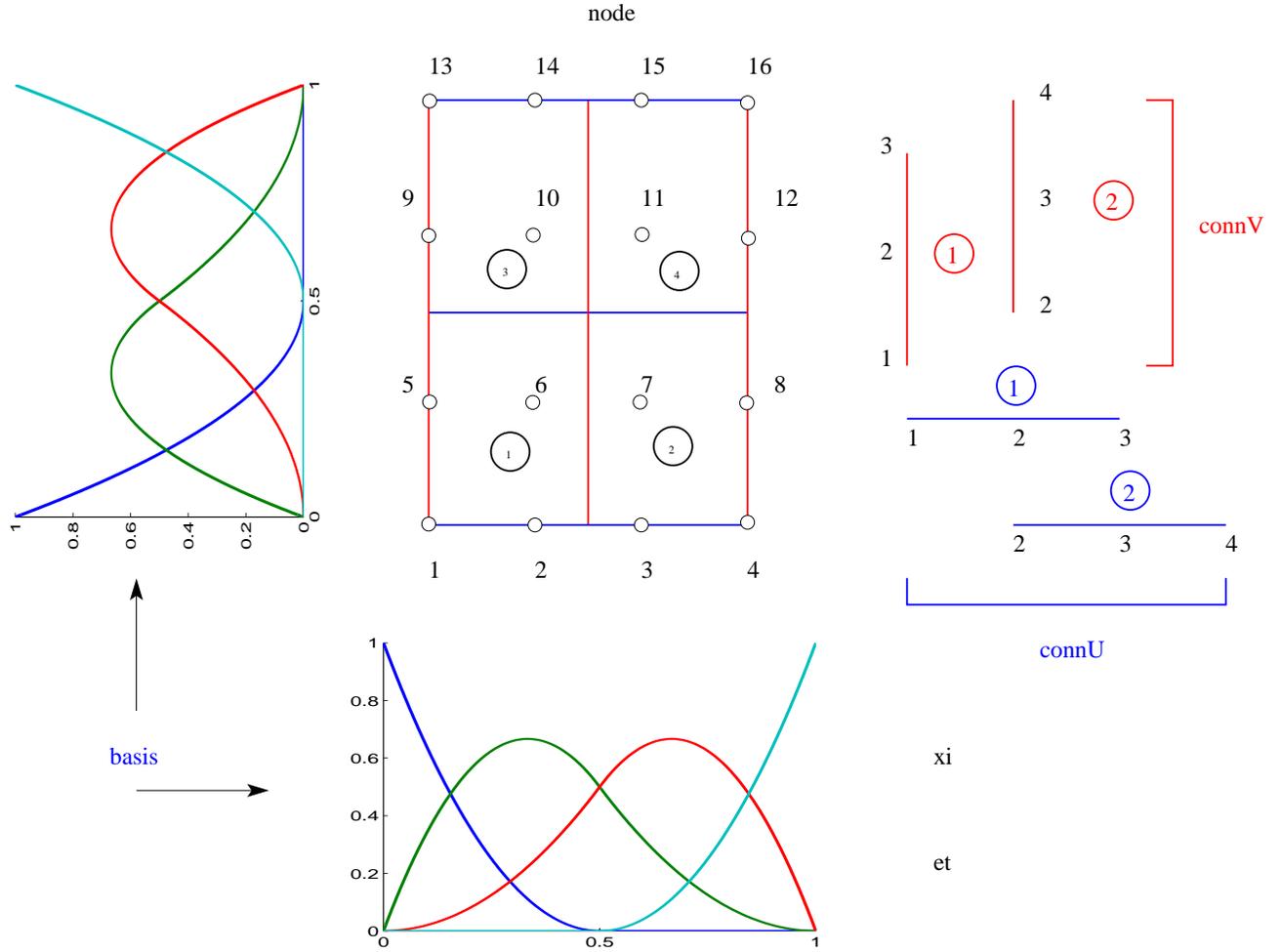}
  \caption{Two dimensional isogeometric analysis: mesh generation for a bi-quadratic
           NURBS surface ($2\times2$ elements). The circles denote the control points.}
  \label{fig:mesh2D} 
\end{figure}

Construction of the \textit{element} matrix is illustrated by a 2D example shown in  Fig. \ref{fig:mesh2D}.
Given the knot vectors $uKnot$ and $vKnot$ together with the orders of the
basis $p$ and $q$, one can compute the number of control points along the $\xi$
and $\zeta$ directions, denoted by $n$ and $m$. Then we define a two
dimensional matrix of dimension $m \times n$ called $node\_pattern$ given by
for the illustrated example shown in Fig. \ref{fig:mesh2D} which 
has 4 control points along $\xi$ and 4 control points along $\eta$ directions 

\begin{equation}
	node\_{pattern} = \begin{bmatrix}
		1 & 2 & 3 & 4\\
		5 & 6 & 7 & 8\\
		9 & 10 & 11 & 12\\
		13 & 14 & 15 & 16\\
	\end{bmatrix},
\label{eq:node-pattern}                
\end{equation}
which is simply an application of Eq.~\eqref{eq:bspline_surface_mapping} for defining global indexes.

The number of elements along the two directions $n_s^1-1,n_s^2-1$ are
$noElemsU=\text{length}(\text{unique}(uKnot))-1$ and
$noElemsV=\text{length}(\text{unique}(uKnot))-1$  (in the Matlab language). 
The connectivity matrix for the $\xi$
direction, denoted by $connU$ which is a $noElemsU \times (p+1)$ matrix and
the connectivity matrix for the $\eta$ direction, denoted by $connV$ which is a $noElemsV \times (q+1)$ matrix
are given by
   
\begin{equation}
	connU = \begin{bmatrix}
		1 & 2 & 3\\
		2 & 3 & 4
	\end{bmatrix}, \quad
	connV = \begin{bmatrix}
		1 & 2 & 3\\
		2 & 3 & 4
	\end{bmatrix},
\end{equation}
for the example under consideration. Matrix $connU$ stores the indices of
the non-zero basis functions at each knot span for the $\xi$ direction.
Having this information and the global indexes of the NURBS basis given in Eq.~\eqref{eq:node-pattern}, 
we are able to define the connectivity matrix
for the whole mesh, called $element$ which is a $noElemsU*noElemsV \times
(p+1)(q+1)$ matrix. For the example being considered, this matrix reads

\begin{equation}
	element = \begin{bmatrix}
		1 & 2 & 3 & 5 & 6 & 7 & 9  & 10 & 11 \\
		2 & 3 & 4 & 6 & 7 & 8 & 10 & 11 & 12 \\
		5 &6 &7 &9 &10 &11 &13 &14 &15\\
		6 &7 &8 &10 &11 &12 &14 &15 &16
	\end{bmatrix}.
\end{equation}
Note that the matrix $element$ is the transpose of the IEN matrix. We decided to use
$element$ instead of IEN to be compatible with the FEM code \cite{chessa-fem} on which MIGEM
is built.

Finally in order to compute the mapping from the parent domain $\tilde{\Omega}$ to
the parametric space, we define the knot intervals in two directions as

\begin{equation}
	rangeU = \begin{bmatrix}
		0 & 0.5\\
		0.5 & 1
	\end{bmatrix}, \quad
	rangeV = \begin{bmatrix}
		0 & 0.5\\
		0.5 & 1
	\end{bmatrix}.
\end{equation}

In order to retrieve the parametric coordinates of a specific element,
the matrix $index$, that is a $noElemsU*noElemsV \times 2$ 
matrix, is used. For a given element $e$, 
its parametric coordinates are determined by \textit{rangeU(index(e,1),:); rangeV(index(e,2),:)}.

The above discussion, together with the illustration given in Fig. \ref{fig:mesh2D} is implemented in Matlab in the file
\textbf{generateIGA2DMesh.m}, located in folder \textbf{meshing}, and shown in Listing \ref{list-mesh2D}.
In that folder, one can find similar M files for generating 1D and 3D meshes.

\begin{snippet}[caption={Mesh generation for two dimensional problems. Input are uKnot, vKnot, and $p,q$.},
	label={material}, label={list-mesh2D},framerule=1pt]
    uniqueUKnots  = unique(uKnot); uniqueVKnots  = unique(vKnot);
    noElemsU      = length(uniqueUKnots)-1;%#of elements xi dir.
    noElemsV      = length(uniqueVKnots)-1;%#of elements eta dir.
    noPtsX        = length(uKnot)-p-1; noPtsY        = length(vKnot)-q-1;
    nodePattern   = zeros(noPtsY,noPtsX);
    count = 1;
    for i=1:noPtsY
        for j=1:noPtsX
            nodePattern(i,j) = count; count            = count + 1;
        end
    end
    % determine our element ranges and the corresponding knot indices along each direction
    [elRangeU,elConnU] = buildConnectivity(p,uKnot,noElemsU);
    [elRangeV,elConnV] = buildConnectivity(q,vKnot,noElemsV);
    noElems = noElemsU * noElemsV;
    element = zeros(noElems,(p+1)*(q+1));
    e = 1;
    for v=1:noElemsV
        vConn = elConnV(v,:);
        for u=1:noElemsU
            c = 1; uConn = elConnU(u,:);
            for i=1:length(vConn)
                for j=1:length(uConn)
                  element(e,c) = nodePattern(vConn(i),uConn(j)); c = c + 1;
                end
            end
            e = e + 1;
        end
    end
    index = zeros(noElems,2);
    count = 1;
    for j=1:size(elRangeV,1)
        for i=1:size(elRangeU,1)
            index(count,1) = i; index(count,2) = j; count = count + 1;
        end
    end
\end{snippet}

\subsection{Shape function routines}

The shape function routines (evaluate NURBS basis functions and derivatives with respect to
parametric coordinates) are implemented using MEX files to improve the performance. They are located in
folder \textbf{C\_files}. Listing \ref{shape-routines} gives some commonly used shape function routines
to compute the NURBS basis function and their first derivatives in 1D, 2D and 3D at a given point. In the 
last line, second derivatives are also  computed.

\begin{snippet}[caption={Shape function routines: uKnot,vKnot,wKnot store $\Xi^1,\Xi^2,\Xi^3$. },
	label={material}, label={shape-routines},framerule=1pt]
    % R, dRdxi: row vectors i.e. dRdxi dimension is (1,3) for p=2    
    [R dRdxi]                = NURBS1DBasisDers(Xi,p,uKnot,weights);
    [R dRdxi dRdeta]         = NURBS2DBasisDers([Xi; Eta],p,q,uKnot,vKnot,weights');
    [R dRdxi dRdeta dRdzeta] = NURBS3DBasisDers([Xi;Eta;Zeta],p,q,r,uKnot,vKnot,wKnot,weights');
    [R dRdxi dRdeta dR2dxi dR2det dR2dxe] = NURBS2DBasis2ndDers([Xi; Eta],p,q,uKnot,vKnot,weights');
\end{snippet}

\begin{snippet1}[caption={Matlab code for 2D spatial derivatives of shape functions. },
                   label={spatial-derivatives},framerule=1pt]
      % input: element e and parameter coordinates of a Gauss point Xi and Eta                   
      sctr   = element(e,:);        %  element scatter vector
      pts    = controlPts(sctr,:);  %  of dimension nn x 2 where nn: number of nodes 
      % derivate of R=Nxi*Neta w.r.t xi and eta
      [dRdxi dRdeta] = NURBS2Dders([Xi;Eta],p,q,uKnot,vKnot,weights'); 
      % Jacobian matrix
      jacob      = pts'*[dRdxi' dRdeta']; 
      invJacob   = inv(jacob);
      % nn x 2 matrix: first column derivatives wrt x ...
      dRdx       = [dRdxi' dRdeta'] * invJacob; 
\end{snippet1}

\subsection{Assembly process}

The assembly of an IGA-FEM code is given in Listing \ref{2d-iga-elasticity} where it
can be seen that the procedure is almost identical to that used in the conventional FEM.
The minor differences lie in (1) the need of the elements in the parameter space (lines 4 to 7)
and (2) the second map (from the parent domain to the parametric domain) in the numerical 
integration of the stiffness matrix (line 35).

\begin{snippet1}[caption={Matlab code for IGA (2D elasticity problems). },
                   label={2d-iga-elasticity},framerule=1pt]
    [W,Q]=quadrature(4, 'GAUSS', 2); % 4x4 point quadrature
    for e=1:noElems % Loop over elements (knot spans)
      idu    = index(e,1);
      idv    = index(e,2);
      xiE    = elRangeU(idu,:); % [xi_i,xi_i+1]
      etaE   = elRangeV(idv,:); % [eta_j,eta_j+1]
      sctr   = element(e,:);    %  element scatter vector
      nn     = length(sctr); nn2 = 2*nn;
      sctrB(1,1:2:2*nn) = 2*sctr-1;
      sctrB(1,2:2:2*nn) = 2*sctr  ;
      B      = zeros(3,2*nn);
      pts    = controlPts(sctr,:); 
      for gp=1:size(W,1)                        % loop over Gauss points 
        pt      = Q(gp,:);                          
        wt      = W(gp);                            
        % compute coords in parameter space
        Xi      = parent2ParametricSpace(xiE,pt(1)); 
        Eta     = parent2ParametricSpace(etaE,pt(2)); 
        J2      = jacobianPaPaMapping(xiE,etaE);
        % derivate of R=Nxi*Neta w.r.t xi and eta
        [dRdxi dRdeta] = NURBS2Dders([Xi;Eta],p,q,uKnot,vKnot,weights'); 
        % Jacobian matrix
        jacob      = pts'*[dRdxi' dRdeta']; 
        J1         = det(jacob); invJacob   = inv(jacob);
        dRdx       = [dRdxi' dRdeta'] * invJacob;
        % B matrix
        B(1,1:2:nn2)  = dRdx(1,:); B(2,2:2:nn2)  = dRdx(2,:);
        B(3,1:2:nn2)  = dRdx(2,:); B(3,2:2:nn2)  = dRdx(1,:);
        % compute elementary stiffness matrix and  assemble it to the global matrix
        K(sctrB,sctrB) = K(sctrB,sctrB) + B'*C*B*J1*J2*wt;
      end
    end
\end{snippet1}

\subsection{Post-processing}\label{sec:visualization}

\begin{figure}[htbp]
	\centering
	\subfloat[NURBS
	mesh]{\includegraphics[width=0.31\textwidth]{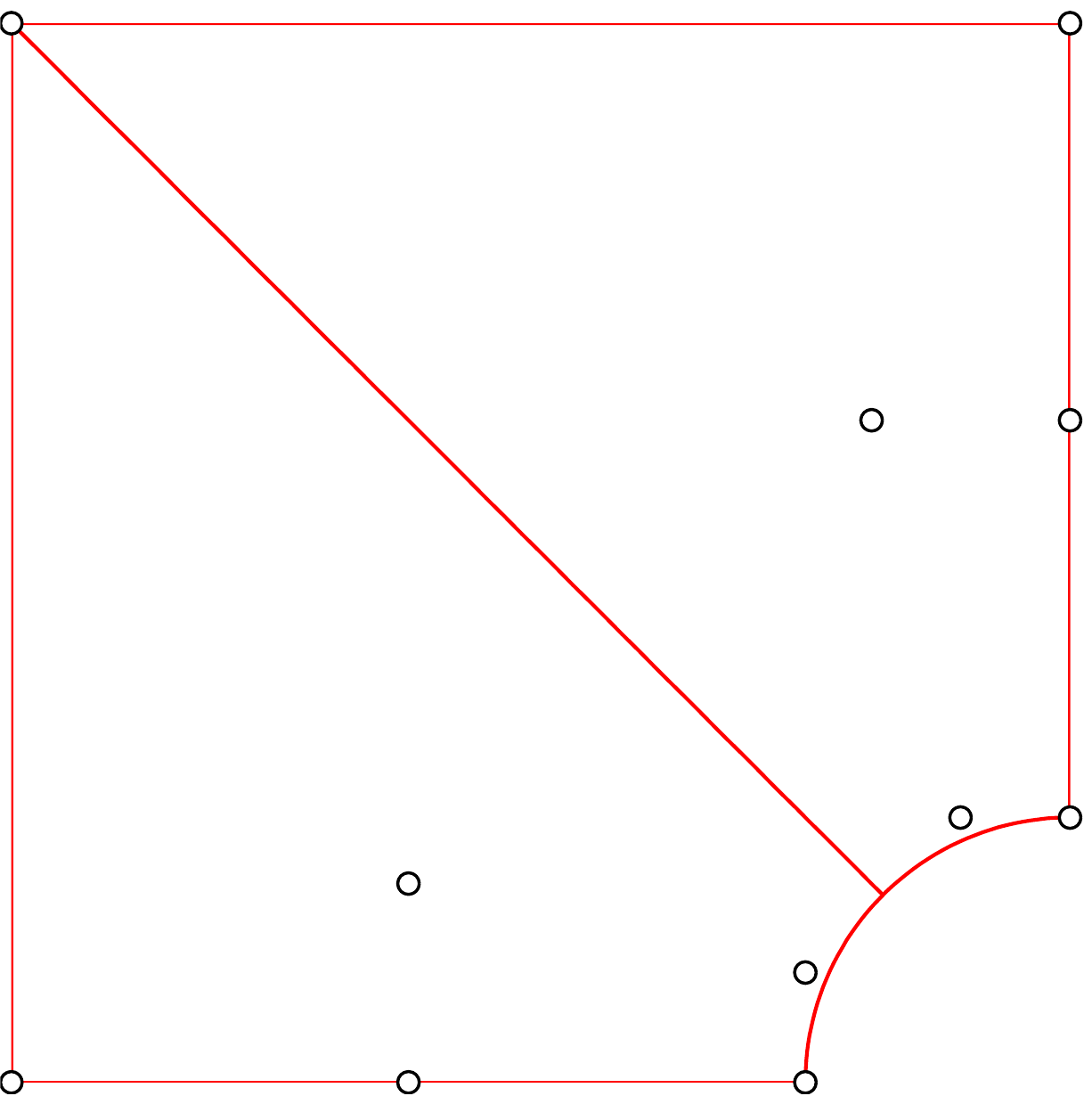}}\;\;\;\;\;
	\subfloat[approximate Q4
	mesh]{\includegraphics[width=0.32\textwidth]{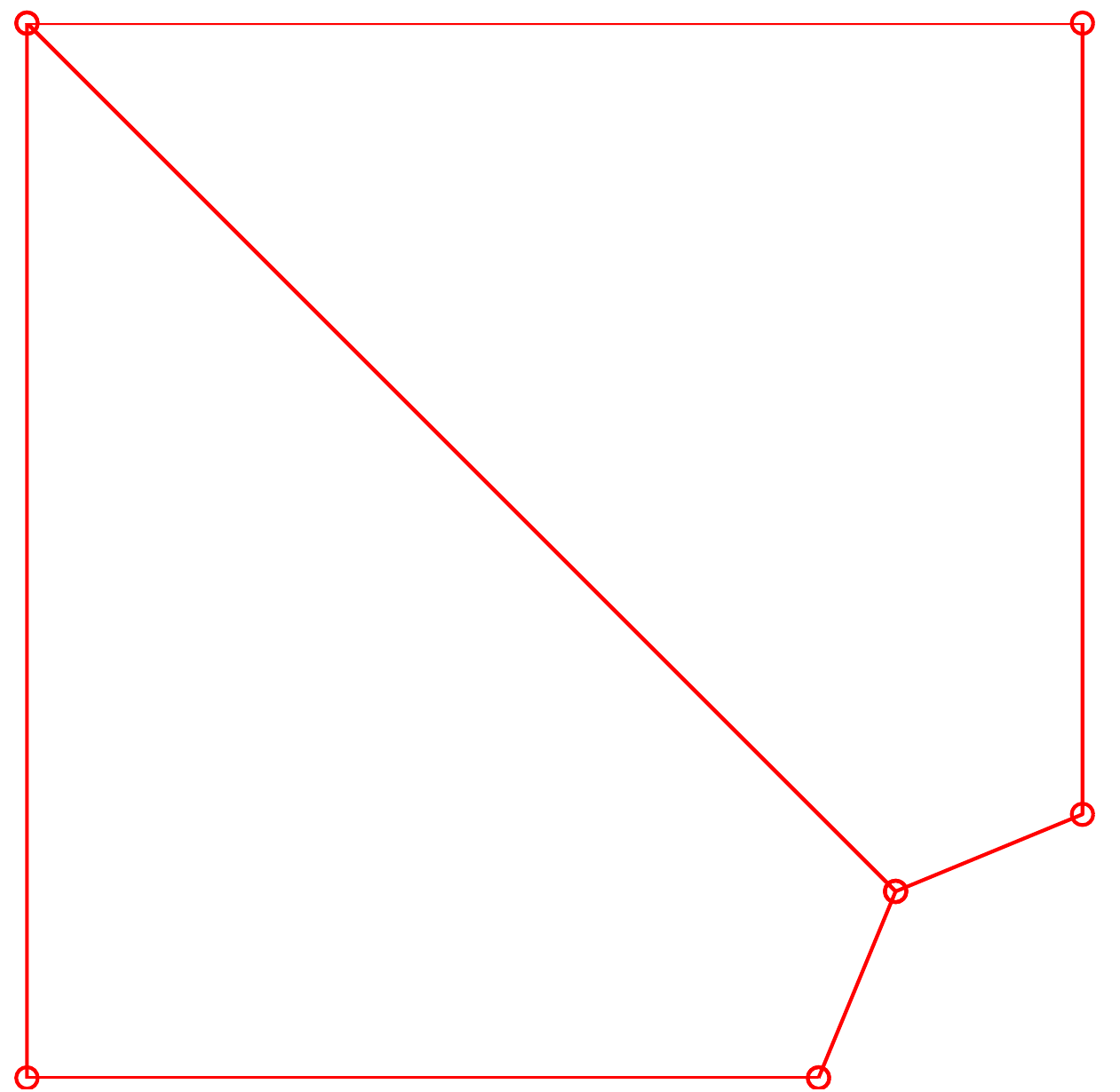}}\\
	\subfloat[stress visualization on a refined mesh]{\includegraphics[width=0.32\textwidth]{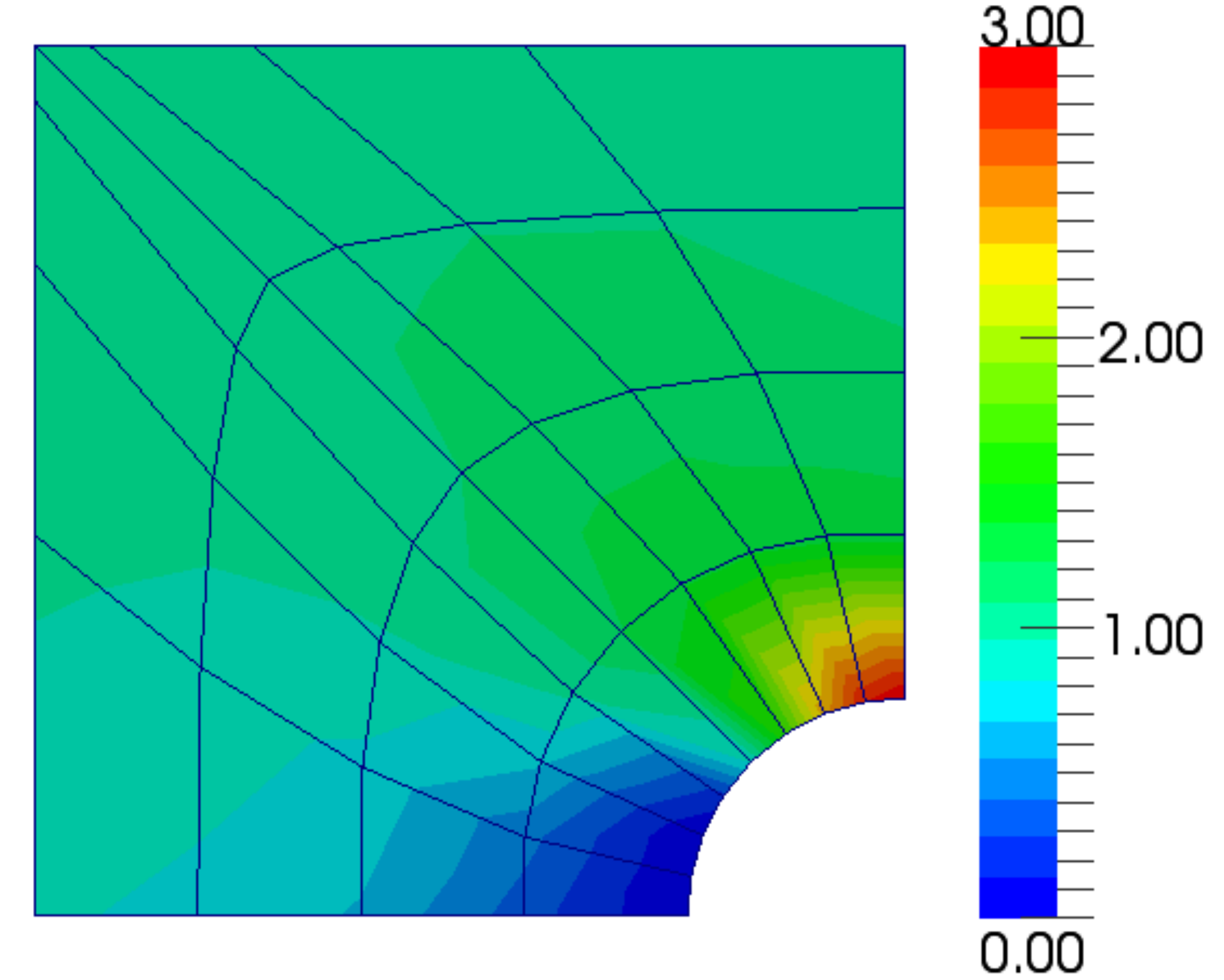}}
	\caption{Exact NURBS mesh (top left) and approximate Q4 mesh (top right) for
	visualization purpose. The nodes in the Q4 mesh are the intersections
	of the $\xi$ and $\eta$ knot lines. The bottom figure shows a contour plot
	of a stress field in Paraview. It should be emphasized that the mesh
	in (b) does not provide a sufficiently smooth contour plot for it to be directly usable. 
        The result given in (c) was
	obtained with a refined NURBS mesh (hence a refined Q4 mesh).}
	\label{fig:visualization1}
\end{figure}

We present here a simple technique to visualize the IGA results that reuse available
visualization techniques for finite elements. To simplify the expos\'{e},
only 2D cases are considered here. In the first step, a mesh
consisting of four-noded quadrilateral (Q4) elements is generated, see Fig.
\ref{fig:visualization1}. We call this mesh the visualization mesh (whose connectivity
matrix is stored in \textit{elementV} and nodal coordinates are stored in \textit{node}) whose nodes
are images of the knots $\xi_i,\eta_j$ in the physical space.
In the second step, quantities of interest \eg
stresses are computed at the nodes of the Q4 mesh. This mesh together with
the nodal values can then be exported to a visualization program such as 
Paraview, see \cite{paraview}, for visualization. It should be emphasized that due to the high
order continuity of the NURBS basis, there is no need to perform nodal averaging
as required in standard $C^0$ finite element analysis to obtain smooth fields. Listing
\ref{2d-iga-postprocessing} gives the Matlab code (we have removed some code due to the similarity
with the assembly code) for building the Q4
visualization mesh, computing the stresses at the nodes of this mesh and exporting
the result to Paraview. The source code can be found in the file
\textbf{plotStress1.m} located in folder \textbf{post-processing}. 
The results can be then visualized directly in Matlab or exported to Paraview, see
Listing \ref{2d-iga-postprocessing1}.

\begin{snippet}[caption={Matlab code for computing stresses and displacements at 
   nodes of the visualization mesh. },
                   label={2d-iga-postprocessing},framerule=1pt]
    buildVisualizationMesh; % build visualization Q4 mesh
    stress = zeros(noElems,4,3);
    disp   = zeros(noElems,4,2);
    for e=1:noElems
        idu    = index(e,1); idv  = index(e,2);
        xiE    = elRangeU(idu,:); etaE = elRangeV(idv,:); 
        sctr   = element(e,:);        % element scatter vector
        sctrB  = [sctr sctr+noCtrPts];%vector scatters B matrix
        uspan = FindSpan(noPtsX-1,p,xiE(1),uKnot);
        vspan = FindSpan(noPtsY-1,q,etaE(1),vKnot);
        % loop over 4 nodes (instead of loop over GPs)
        gp = 1;
        for iv=1:2
          % Q4 elements, the nodes are numbered counter-clockwise
          if (iv==2) xiE = sort(xiE,'descend'); end
          for iu=1:2
            Xi  = xiE(iu); Eta = etaE(iv);
            [N dRdxi dRdeta]= NURBS2DBasisDersSpecial([Xi;Eta],p,q,uKnot,vKnot,weights',[uspan;vspan]);
            % B matrix as usual
            strain          = B*U(sctrB);
            stress(e,gp,:)  = C*strain;
            gp = gp +1;
         end
      end
    end
\end{snippet}

\begin{snippet}[caption={Matlab code for post-processing. },
                   label={2d-iga-postprocessing1},framerule=1pt]
    % plot sigma_xx contour directly in Matlab
    figure
    plot_field(node,elementV,'Q4',stress(:,:,1));
    % export to VTK format to plot in Mayavi or Paraview
    sigmaXX = zeros(size(node,1),1);
    sigmaYY = zeros(size(node,1),1);
    sigmaXY = zeros(size(node,1),1);
    dispX = zeros(size(node,1),1);
    dispY = zeros(size(node,1),1);
    for e=1:size(elementV,1)
        connect = elementV(e,:);
        for in=1:4
            nid = connect(in);
            sigmaXX(nid) = stress(e,in,1);
            sigmaYY(nid) = stress(e,in,2);
            sigmaXY(nid) = stress(e,in,3);
        end
    end
    VTKPostProcess(node,elementV,2,'Quad4','result.vtu',...
                 [sigmaXX sigmaYY sigmaXY],[dispX dispY]);
\end{snippet}

For three-dimensional problems, the same procedure is used where
a mesh of tri-linear brick elements is created and the values of interest
are computed at the nodes of this mesh (see file \textbf{plotStress3d.m}). 
The results are then exported to Paraview under a structured grid format
(*.vts files), see the file \textbf{mshToVTK.m}).\\

\subsection{$h,p,k$-refinement}\label{sec:hpk-refinement}

\begin{figure}[htbp]
  \centering 	
  \subfloat[Initial NURBS mesh]{\includegraphics[width=0.31\textwidth]{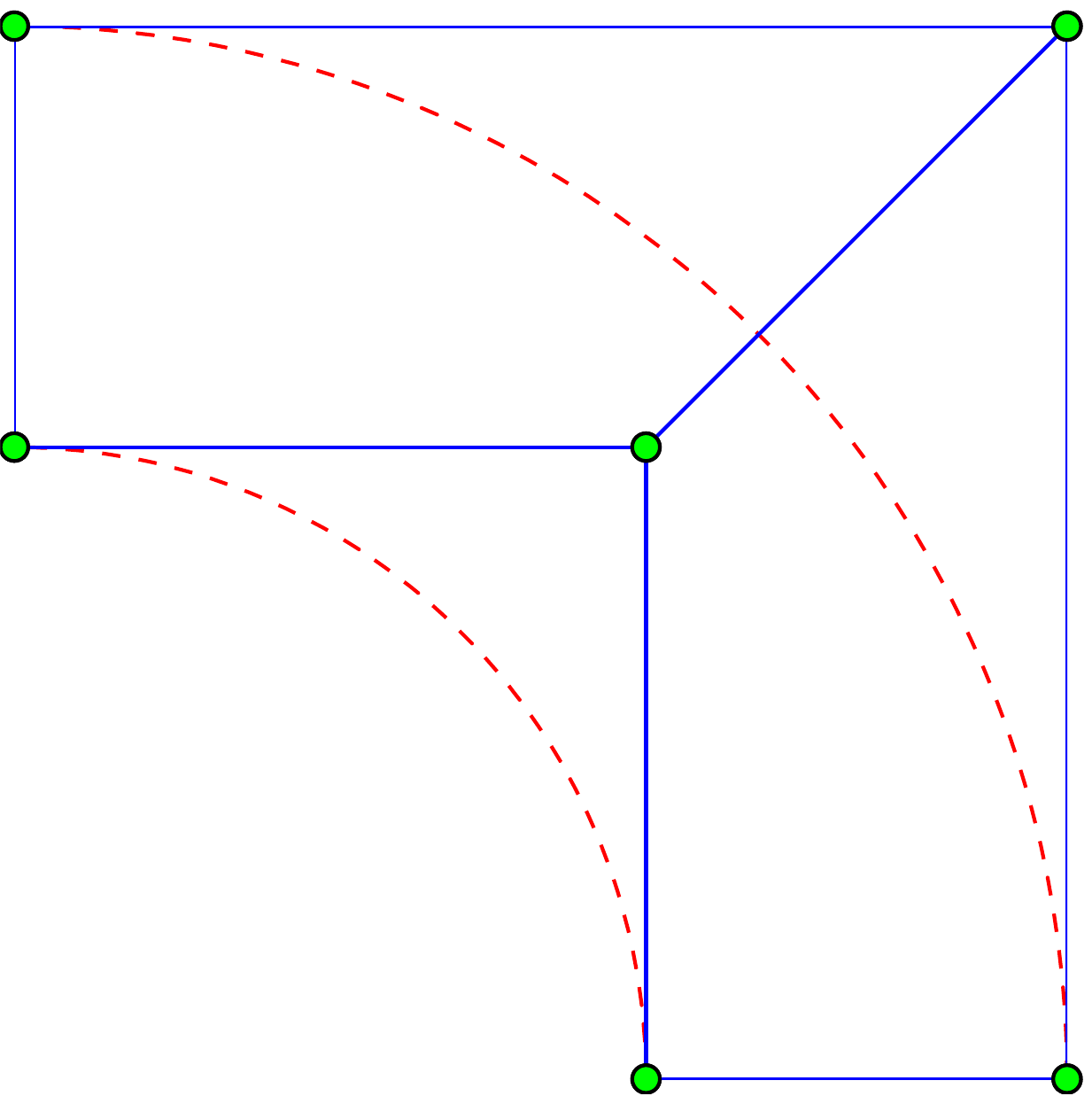}}\;\;\;\;\;
  \subfloat[After h-refinement]{\includegraphics[width=0.31\textwidth]{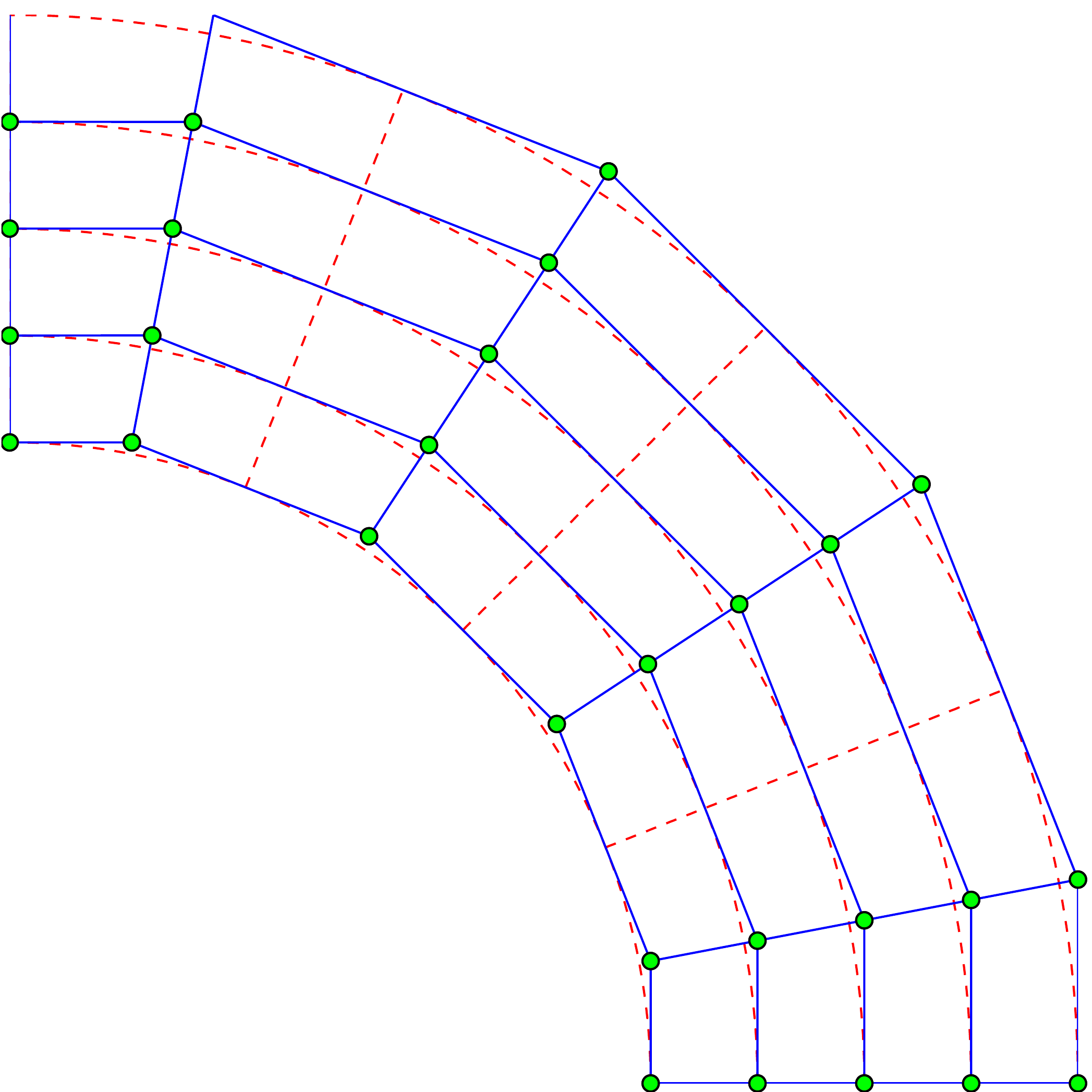}}\;\;\;\;\;
  \subfloat[After k-refinement]{\includegraphics[width=0.31\textwidth]{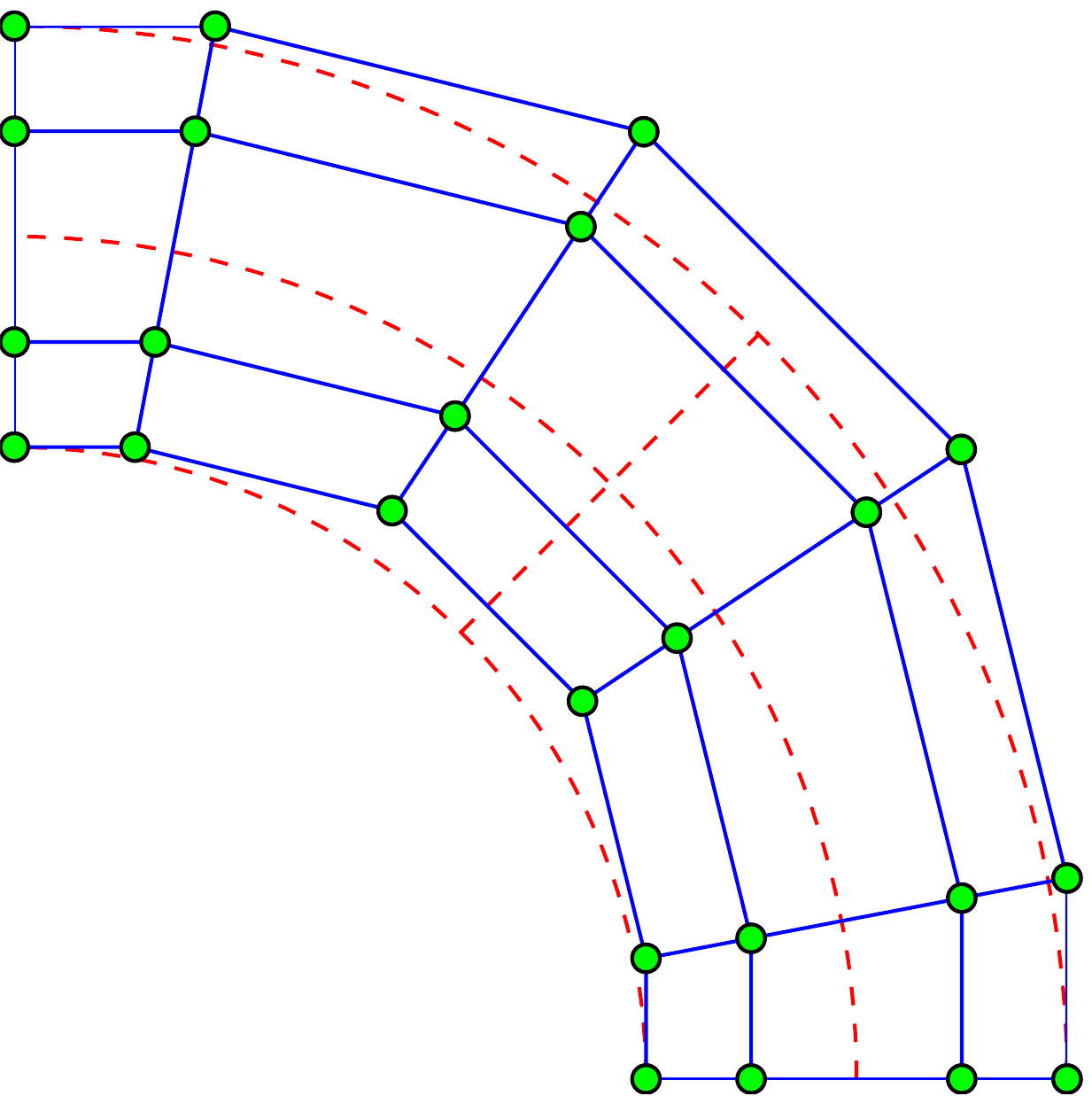}}\;\;\;\;\;
  \subfloat[After hp-refinement]{\includegraphics[width=0.31\textwidth]{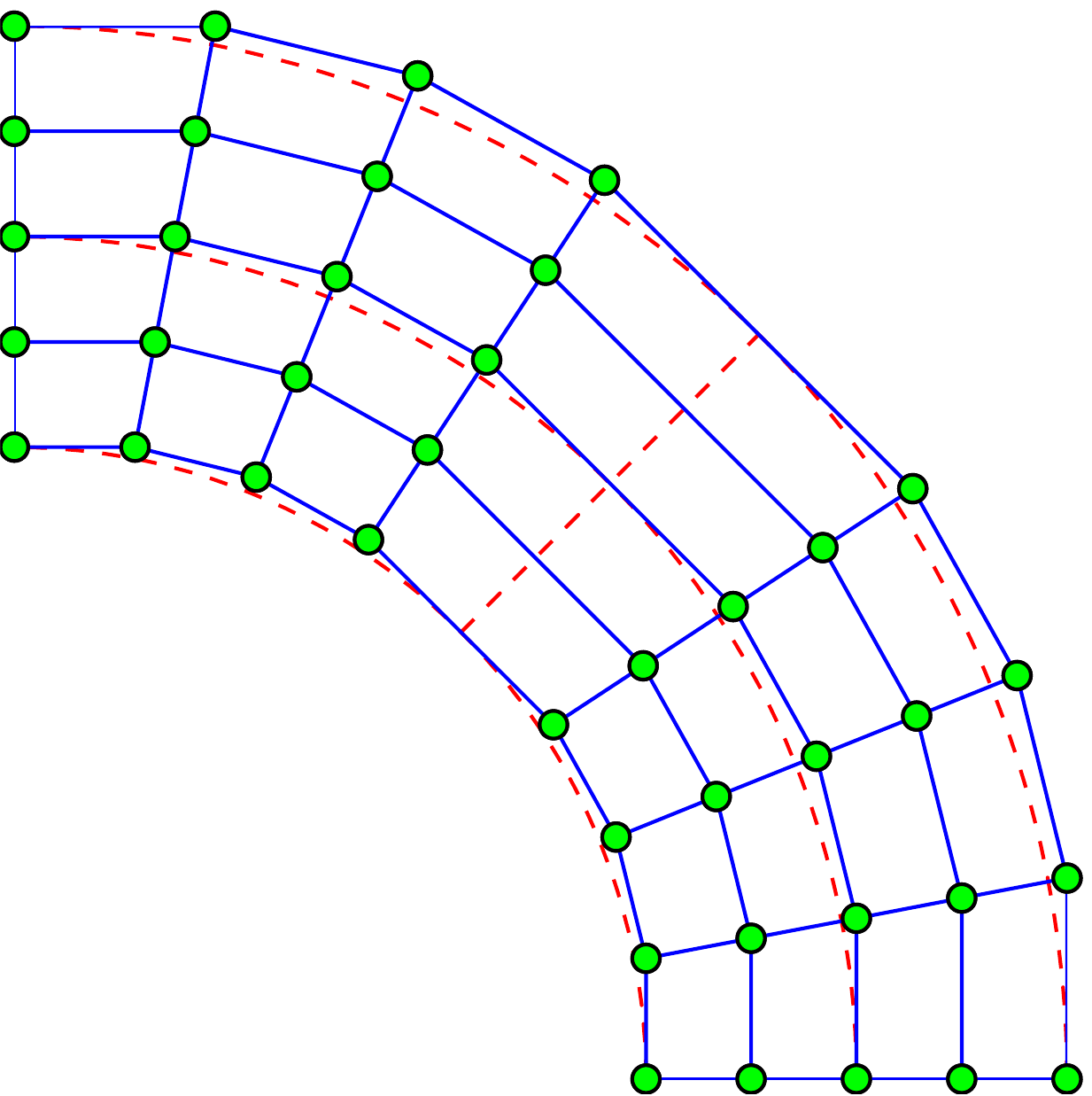}}\;\;\;\;\;
  \caption{Illustration of the utilization of the NURBS toolbox in building NURBS object (a).
From the initial mesh, different meshes can be obtained using either $h-$refinement, $p-$refinement or
combination thereof. As can be seen from (c) and (d), $k-$refinement (c) is more efficient than
$hp-$refinement (d). The function \textbf{plotMesh.m} in folder \textbf{meshing}
is used to plot NURBS mesh and 
control polygon.}
  \label{fig:refinement}
\end{figure}

For the refinement of NURBS, we reuse the NURBS Toolbox described in \cite{falcao_geopde_2011}.
We construct a NURBS surface as shown in Fig.\ref{fig:refinement}a. The corresponding Matlab
code is given in Listing \ref{2d-iga-geometry}. Using a uniform $h$-refinement (Listing \ref{2d-iga-hrefine})
that divides a knot span into two, 
we obtain the mesh given in Fig.\ref{fig:refinement}b.
If one needs to use $k-$refinement ($p-$refinement followed by $h-$refinement), 
then the code in Listing \ref{2d-iga-krefine} can be used (see Fig.\ref{fig:refinement}c).
Finally, Fig.\ref{fig:refinement}d gives the mesh which is obtained by the process in which
$h-$refinement is employed first and then $p-$refinement is performed (see Listing \ref{2d-iga-prefine}). 
After using the NURBS toolbox, the NURBS object is then converted to MIGFEM data structures
using the function \textbf{convert2DNurbs} located in folder \textbf{nurbs-util}.\\

\begin{snippet}[caption={Construct a NURBS surface using the NURBS toolbox. },
                   label={2d-iga-geometry},framerule=1pt]
    a = 0.3; b = 0.5; % inner/outer radius
    %% knots, control points
    uKnot = [0 0 0 1 1 1]; %quadratic basis
    vKnot = [0 0 1 1];     %linear basis
    % homogeneous coords (x*w,y*w,z*w,w)
    % 3 pts in u-dir, 2 pts in v-dir 
    controlPts          = zeros(4,3,2); 
    
    controlPts(1:2,1,1) = [a;0];
    controlPts(1:2,2,1) = [a;a;];
    controlPts(1:2,3,1) = [0;a];
    
    controlPts(1:2,1,2) = [b;0];
    controlPts(1:2,2,2) = [b;b];
    controlPts(1:2,3,2) = [0;b];
    
    controlPts(4,:,:)   = 1;
    controlPts(4,2,1) = 1/sqrt(2);
    controlPts(4,2,2) = 1/sqrt(2);
    % homogenous coordinates (x*w,y*w,z*w)
    controlPts(1:2,2,1) = controlPts(1:2,2,1)*fac;
    controlPts(1:2,2,2) = controlPts(1:2,2,2)*fac;
    %% build NURBS object
    solid = nrbmak(controlPts,{uKnot vKnot});
\end{snippet}

\begin{snippet1}[caption={h-refinement using the NURBS toolbox. },
                   label={2d-iga-hrefine},framerule=1pt]
    refineLevel = 2;
    for i=1:refineLevel
      uKnotVectorU = unique(uKnot);
      uKnotVectorV = unique(vKnot);
      % new knots along two directions
      newKnotsX =uKnotVectorU(1:end-1)+0.5*diff(uKnotVectorU);
      newKnotsY =uKnotVectorV(1:end-1)+0.5*diff(uKnotVectorV);
      
      newKnots  = {newKnotsX newKnotsY};
      solid     = nrbkntins(solid,newKnots);
      uKnot     = cell2mat(solid.knots(1));
      vKnot     = cell2mat(solid.knots(2));
    end
    convert2DNurbs % convert to the MIGFEM format
    plotMesh(controlPts,weights,uKnot,vKnot,p,q,res,'r-','try.eps');
\end{snippet1}

\begin{snippet1}[caption={k-refinement using the NURBS toolbox. },
                   label={2d-iga-krefine},framerule=1pt]
    % code from Listing 5
    % order elevation, p=p+2, q=q+1
    solid = nrbdegelev(solid,[2 1]);
    % then, knot insertion using the code from Listing 6
    convert2DNurbs % convert to the MIGFEM format
\end{snippet1}

\begin{snippet1}[caption={h-refinement followed by p-refinement using the NURBS toolbox. },
                   label={2d-iga-prefine},framerule=1pt]
    % code from Listing 5
    % then, knot insertion using the code from Listing 6
    % order elevation, p=p+2, q=q+1
    solid = nrbdegelev(solid,[2 1]);
    convert2DNurbs % convert to the MIGFEM format
\end{snippet1}

\begin{snippet1}[caption={A typical input file for MIGFEM. },
                   label={2d-iga-input},framerule=1pt]
    % code from Listing 5
    % code from Listing 6 (if only h-refinement)
    convert2DNurbs%convert to the MIGFEM format (if not done yet)
    generateIGA2DMesh % build the mesh (element connectivity)
\end{snippet1}

\subsection{Input file for MIGFEM} A typical input file is given in 
Listing \ref{2d-iga-input}. This file replaces the standard FE mesh file.
Note that, for backward compatibility with older versions of the MIGFEM code, 
some input files do not use the NURBS toolbox to create the NURBS object.
For those input files, it is however impossible to perform order elevation
and hence $k-$refinement. It is, therefore, recommended to create the
NURBS objects using the NURBS toolbox \cite{falcao_geopde_2011} which, 
besides the aforementioned refinement functionalities,
also supports many useful operations such as extrusion, rotation \etc

Incorporating NURBS into an existing FE code cannot be considered a trivial task due to the presence
of various spaces (index, parameter) that are not present in conventional FE codes,
B\'{e}zier extraction \cite{borden_isogeometric_2011,scott_isogeometric_2011} provides
a FE data structure that allows for a straightforward implementation of NURBS/T-splines into
any FE codes. \ref{sec:extraction} briefly presents this concept and its implementation in MIGFEM.

\subsection{XFEM implementation}\label{sec:xfem-implementation}

There are certainly different ways to implement XIGA. 
We present a way that reuse most of the tools present in an existing XFEM code. We use the approximate Q4 mesh 
used for visualization as discussed in 
Section \ref{sec:visualization} for selection of enriched control points.
The level set values defining the crack at the vertices of this mesh are then
computed. Based on these level sets, elements cut by the crack and elements
containing the crack tip can be determined \cite{Sukumar1}. For example,
element ABCD in Fig. \ref{fig:enriched-nodes} is cut by the crack. In a FEM
context, its four nodes are then enriched using the Heaviside function. In an isogeometric
framework, however, the control points associated to this element are
enriched. Listing \ref{list-selection-enriched-nodes} details the Matlab\textsuperscript{\textregistered{}}
implementation of this process. Note that this listing is taken directly from
our XFEM code \cite{nguyen_meshless_2008} with only one small modification and thus the proposed technique is 
considered simpler than that adopted in \cite{ghorashi_extended_2012}. We emphasize that the crack geometry is 
defined in the physical space to keep the usual XFEM notation.

\begin{snippet}[caption={Selection of enriched control points/nodes. },
                   label={list-selection-enriched-nodes},framerule=1pt]
    enrich_node = zeros(noCtrPts,1);
    count1 = 0;
    count2 = 0;
    for iel = 1 : numelem
        sctr    = elementV(iel,:);
        sctrIGA = element(iel,:);
        phi  = ls(sctr,1); % normal level set
        psi  = ls(sctr,2); % tangent level set
        if ( max(phi)*min(phi) < 0 )
            if max(psi) < 0
                count1 = count1 + 1 ; % one split element
                split_elem(count1) = iel;
                enrich_node(sctrIGA)   = 1;
            elseif max(psi)*min(psi) < 0
                count2 = count2 + 1 ; % one tip element
                tip_elem(count2) = iel;
                enrich_node(sctrIGA)   = 2;
            end
        end
    end
    split_nodes = find(enrich_node == 1);
    tip_nodes   = find(enrich_node == 2);		   
\end{snippet}

\begin{rmk} 
Note that for simple geometries as those tackled in this paper, the use of
level sets is certainly not necessary. More generally, describing open
surfaces (cracks) with level sets requires two level sets functions that, for
crack growth simulations, must be reinitialised for stability (this decreases
accuracy) and reorthogonalised every few time steps. This is particularly
cumbersome and probably explains along with the difficulties in dealing intersecting
and branching cracks the recent trend of research efforts in the area of phase
field models of fracture, see \eg \cite{NME:NME2861,Borden201277}, and the thick level set
method \cite{NME:NME3069}.
\end{rmk}

\noindent \textbf{Crack visualization} 
We have $H(\vm{x}^+)-H(\vm{x}^-)=2$ and
$B_1(\vm{x}^+)-B_1(\vm{x}^-)=2\sqrt{r}$, therefore the displacement jump
at a point $\vm{x}$ on the crack face is given by

\begin{equation}
	\jump{\vm{u}}({\rm {\bf x}})= 2\sum\limits_{J\in
\mathcal{S}^{c}} {R_J ({\rm {\bf x}}) {\rm {\bf
a}}_J } + 2\sqrt{r}\sum\limits_{K\in \mathcal{S}^{f}} {R_K ({\rm {\bf x}}) {{\rm {\bf b}}_{K}^1}
}.
\label{eq:jump}
\end{equation}

\noindent Note that the other branch functions $B_\alpha$, $\alpha=2,3,4$ are continuous functions and thus
do not contribute to the displacement jump.

\begin{figure}[htbp]
  \centering 
  \includegraphics[width=0.65\textwidth]{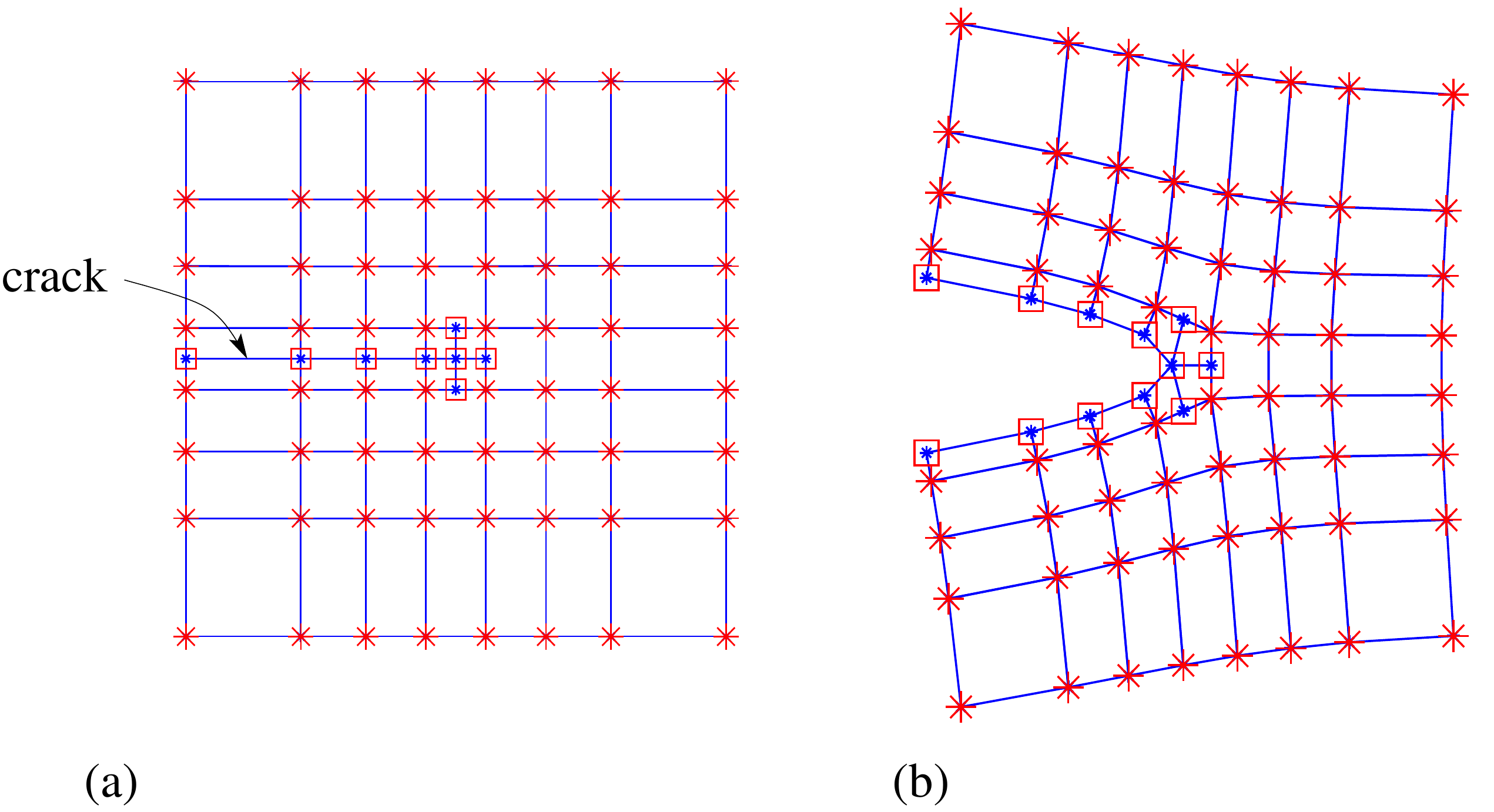}
  \caption{Crack visualization in XIGA: (a) build a mesh that is compatible
  to the crack by introducing double nodes along the crack (square nodes) and (b) 
  assign the displacement jumps to these new nodes. Note that the square
  nodes in the tip element are used only for compatibility purposes. Exact mode I
  displacements are imposed on the bottom, right and top edge using the
  Lagrange multiplier method while Neumann BCs from the exact stress field are enforced on the left edge.
  We refer to \cite{nguyen_meshless_2008} for a detailed description of this standard problem. }
  \label{fig:crack-visualization}
\end{figure}

\begin{figure}[htbp]
  \centering 
  \includegraphics[width=0.25\textwidth]{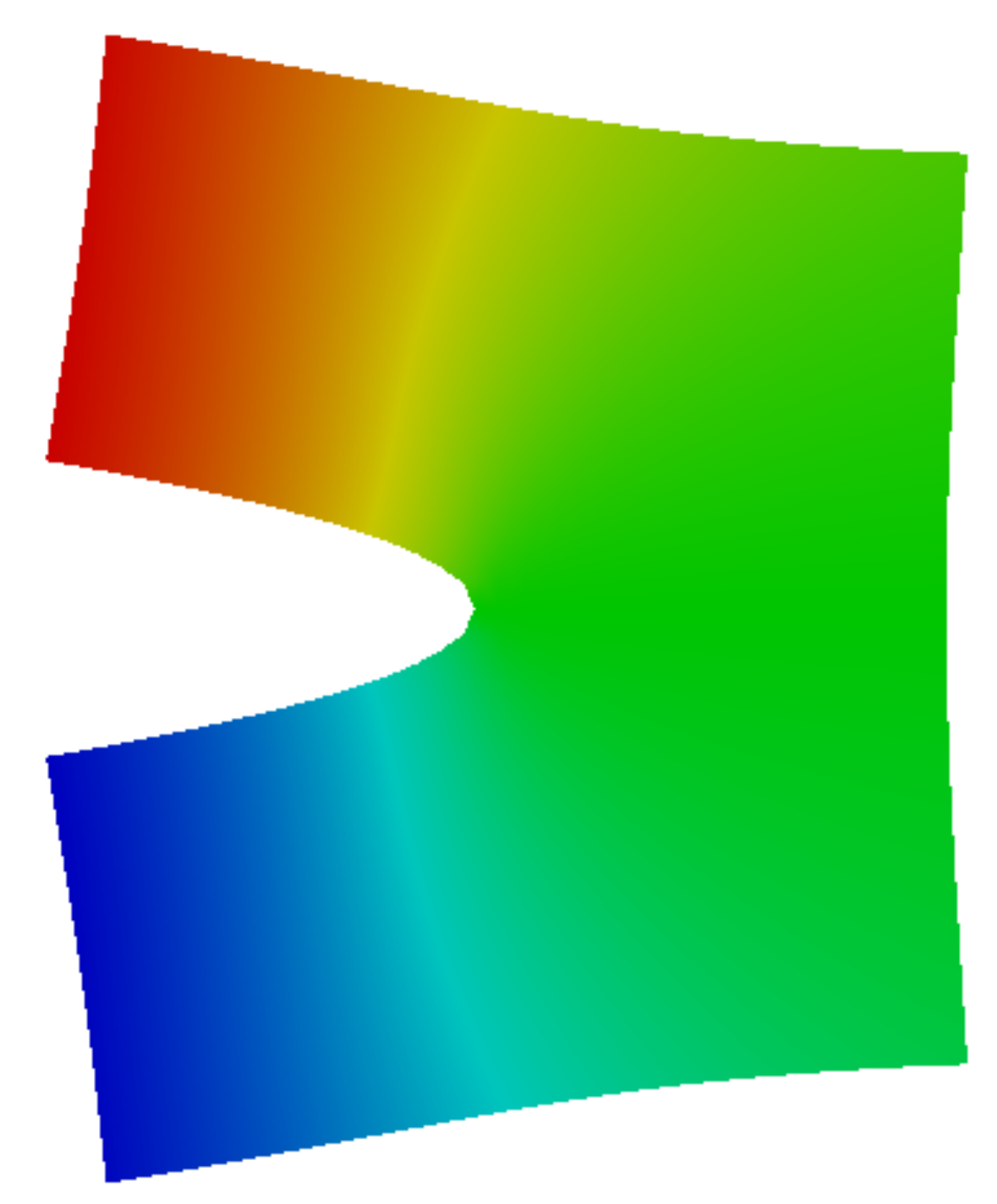}\;\;\;\;\;
  \includegraphics[width=0.23\textwidth]{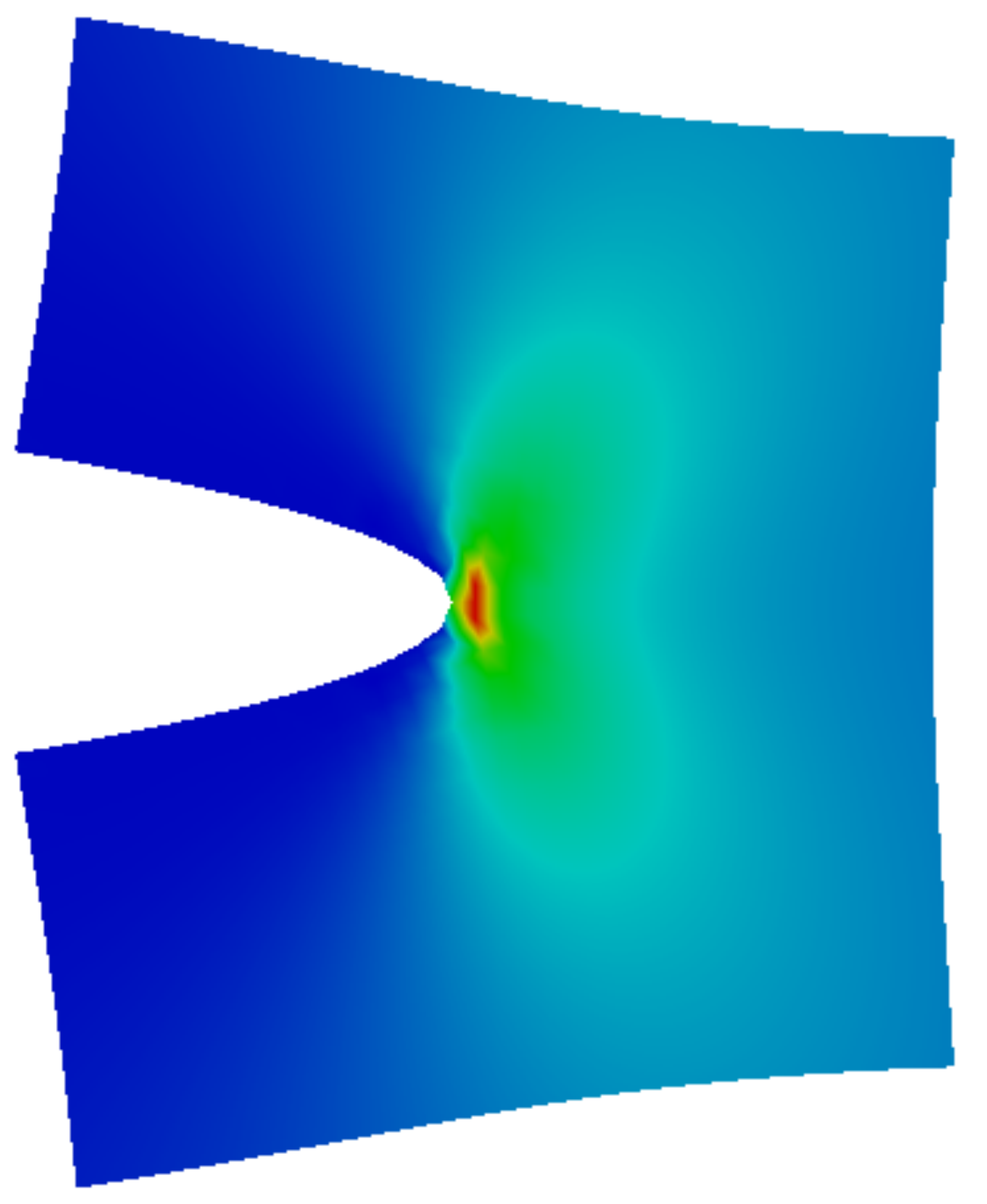}
  \caption{Contour plots on a cracked mesh: (a) vertical displacement and (b)
  normal stress in the vertical direction. }
  \label{fig:crack-visualization-result}
\end{figure}

Fig. \ref{fig:crack-visualization} illustrates the idea for crack
visualization with the script $\textbf{post-processing/crackedMeshNURBS.m}$ providing implementation details. The contour plots of the displacement and stress field
of a mode I cracked sample are given in Fig. \ref{fig:crack-visualization-result}. Note that the stresses at points on
the crack are simply set to zero (traction-free cracks) and the stresses of
the new nodes of the tip element are interpolated from the values of the four
nodes of this Q4 element.

\begin{rmk}
   As is the case for XFEM, integration over elements cut by the cracks
   usually requires subdivision of the elements into integration subcells. We refer to
   \cite{sunda-xfem-no-subcells} for a recent discussion on this issue. In addition to this
   popular technique, a simple integration rule is also provided--elements
   crossed by the crack and tip-enriched elements are numerically integrated
   using a regular Gauss-Legendre quadrature with a large number of Gauss points as
   done in \cite{Elguedj2006501}.
\end{rmk}

\section{Structural mechanics}\label{sec:structural-mechanics}

Thanks to the high order continuity provided by NURBS/T-splines, the implementation of rotation-free
thin beam/plate/shell elements becomes direct and simple. 
In this section, we are going to present the implementation
of a rotation-free IGA Kirchhoff plate formulation (rotation free shell elements
can be found in folder \textbf{structural-mechanics}). The 
plate geometry and the deflection are both approximated by NURBS. At control points there is only
one unknown-the deflection or transverse displacement.
For simplicity only isotropic elastic plates are considered. We refer to \cite{reddy_book}
for a treatment of plate theories.

The element stiffness matrix is defined as

\begin{equation}
\vm{K}_e = \int_{\Omega_e} \vm{B}_e\trans \vm{D} \vm{B}_e \di \Omega,
\end{equation}

\noindent where the constitutive matrix $\vm{D}$ reads

\begin{equation}
\vm{D}=\frac{Eh^3}{12(1-\nu^2)}\begin{bmatrix}
1 & \nu & 0\\
\nu & 1 & 0\\
0 & 0 & 0.5(1-\nu)
\end{bmatrix},
\end{equation}

\noindent where $E,\nu$ are the Young's modulus and Poisson's ratio, respectively;
$h$ denotes the plate thickness and the element displacement-curvature matrix $\vm{B}_e$ that contains
the second derivatives of the shape functions is given by

\begin{equation}
\vm{B}_e = \begin{bmatrix} 
R_{1,xx} & R_{2,xx} & \cdots R_{n,xx}\\
R_{1,yy} & R_{2,yy} & \cdots R_{n,yy}\\
2R_{1,xy} & 2R_{2,xy} & \cdots 2R_{n,xy}\\
\end{bmatrix},
\end{equation}
where $n$ denotes the number of basis functions of element $e$ and  $R_{A,xx}\equiv d^2R_A/dx^2$.

\subsection{Boundary conditions}\label{sec-bcs}

\begin{figure}[h!]
  \centering 
  \includegraphics[width=0.26\textwidth]{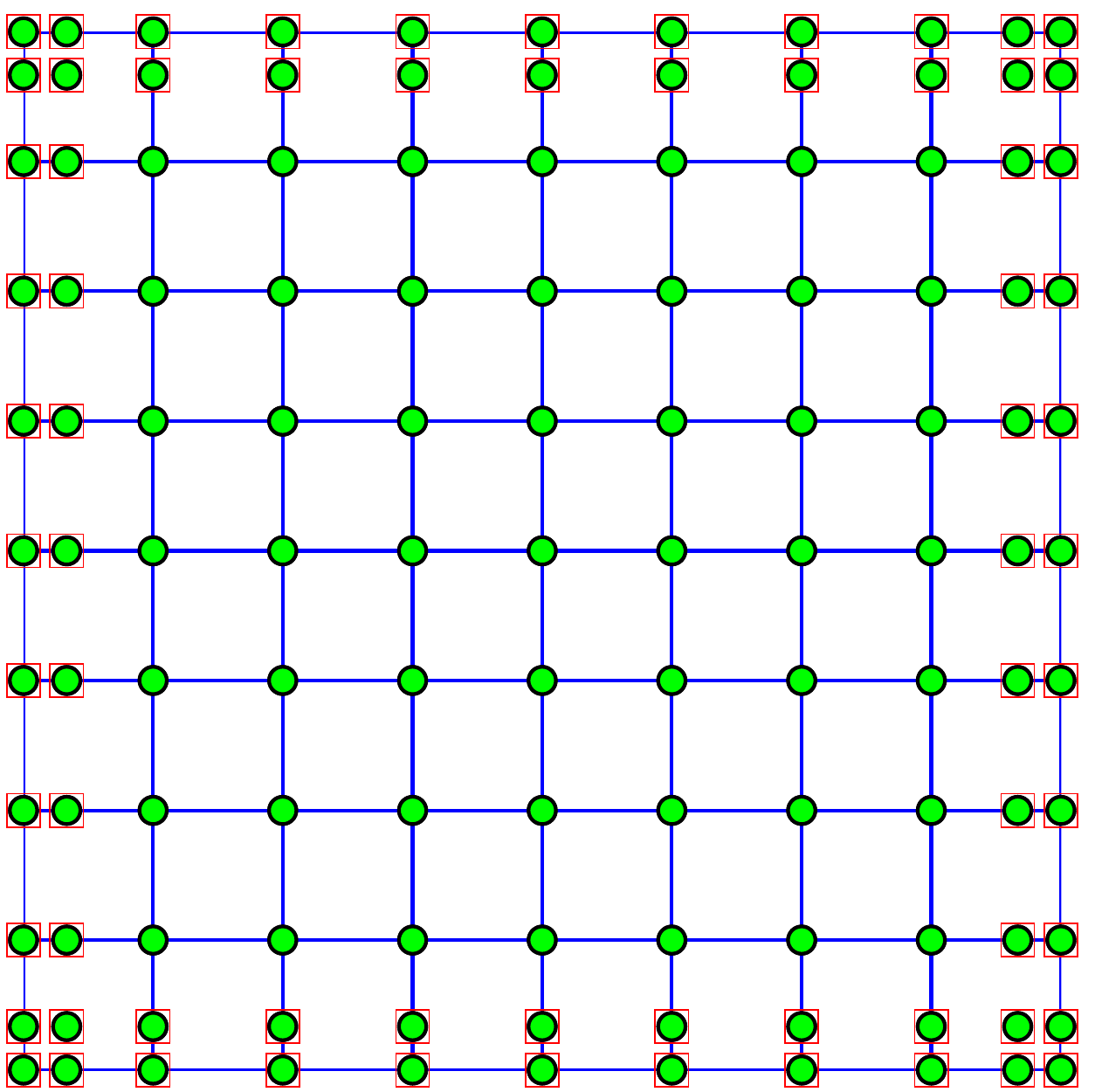}
  \caption{Enforcing BCs for a fully clamped plate:
          simply fixing the deflections of two rows of control points around the clamped boundary.
          Note that the set of CPs next to the boundary CPs are not artificially added to impose
          the rotations. They are simply the CPs defining the geometry of the plate.}
  \label{fig:plate-bc} 
\end{figure}

For clamped BCs one needs to fix the rotations. The nodal unknowns are, however, only the transverse
displacements $\vm{w}$. To fix the rotation of
a boundary, we simple fix two row of control points at the boundary
\cite{kiendl_isogeometric_2009} because these control points define the tangent
of the surface at the boundary (cf. Fig.~\ref{fig:splinecurve1}), see Fig.~\ref{fig:plate-bc}. 

\subsection{Symmetry boundary conditions}\label{sec-sym} 

\begin{figure}[h!]
  \centering 
  \psfrag{fix}{fix}
  \psfrag{fix}{fix}
  \psfrag{sym}{sym}
  \psfrag{sym}{sym}
  \psfrag{cou}{coupled nodes}
  \includegraphics[width=0.55\textwidth]{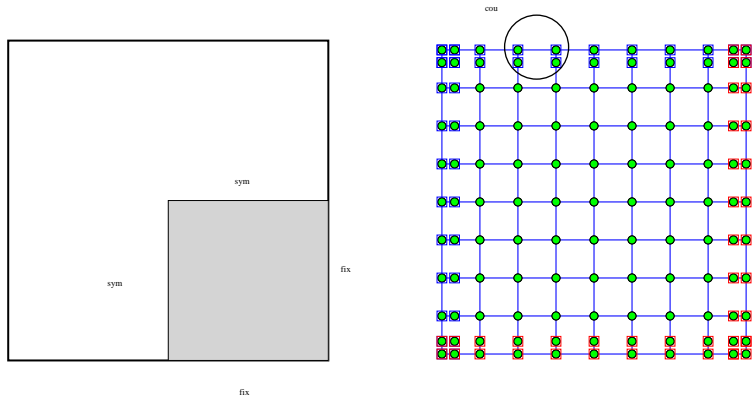}
  \caption{A fully clamped square plate: 1/4 model is analyzied using appropriate symmetry BCs. Along the
 symmetry lines, the rotation is fixed which can be achieved by enforcing the deflection of two rows of control
 points that define the tangent of the plate to have the same value.}
  \label{fig:plate-bc-sym} 
\end{figure}

Fig. \ref{fig:plate-bc-sym} illustrates 
the use of symmetry boundary conditions when
only 1/4 of the plate is modelled. Along the symmetry lines, the rotation should be zero which
can be enforced by constraining the deflection ($w$) of two rows of control points
along these lines together. These contraints can be implemented using a simple penalty
technique as shown in Listing \ref{list-symmetry}.

\begin{snippet}[caption={Enforcing symmetry BCs.},
	label={material}, label={list-symmetry},framerule=1pt]
    w     = 1e7;
    penaltyStiffness = w*[1 -1;-1 1];
    for i=1:length(topNodes)
        sctr  = [topNodes(i) nextToTopNodes(i)];
        K(sctr,sctr) = K(sctr,sctr) + penaltyStiffness;
    end
    % the same for the left two rows    
\end{snippet}

For completeness, we give the implementation of the rotation free Kirchhoff plate elements 
in Listing \ref{Ke-Kirchhoff-plate}. 
Note that we have skipped codes that are common with IGA code for 2D continua. 

\begin{snippet1}[caption={Computation of stiffness matrix for the rotation free Kirchhoff plate.},
	label={material}, label={Ke-Kirchhoff-plate},framerule=1pt]
    for gp=1:size(W,1)
        pt      = Q(gp,:); wt      = W(gp);
        % shape functions, first and second derivatives w.r.t natural coords
        [R dRdxi dRdeta dR2dxi dR2det dR2dxe]=NURBS2DBasis2ndDers([Xi; Eta],p,q,uKnot,vKnot,weights');
        jacob  = [dRdxi; dRdeta]          * pts; % 2x2 matrix
        jacob2 = [dR2dxi; dR2det; dR2dxe] * pts; % 3x2 matrix
        J1    = det(jacob);
        dxdxi = jacob(1,1); dydxi = jacob(1,2);
        dxdet = jacob(2,1); dydet = jacob(2,2);
        j33   = [dxdxi^2     dydxi^2     2*dxdxi*dydxi;
                 dxdet^2     dydet^2     2*dxdet*dydet;
                 dxdxi*dxdet dydxi*dydet dxdxi*dydet+dxdet*dydxi];
        % Jacobian inverse and spatial 1st and 2nd derivatives
        invJacob   = inv(jacob);
        dRdx       = invJacob*[dRdxi;dRdeta];
        dR2dx      = inv(j33)*([dR2dxi; dR2det; dR2dxe]-jacob2*dRdx);
        % B matrix
        B          = dR2dx;
        B(3,:)     = B(3,:)*2;
        % compute elementary stiffness matrix, sctr=element(e,:)
        K(sctr,sctr) = K(sctr,sctr) + B' * C * B * J1 * J2 * wt;
    end
\end{snippet1}

\section{Verification examples}\label{examples}

In this section, numerical examples in linear elasticity and linear elastic
fracture mechanics in 2D and 3D are presented with the purpose to 
serve as a set of verification examples for MIGFEM. They include  
an infinite plate with a circular hole under constant in-plane tension, 
the pinched cylinder, an edge cracked plate in tension, a
three-dimensional mode I fracture problem and a large deformation thin shell problem. 
Unless otherwise stated, standard direct imposition of Dirichlet BCs is used.
Units are standard International System (SI) units.

\subsection{Two and three dimensional solid mechanics}

\subsubsection{Infinite plate with a circular hole}

\begin{figure}[htbp]
  \centering 
  \psfrag{ex}[c]{Exact traction}
  \psfrag{ex}[c]{Exact traction}
  \psfrag{sy}{Symmetry}
  \psfrag{sy}{Symmetry}
  \psfrag{R}{$R$}
  \psfrag{r}{$r$}
  \psfrag{l}{$L$}
  \psfrag{theta}{$\theta$}
  \psfrag{e}{$E$}
  \psfrag{nu}{$\nu$}
  \psfrag{eq}{$=$} \psfrag{eq}{$=$} \psfrag{eq}{$=$} \psfrag{eq}{$=$}
  \psfrag{v1}{1000}\psfrag{v2}{0.3}\psfrag{v3}{1}\psfrag{v4}{4}
  \psfrag{A}{A}\psfrag{B}{B}\psfrag{C}{C}\psfrag{D}{D}\psfrag{E}{E}
  \psfrag{x}{$x$}\psfrag{y}{$y$}
  \includegraphics[width=0.5\textwidth]{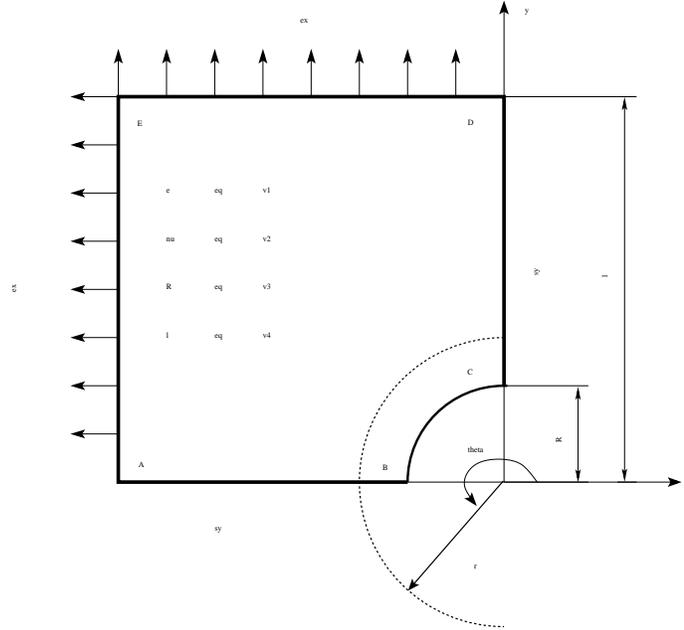}
  \caption{Infinite plate with a circular hole under constant in-plane tension: quarter
  model. Boundary conditions include:  $u_y=0$ (AB),  $u_x=0$ (CD), 
   $\overline{\mathbf t}\trans$ = $(-\sigma_{xx}, -\sigma_{xy})$ (AE),
    $\overline{\mathbf t}\trans =(\sigma_{xy},\sigma_{yy})$ (ED).}
  \label{fig:plateHole-geo}
\end{figure}

The problem considered is that of an infinite plate with a circular hole in the centre under
constant in-plane tension at infinity as shown in Fig. 
\ref{fig:plateHole-geo} where, due to symmetry, only a quarter of the plate is modeled. 
The plate dimension is taken to be $L \times L$ and the circular hole has a radius $R$.
The exact stress field in the plate is given by

\begin{subequations}

\begin{equation}
	\sigma_{xx}(r,\theta) = 1 - \frac{R^2}{r^2}\biggl ( \frac{3}{2}\cos
2\theta+\cos 4 \theta  \biggr ) + \frac{3}{2}\frac{R^4}{r^4}\cos
4\theta
\end{equation}

\begin{equation}
	\sigma_{yy}(r,\theta) = - \frac{R^2}{r^2}\biggl ( \frac{1}{2} \cos
2\theta-\cos 4 \theta  \biggr ) - \frac{3}{2}\frac{R^4}{r^4}\cos
4\theta
\end{equation}

\begin{equation}
\sigma_{xy}(r,\theta) = - \frac{R^2}{r^2}\biggl ( \frac{1}{2}\sin
2\theta+\sin 4 \theta  \biggr ) + \frac{3}{2}\frac{R^4}{r^4}\sin
4\theta,
\end{equation}
\label{hole_exact}
\end{subequations}

\noindent where $r,\theta$ are the usual polar coordinates centered
at the center of the hole.

\begin{figure}[htbp]
  \centering 
  \includegraphics[width=0.3\textwidth]{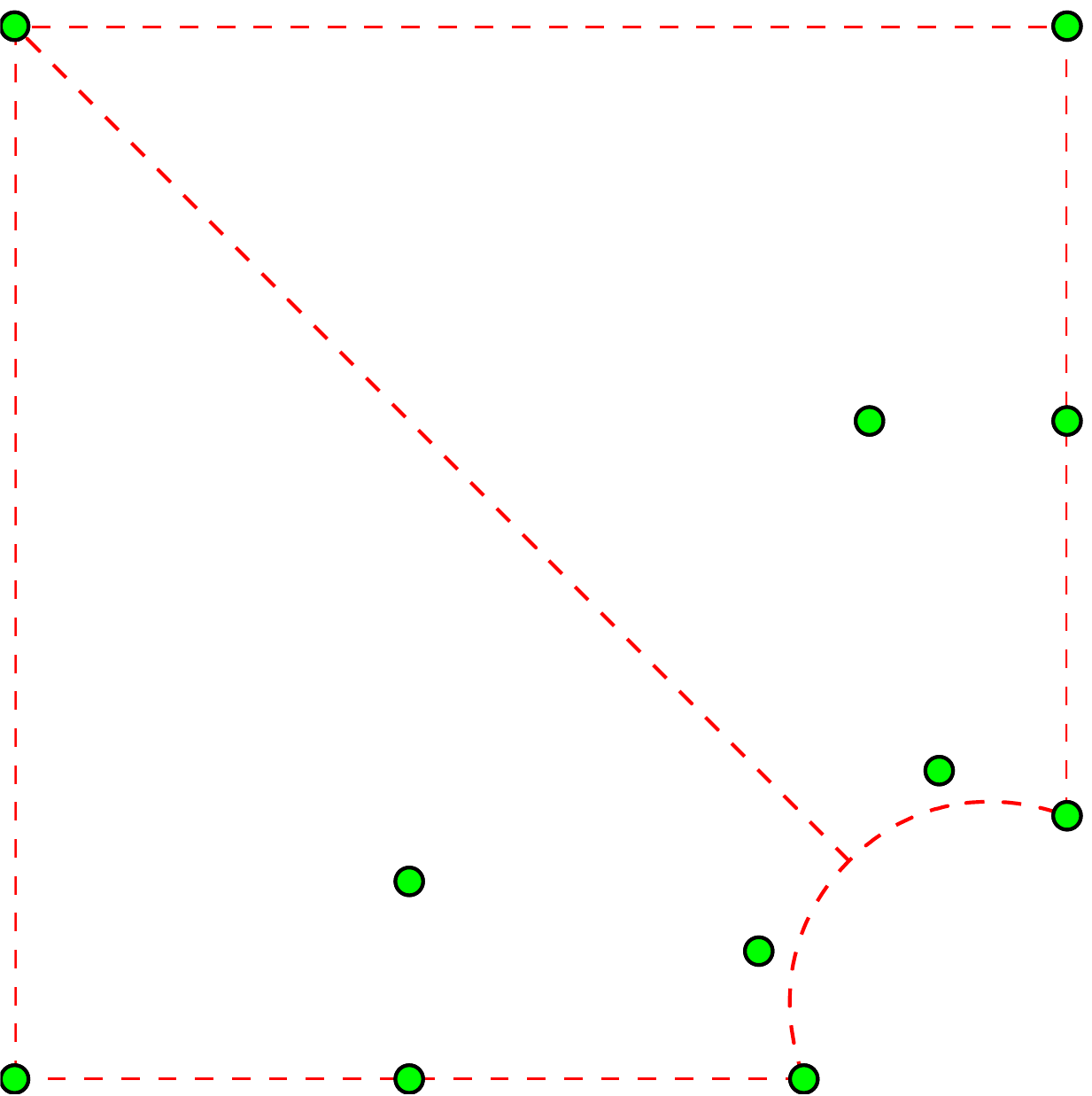}\;\;\;\;
  \includegraphics[width=0.3\textwidth]{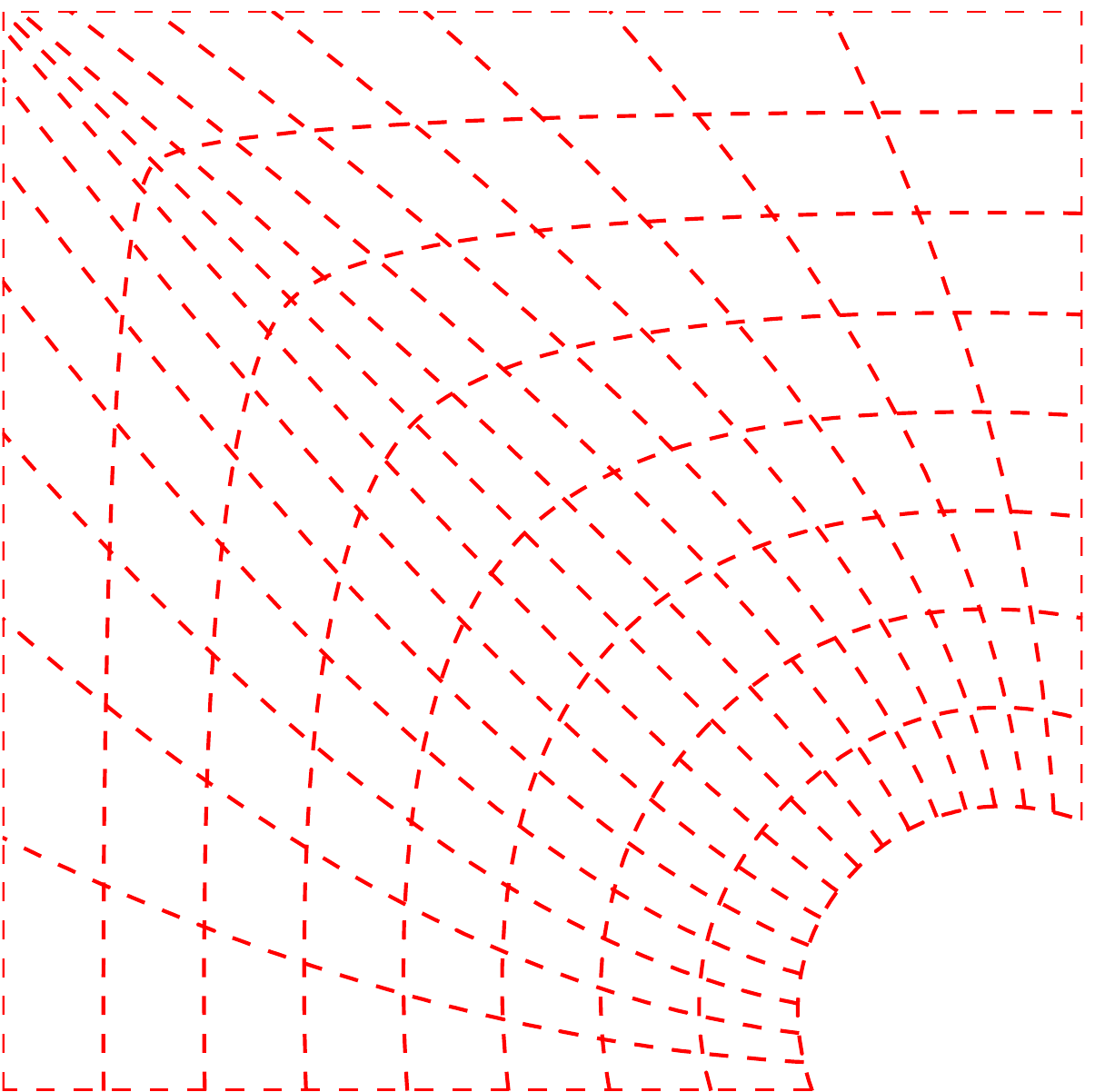}
  \caption{Plate with a hole: coarse mesh of 2 bi-quadratic elements (left) and refined mesh (right).}
  \label{fig:hole-mesh}
\end{figure}

\begin{figure}[htbp]
  \centering 
  \includegraphics[width=0.4\textwidth]{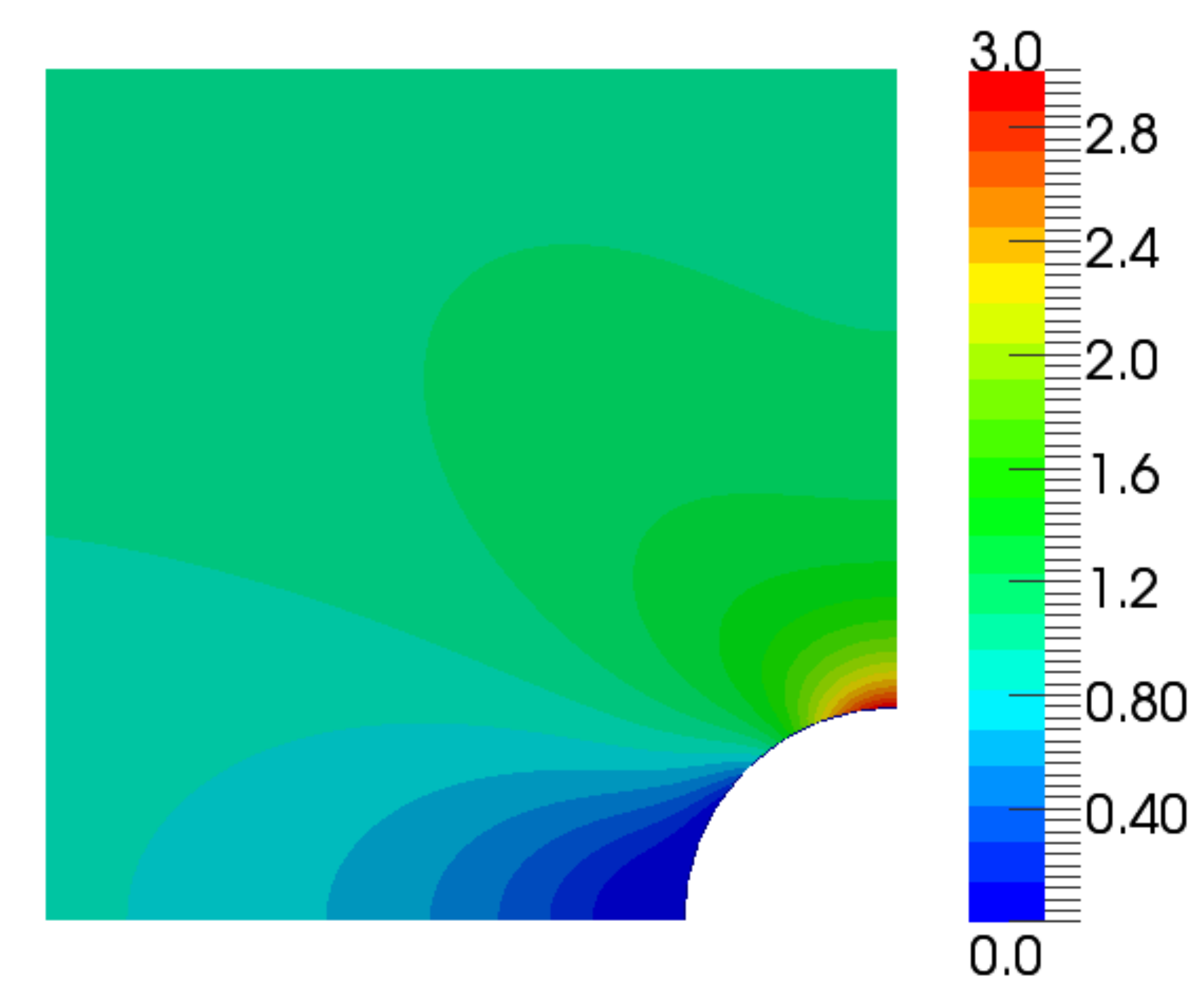}
  \caption{Plate with a hole: distribution of numerical $\sigma_{xx}$ obtained
  with a $32\times16$ quadratic mesh having 4488 dofs.}
  \label{fig:hole-stress}
\end{figure}

The material
properties are specified as $E=10^3$,  Poisson's ratio  $\nu=0.3$ and the
geometry is such that $L=4$, $R=1$. A plane stress condition is assumed.
The problem is solved with quadratic NURBS meshes such as those shown in Fig. \ref{fig:hole-mesh}.
The control points and weights for the coarsest mesh can be found in \cite{hughes_isogeometric_2005} or
file \textbf{plateHoleCkData.m}.
Fig. \ref{fig:hole-stress}, generated in Paraview, 
illustrates the contour plot of numerical $\sigma_{xx}$. Note that the stress
concentration at point ($R,3\pi/2$) is well captured and a smooth stress
field is obtained throughout.

\begin{rmk}
Using the visualization technique described in Section \ref{sec:visualization}
for this problem, a note should be made on the evaluation of the stress field at the top left corner
where there are two coincident control points. This causes a singular Jacobian
matrix. Therefore at this corner, the stresses at a point
slightly shifted from the original position are used.
\end{rmk}

\subsubsection{Pinched cylinder}

\begin{figure}[htbp]
\centering
\psfrag{x}{$x$}\psfrag{y}{$y$}\psfrag{z}{$z$}
\psfrag{E}{$E$}\psfrag{eq}{=}\psfrag{valE}{10}
\psfrag{nu}{$\nu$}\psfrag{eq}{=}\psfrag{valNu}{10}
\psfrag{P}{$P$}\psfrag{eq}{=}\psfrag{valP}{10}
\psfrag{R}{$R$}\psfrag{eq}{=}\psfrag{valR}{10}
\psfrag{L}{$L$}\psfrag{eq}{=}\psfrag{valL}{10}
\psfrag{t}{$t$}\psfrag{eq}{=}\psfrag{valT}{10}
\includegraphics[width=0.5\textwidth]{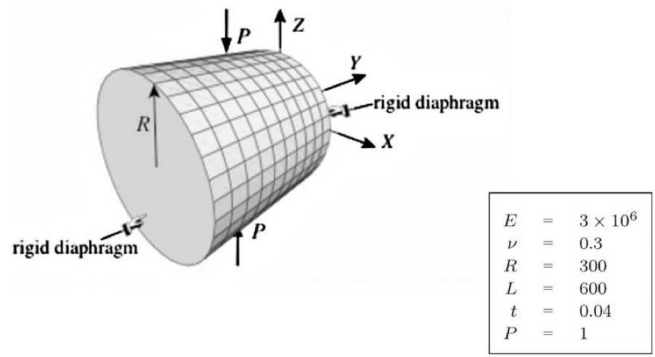}
\caption{Pinched cylinder. Problem description and data \cite{felippa:note2}.}
\label{fig:pinchedCylinderData}
\end{figure}

\begin{figure}[htbp]
  \centering 
  \includegraphics[width=0.75\textwidth]{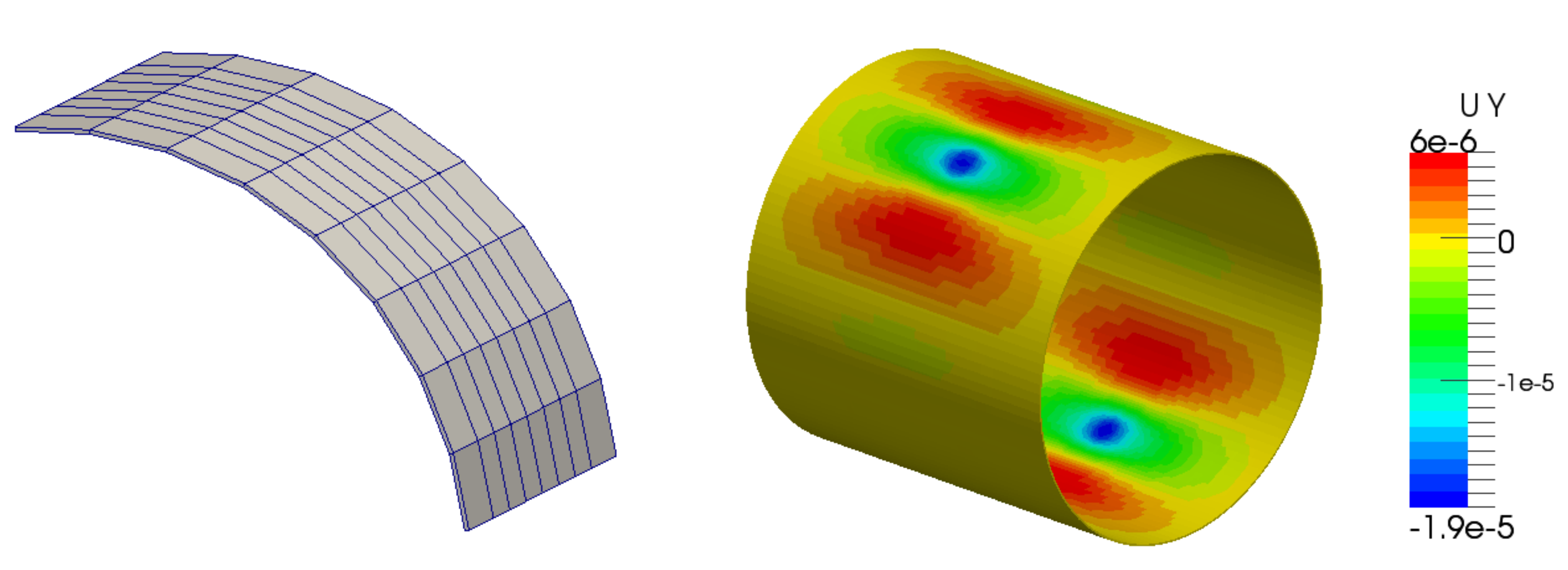}
  \caption{Pinched cylinder: (a) mesh of one octant of the cylinder and (b) 
  contour plot of the displacement in the direction of the point load.}
   \label{fig:pinchedCylinderResult}
\end{figure}

In order to demonstrate the performance of the 3D IGA implementation, we consider the pinched cylinder
problem as shown in Fig. \ref{fig:pinchedCylinderData}. Note that we discretize the shell with solid
NURBS elements. Due to symmetry, only 1/8 of the model is analyzed. 
A tri-quadratic NURBS mesh
($p=q=r=2$) was used for the computation. Details can be found in the
file \textbf{igaPinchedCylinder.m}. Fig. \ref{fig:pinchedCylinderResult}
shows the mesh and the contour plot of the displacement in the point load
direction. Post-processing is done in Paraview and we refer to Section 
\ref{sec:visualization} for details.  We recognise that the problem under
consideration is a shell like structure that would be more accurately
modelled using appropriate shell elements, but the example is merely intended
to illustrate the ability of the method to analyse 3D geometries.

\subsection{Two and three dimensional fracture mechanics}

\subsubsection{Edge cracked plate in tension}

\begin{figure}[htbp]
  \centering 
  \psfrag{si}{$\sigma$}\psfrag{h}{$h$}\psfrag{a}{$a$}\psfrag{b}{$b$}
  \includegraphics[width=0.2\textwidth]{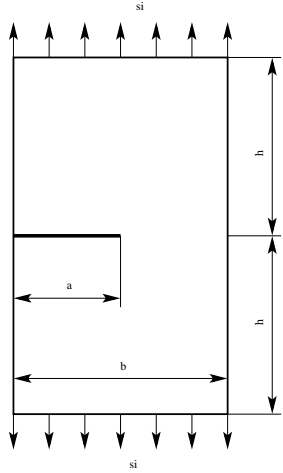}
  \caption{Edge cracked plate in tension: geometry and loading.}
  \label{fig:edge-crack-geo}
\end{figure}

A plate of dimension $b \times 2h$ is loaded by a tensile stress $\sigma = 1.0$ 
along the top edge and bottom edge as shown in Fig. \ref{fig:edge-crack-geo}. 
In the computation, the displacement along the $y$-axis 
is fixed at the bottom edge and the bottom left corner is fixed in
both $x$ and $y$ directions. The material parameters are $E=10^3$ and $\nu=0.3$. A plane
strain condition is assumed.
The reference mode I stress intensity factor (SIF) for this problem is given in
\cite{tada_1985} and is calculated as

\begin{equation}
	K_I=F\left(\frac{a}{b}\right)\sigma\sqrt{\pi a},
\end{equation}

\noindent where $a$ is the crack length, $b$ is the plate width, and $F(a/b)$
is an empirical function. For $a/b \le 0.6$, function $F$ is given by

\begin{equation}
	F\left(\frac{a}{b}\right)=1.12-0.23\left(\frac{a}{b}\right)+
	10.55\left(\frac{a}{b}\right)^2-21.72\left(\frac{a}{b}\right)^3
	+30.39\left(\frac{a}{b}\right)^4.
\end{equation}

\noindent In the present implementation, this problem is solved using
both XFEM and an extended
isogeometric formulation. The SIF is computed using an interaction integral.
We refer to \cite{mos_finite_1999} for details.

\begin{figure}[htbp]
  \centering 
  \subfloat[XFEM mesh (Q4)]{\includegraphics[width=0.2\textwidth]{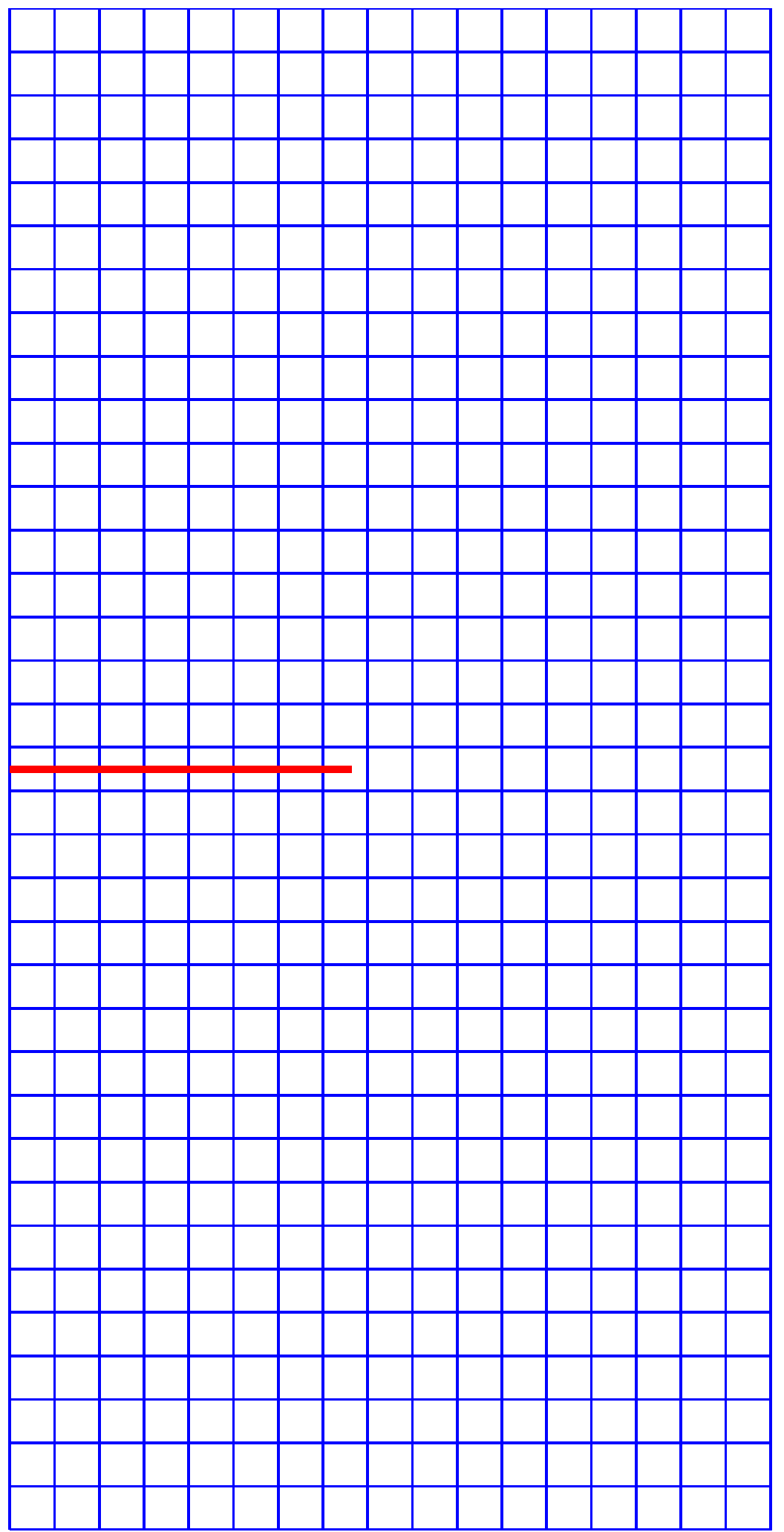}}\;\;\;\;
  \subfloat[XIGA mesh (cubic
  Bsplines)]{\includegraphics[width=0.204\textwidth]{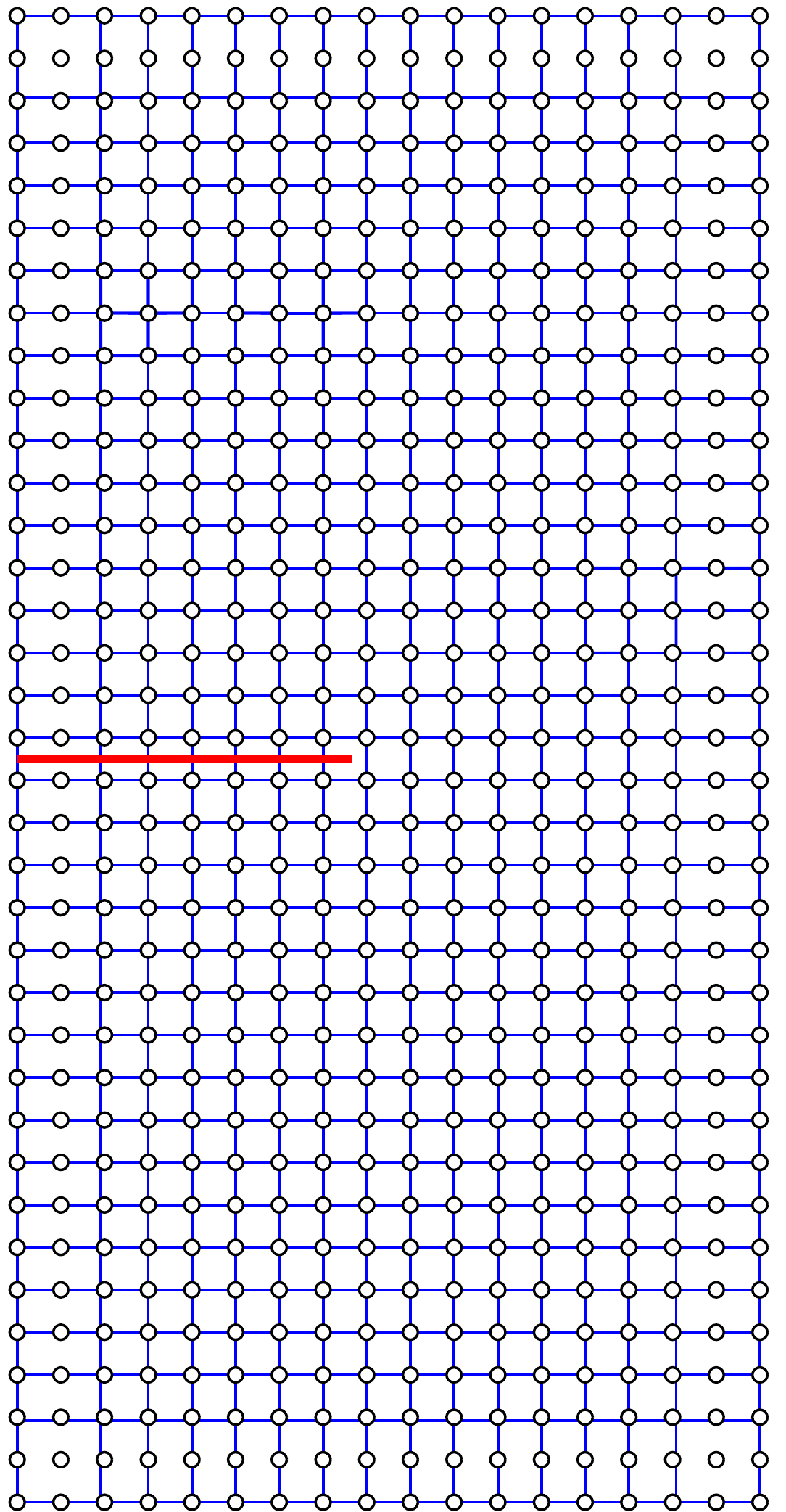}}\;\;
  \subfloat[XIGA mesh, enriched nodes]{\includegraphics[width=0.4\textwidth]{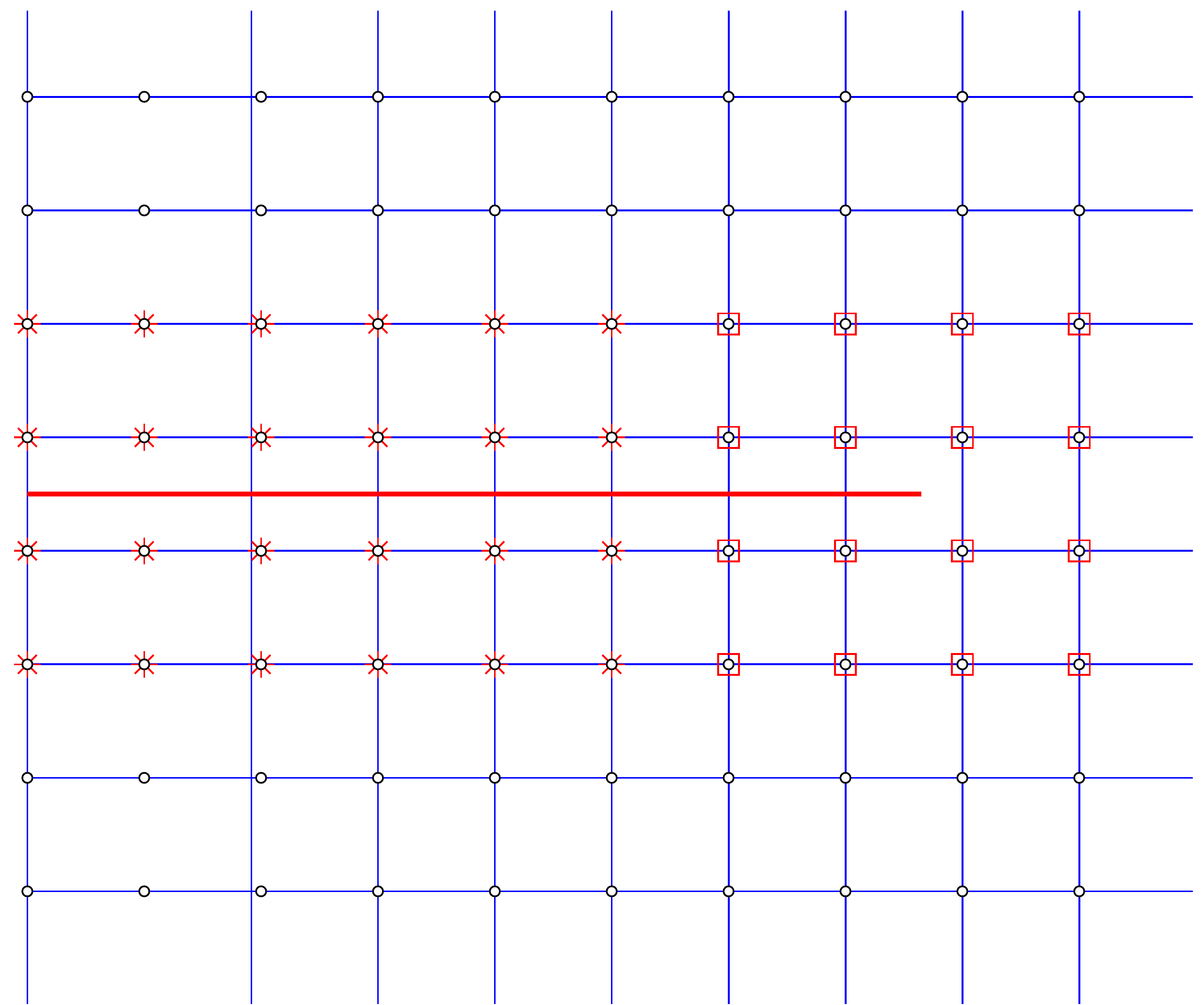}}
  \caption{Edge cracked plate: XFEM and XIGA meshes. Both have the same number
  of displacement dofs of 1296. The thick line denotes the crack. Square nodes
  denote tip enriched nodes whereas star nodes represent Heaviside enriched
  nodes.}
  \label{fig:edge-crack-meshes} 
\end{figure}

\begin{figure}[htbp]
  \centering 
  \subfloat[XFEM]{\includegraphics[width=0.25\textwidth]{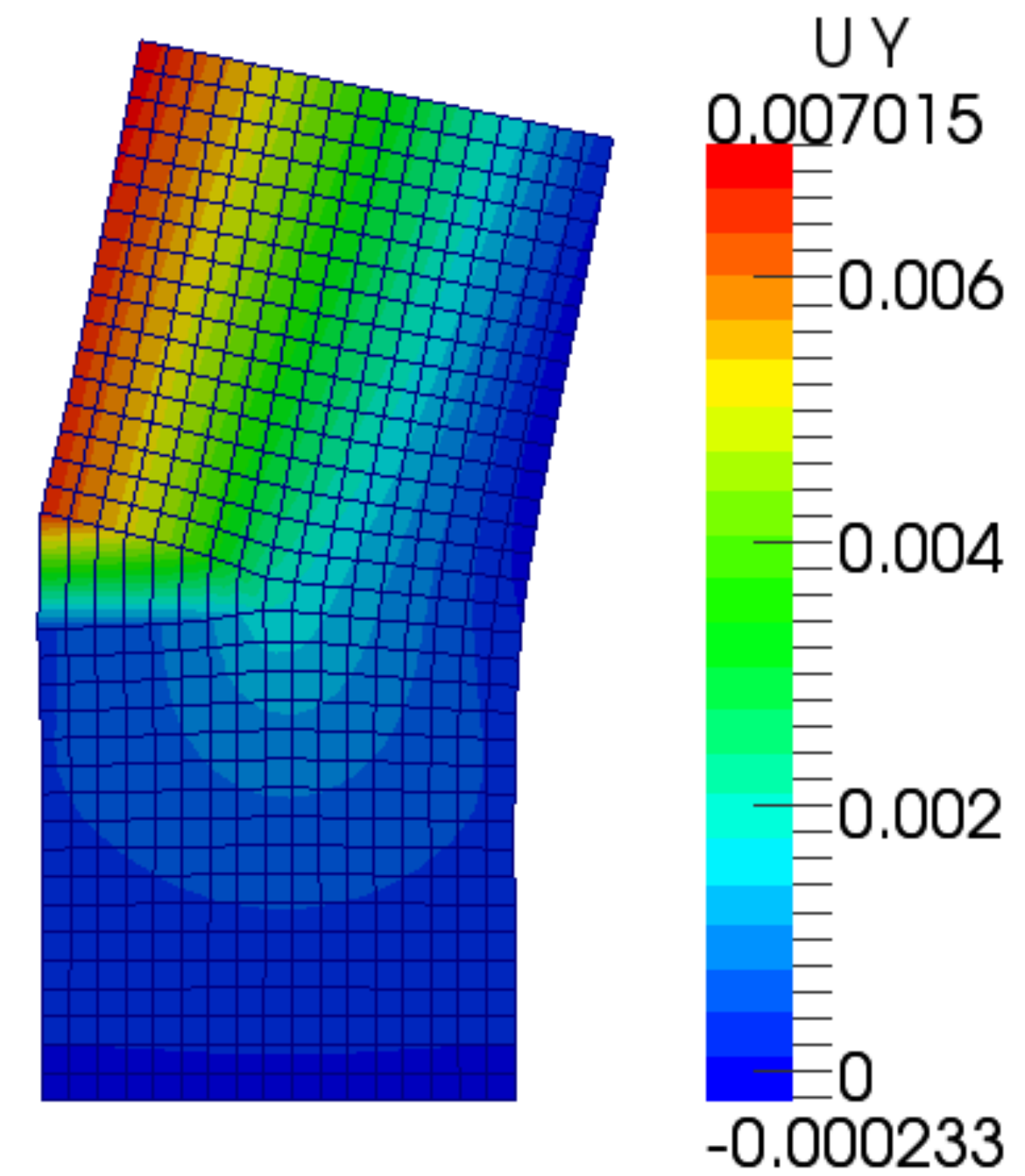}}\;\;\;\;
  \subfloat[XIGA]{\includegraphics[width=0.25\textwidth]{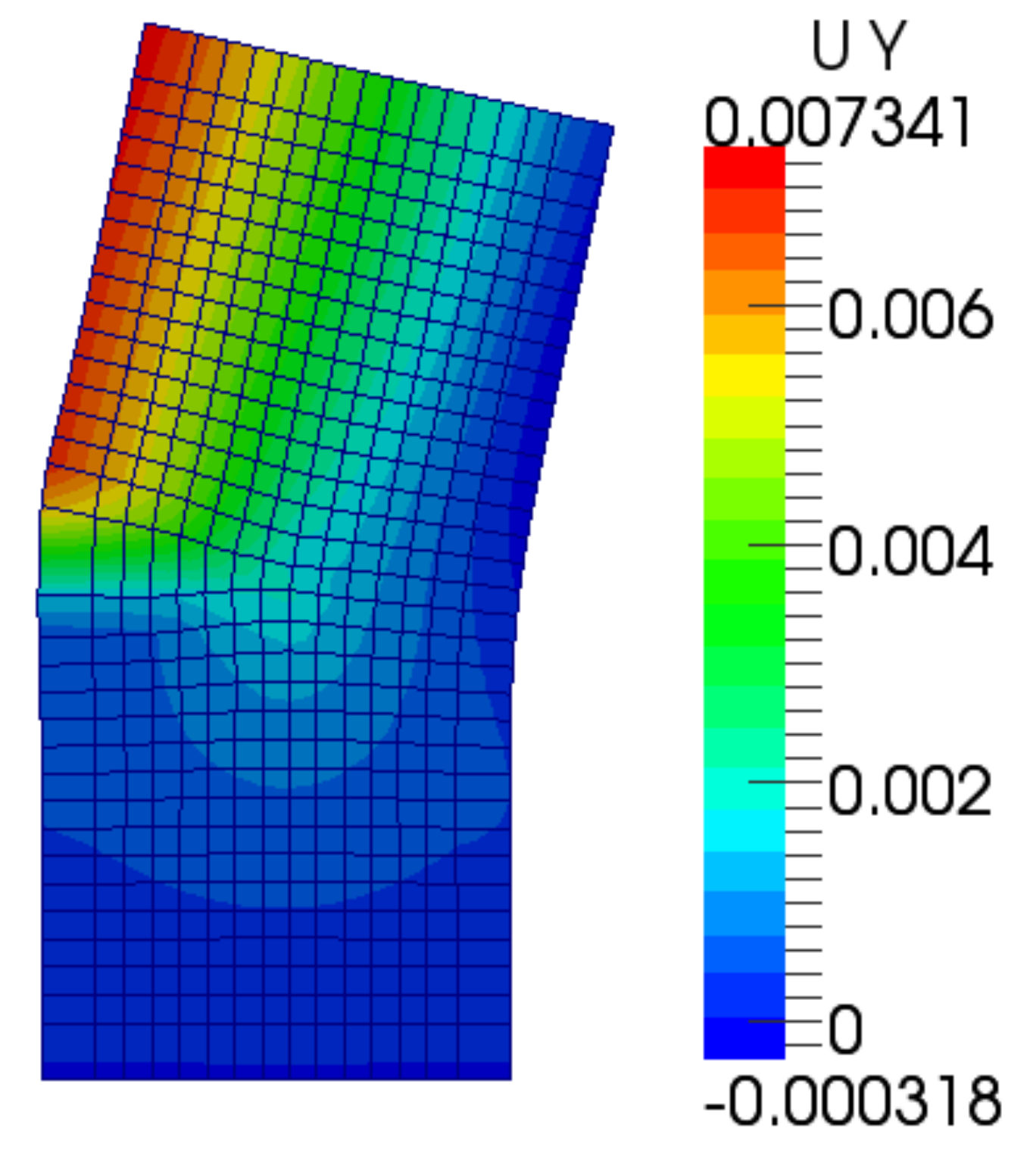}}
  \caption{Edge cracked plate: $u_y$ contour plots on deformed configuration
  (enlargement factor of 30 used).}
  \label{fig:edge-crack-deforms} 
\end{figure}

We first verify the implementation of the XIGA code by comparing the XIGA
result with the XFEM result for the case of $a=0.45$, $b=1$ and $h=1$. The
XFEM and XIGA meshes are given in Fig. \ref{fig:edge-crack-meshes}.
Both meshes have the same uniform distribution of nodes/control points. 
For the XFEM mesh, bilinear Q4 elements are used. For the XIGA mesh,
cubic ($p=q=3$) B-spline basis functions are adopted.
Fig. \ref{fig:edge-crack-deforms} shows the contour plots of the vertical
displacement obtained with XFEM and XIGA.

\begin{figure}[htbp]
  \centering 
  \includegraphics[width=0.4\textwidth]{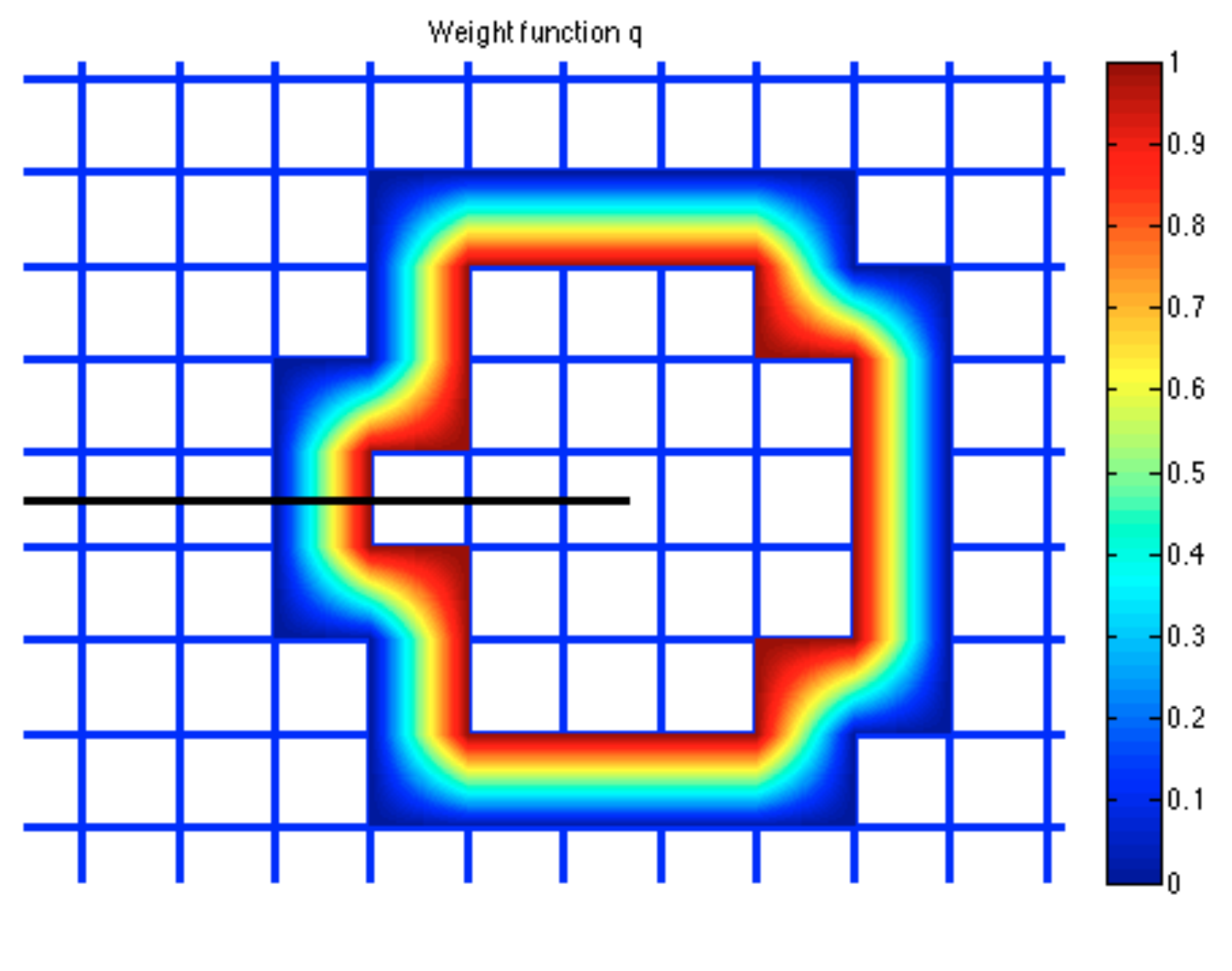}
  \caption{Distribution of weight function used in the computation of 
          the interaction integral. Four-noded quadrilateral elements
	  with bilinear Lagrange shape functions are used to
         interpolate the weight
	  function.}
  \label{fig:weight-function}
\end{figure}

We now consider the computation of the mode I SIF for a crack of
length $a=0.3$. The reference SIF for this problem is $K_I^\text{ref}=1.6118$. Both a linear and
a cubic B-spline basis are used for three different meshes. The results are given in 
Table \ref{tab:edge-crack}. It should be emphasized that in the computation of
the interaction integral, we use bilinear Lagrange shape functions \ie shape
functions of Q4 elements to compute the 
derivatives of the weight function. This guarantees that the weight function 
takes a value of unity on an open set containing the crack tip and 
vanishes on an outer contour as shown in Fig. \ref{fig:weight-function}.
B-splines functions are not interpolatory and therefore cannot be used to approximate
the weight functions.

\begin{table}[htbp]
  \centering
  \begin{tabularx}{\textwidth}{CCCCCC}
    \toprule
    mesh & disp. dofs & $K_I$ (linear) & Error (\%) & $K_I$ (cubic) & Error (\%)\\
    \midrule
    $9\times18$  & 324  & 1.4997 & 6.96 & 1.5560 & 3.46\\
    $18\times36$ & 1296 & 1.5823 & 1.83 & 1.6150 & 0.20\\
    $36\times72$ & 5184 & 1.5968 & 0.93 & 1.6117 & 0.01\\
    \bottomrule
  \end{tabularx}
  \caption{Edge cracked plate: SIFs results. The reference SIF is $K_I^\text{ref}=1.6118$.
  Note that linear NURBS are equivalent to the conventional bilinear finite elements.}
  \label{tab:edge-crack}
\end{table}

\subsubsection{Three-dimensional mode I fracture problem}

This example aims to show the capability of MIGFEM
for solving three-dimensional (3D) fracture problems. For 3D cracks,
the polar coordinates in the branch functions given in Eq. (\ref{branch}) 
are defined in terms of the level sets as \cite{xfem-3d-abaqus}

\begin{equation}
	r = \sqrt{\varphi(\xi,\eta,\zeta)^2 + \psi(\xi,\eta,\zeta)^2}, \quad
	\theta =
	\text{atan}\left(\frac{\varphi(\xi,\eta,\zeta)}{\psi(\xi,\eta,\zeta)}\right),
\end{equation}
\noindent where the level set field $\Phi=(\varphi,\psi)$ is interpolated as

\begin{equation}
	\Phi(\xi,\eta,\zeta) = \sum_I R_I(\xi,\eta,\zeta) \Phi_I. 
\end{equation}
We refer to \cite{xfem-3d-abaqus} for details concerning the derivatives of 
the branch functions with respect to the parametric coordinates
$(\xi,\eta,\zeta)$.

\begin{figure}[htbp]
  \centering 
  \psfrag{10}{10}
  \psfrag{10}{10}
  \psfrag{x}{$x$}
  \psfrag{y}{$y$}
  \psfrag{z}{$z$}
  \includegraphics[width=0.4\textwidth]{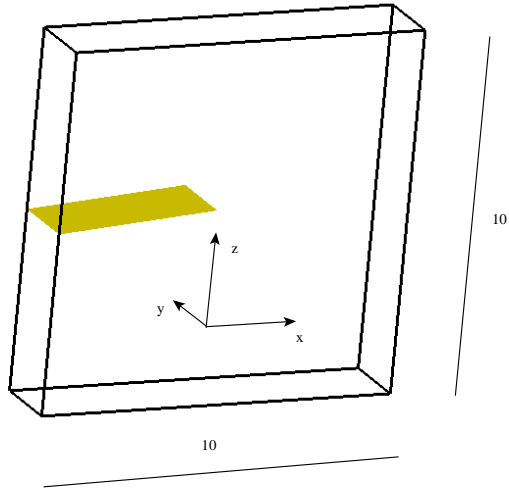}
  \caption{Three-dimensional mode I fracture problem: infinite plate
  with a center planar crack. The plate thickness is 2, the crack length is 5
  and the crack width is 2. The crack is located in the mid-plane of the plate.}
	  \label{fig:infiniteCrack3d-geo}
\end{figure}

The mode I 3D fracture problem we are solving is given in Fig.
\ref{fig:infiniteCrack3d-geo}. The exact displacement field is given by

\begin{equation}
\begin{split}
 u_x(r,\theta) &=\D{\frac{2(1+\upsilon)}{\sqrt{2\pi}}}\D {\frac{K_I}{E}}\sqrt{r}\cos\frac{\D \theta}{\D 2}
 \biggl({2-2\upsilon - \cos
^2\D{\frac{\theta}{2}}}\biggr)\\
u_y(r,\theta) &= 0\\
 u_z(r,\theta) &=\D{\frac{2(1+\upsilon)}{\sqrt{2\pi}}}\D {\frac{K_I}{E}}\sqrt{r}\sin\frac{\D \theta}{\D 2}
 \biggl({2-2\upsilon - \cos
^2\D{\frac{\theta}{2}}}\biggr),
 \end{split}
 \label{eq:exactDisp3D}
\end{equation}
 \noindent
where $K_I=\sigma\sqrt{\pi a}$ is the stress intensity factor,
$\upsilon$ is Poisson's ratio and $E$ is Young's modulus. 
In our example,  $a=100$ mm; $E = 10^7$ $\textrm{N/mm}^2$, $\upsilon = 0.3$,
$\sigma = 10^4$ $\textrm{N}/\textrm{mm}^2$. On the bottom, right and top
surfaces, essential BCs taken from Eq. (\ref{eq:exactDisp3D}) are imposed
using the penalty method. A penalty parameter of $1e10$ was used.
We note that this problem can be
more effectively solved with two-dimensional elements. This example however
aims at presenting how 3D extended IGA can be implemented. Furthermore, it
also illustrates how Dirichlet BCs are enforced on surfaces rather than the
usual case of line boundaries. To this end, a two-dimensional NURBS mesh for a
given surface is generated from the set of control points that define this surface
(see the file \textbf{surfaceMesh.m}).

\begin{figure}[htbp]
  \centering 
  \includegraphics[width=0.6\textwidth]{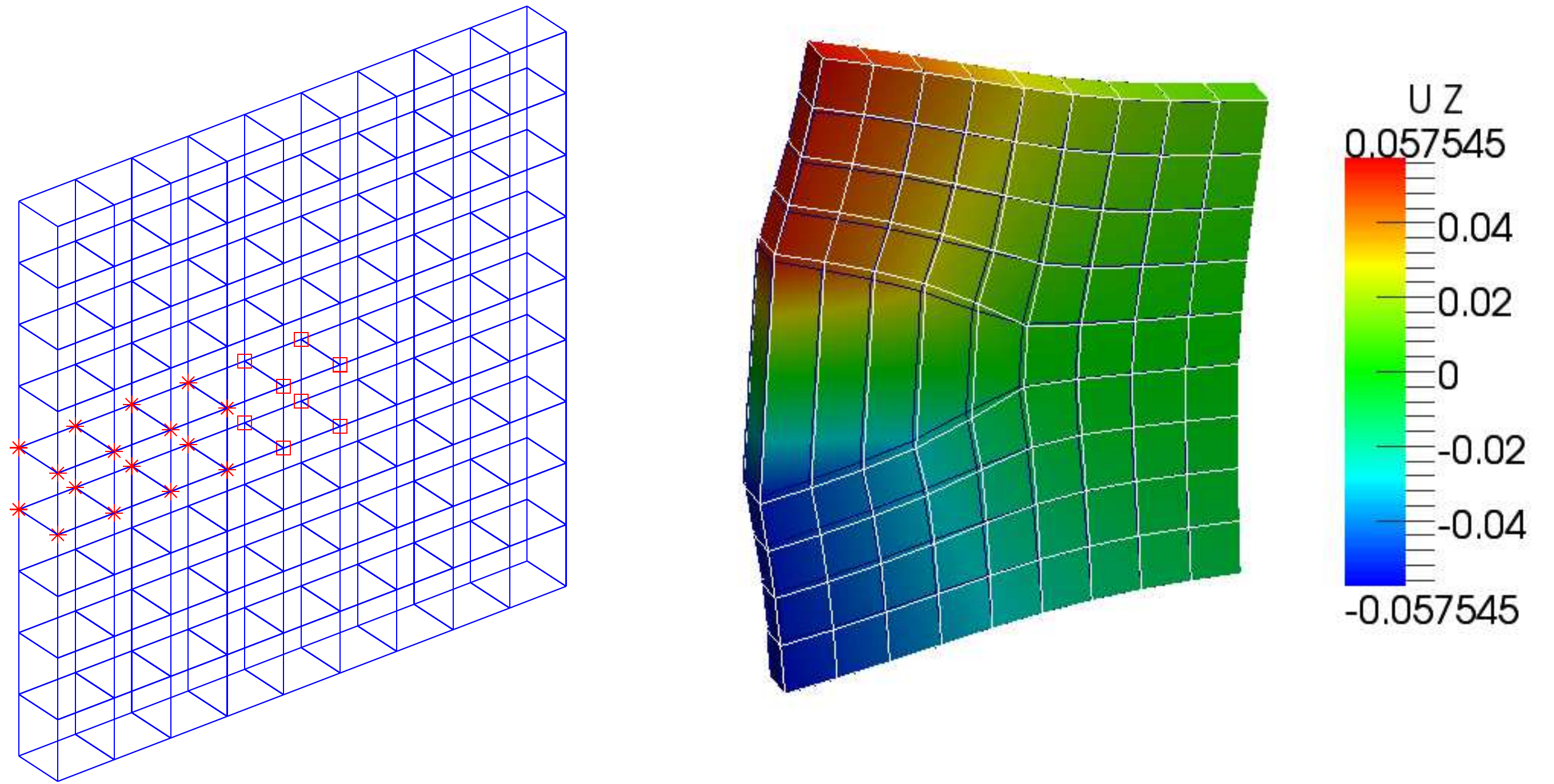}
  \caption{Three-dimensional mode I fracture problem: mesh of linear 
  B-spline elements and enriched control points (left); numerical deformed configuration
  (magnification factor of 40) superimposed on the exact deformed
  configuration (right).}
	  \label{fig:infiniteCrack3dC0}
\end{figure}

\begin{figure}[htbp]
  \centering 
  \includegraphics[width=0.7\textwidth]{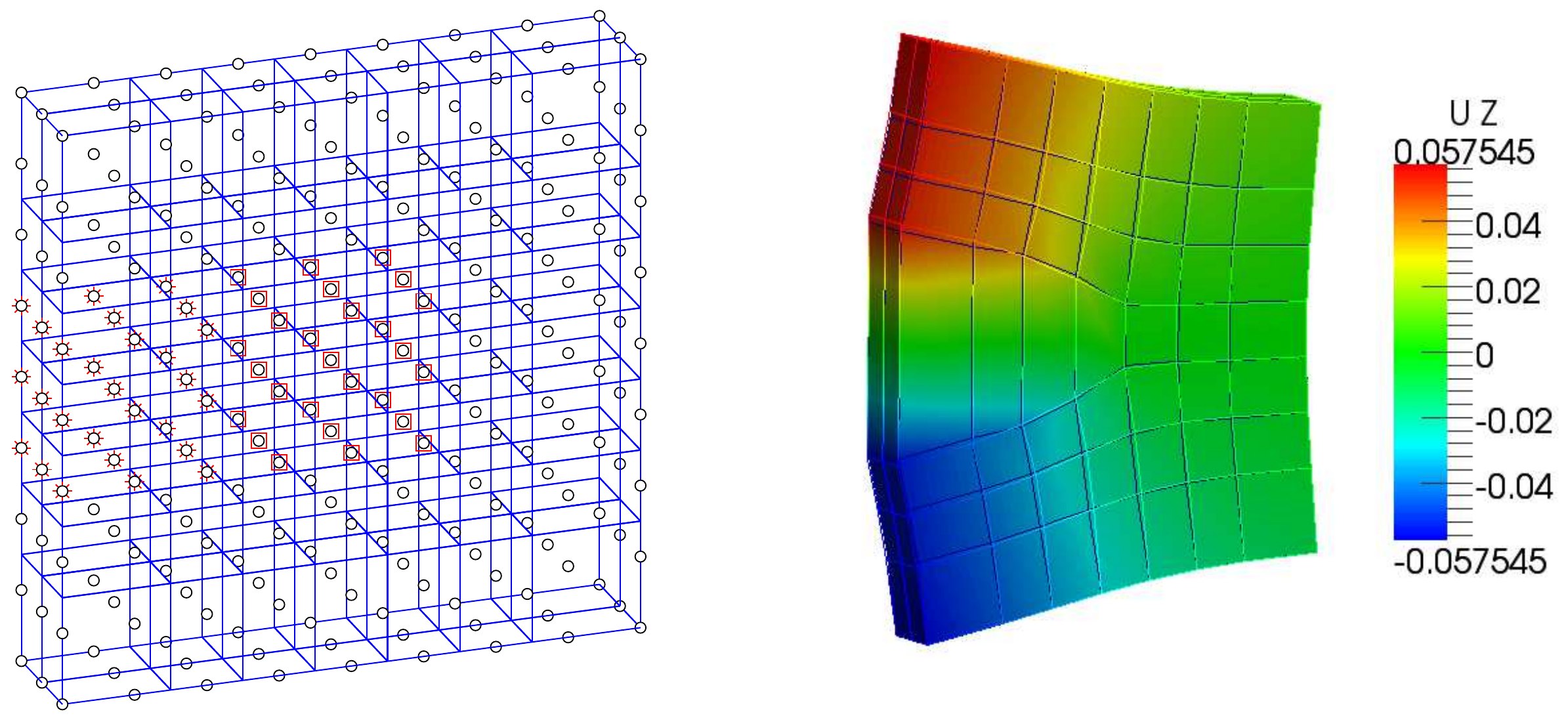}
  \caption{Three-dimensional mode I fracture problem: mesh of quadratic 
  B-spline elements and enriched control points (left); numerical deformed configuration
  (magnification factor of 40) superimposed on the exact deformed
  configuration (right).}
	  \label{fig:infiniteCrack3dC1}
\end{figure}

The problem is first solved using a linear B-spline basis. A mesh of
$9\times9\times1$ elements is used. The mesh, enriched nodes and comparison
of the numerical deformed configuration against the exact profile are given in 
Fig. \ref{fig:infiniteCrack3dC0}. Next, a mesh of $7\times7\times2$ elements is used 
where, in the through-thickness direction, there are two linear elements ($q=1$) and for the two
other directions, a quadratic basis ($p=r=2$) is used. The result is given
in Fig. \ref{fig:infiniteCrack3dC1} and we note a good qualitative agreement
between Figs. \ref{fig:infiniteCrack3dC0} and \ref{fig:infiniteCrack3dC1}. We 
did not perform a SIF computation for this problem as 3D SIF computation is not yet implemented
at present.

\subsection{Structural mechanics}

\begin{figure}[h!]
  \centering 
  \includegraphics[width=0.55\textwidth]{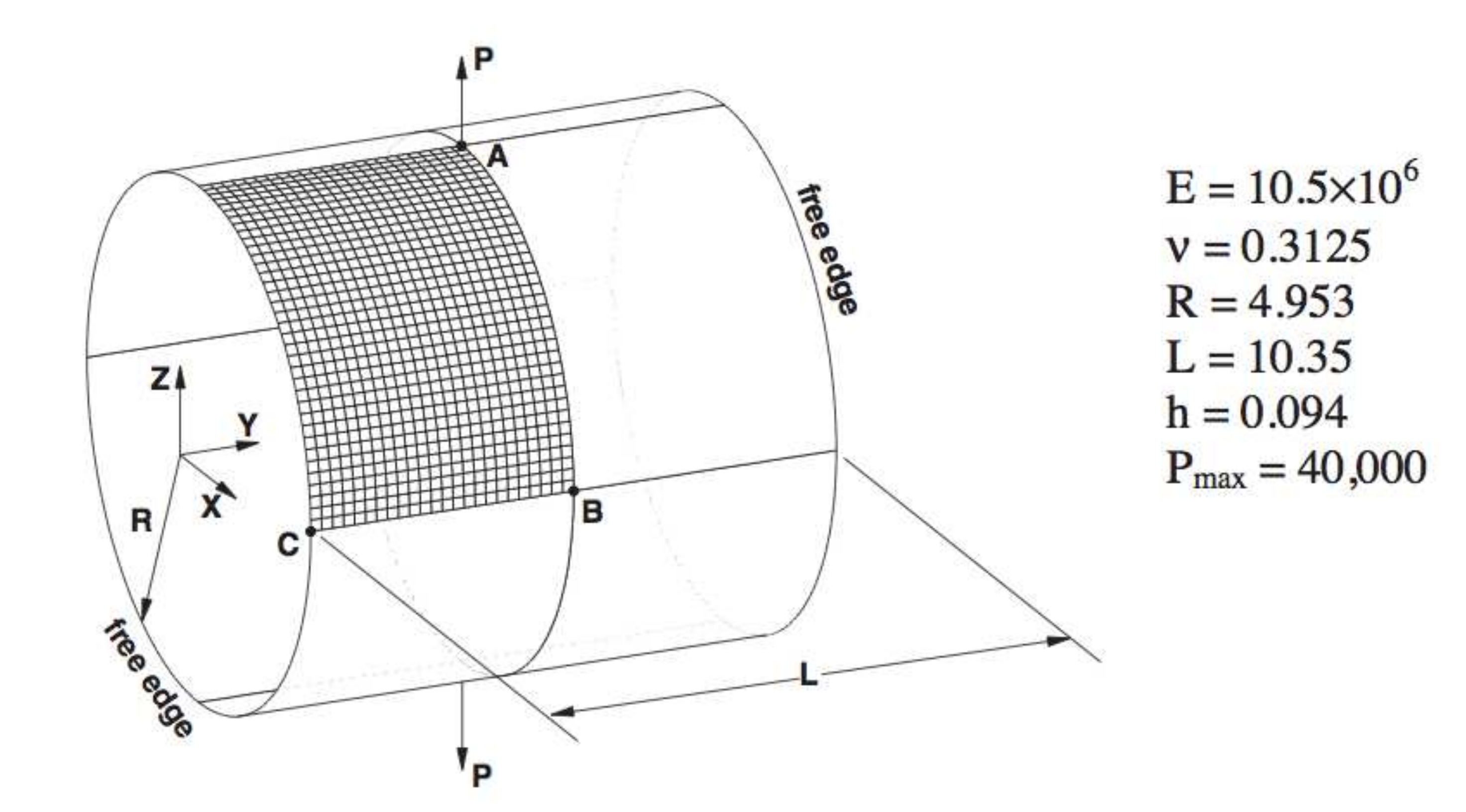}
  \caption{The open-end cylindrical shell subjected to radial pulling forces: problem setup \cite{Sze20041551}.}
  \label{fig:freeEndsGeo} 
\end{figure}

\begin{figure}[h!]
  \centering 
  \psfrag{sym1}[c]{sym ($u_z=0$, zero rotation)}
  \psfrag{sym2}[c]{sym ($u_y=0$, zero rotation)}
  \psfrag{sym3}[c]{sym ($u_x=0$, zero rotation)}
  \psfrag{free}{free}
  \psfrag{f}{$P/4$}
  \psfrag{ux}{$u_y=0$}
  \psfrag{uxx}{$u_x=0$}
  \psfrag{uxxx}{$u_z=0$}
  \psfrag{ui}{$\vm{u}_I$}
  \psfrag{ui1}{$\vm{u}_{I1}$}
  \includegraphics[width=0.36\textwidth]{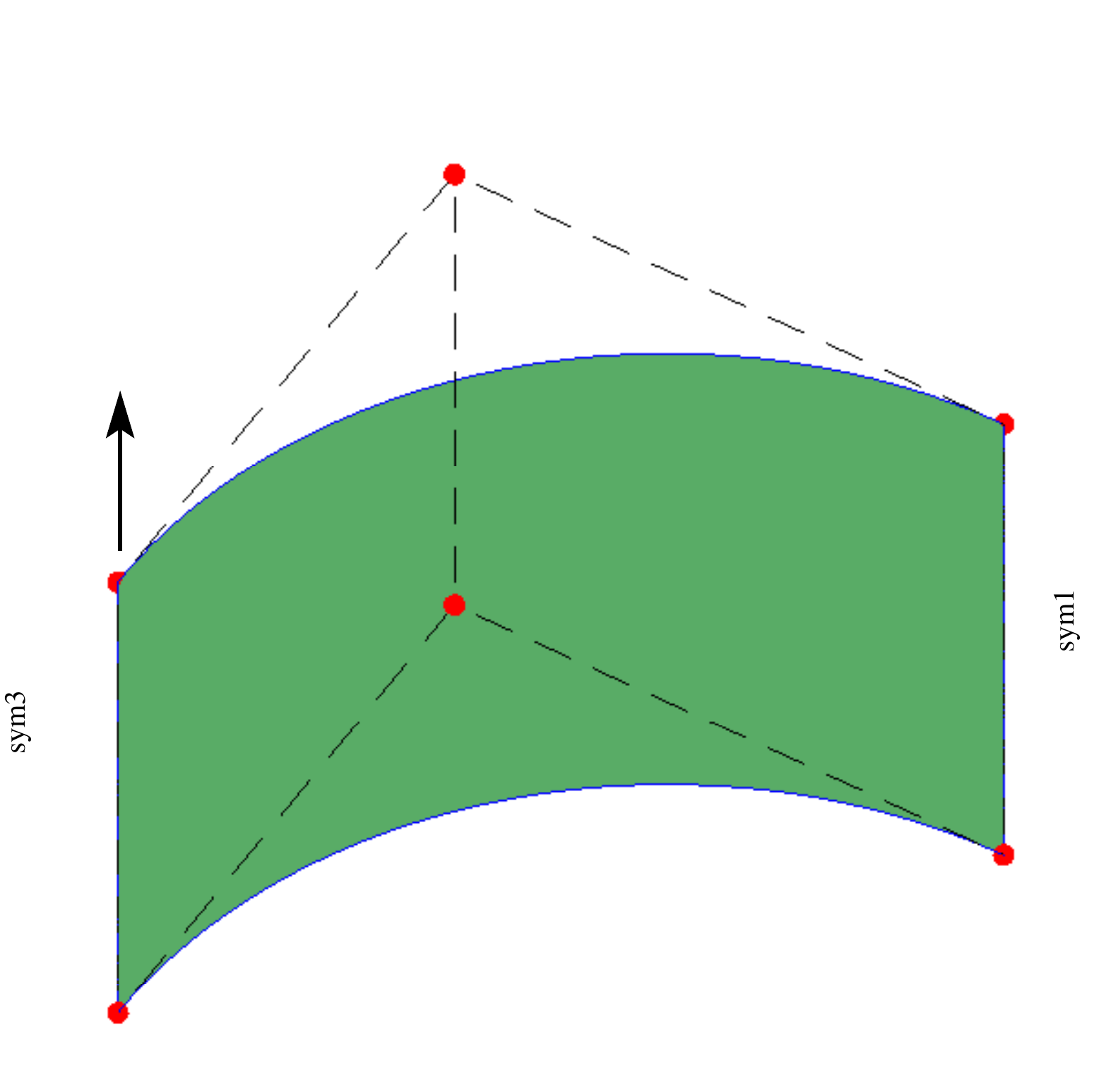}\;\;\;
  \includegraphics[width=0.36\textwidth]{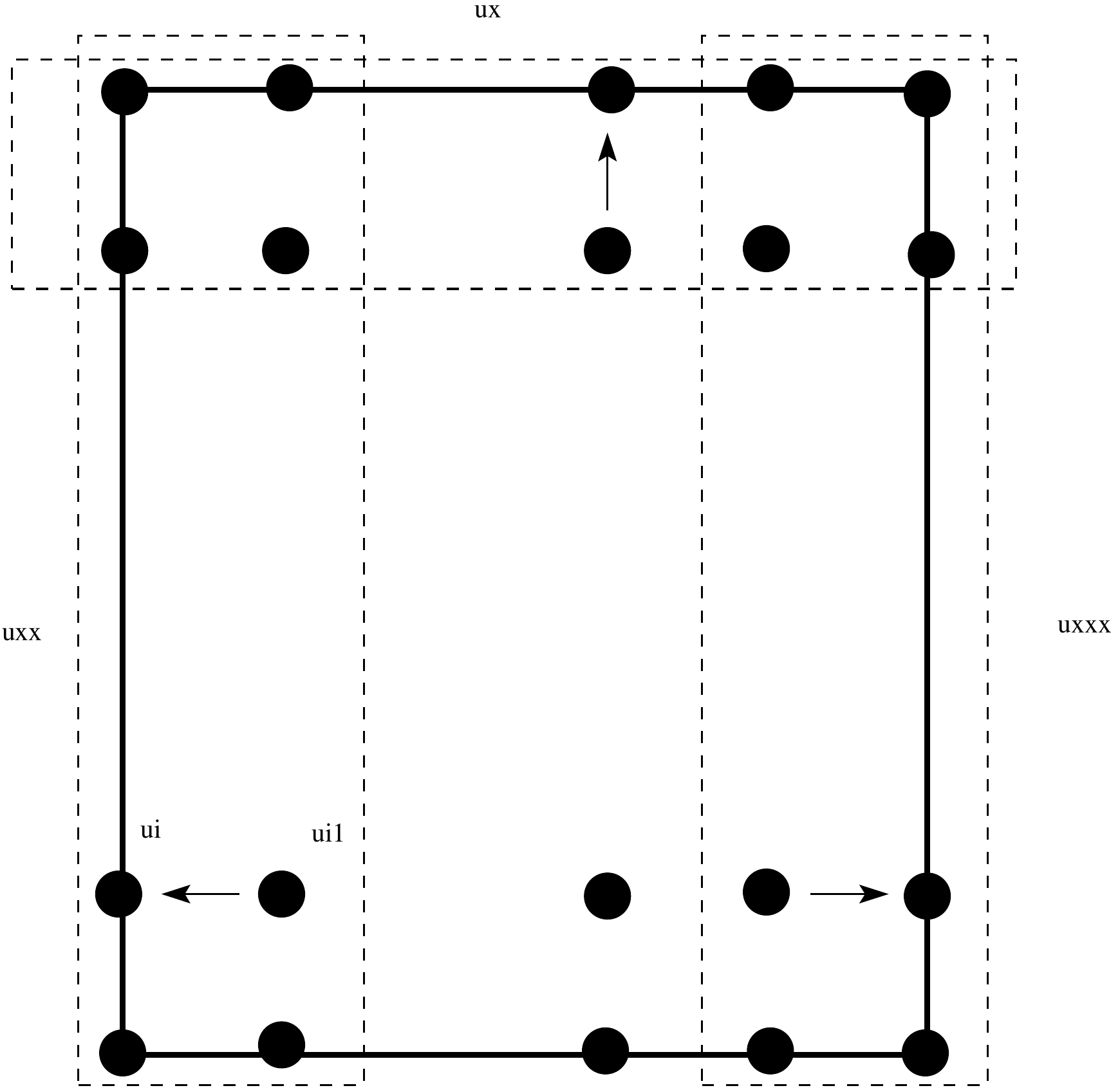}
  \caption{Open ended cylindrical shell: 1/8 of the model with one bi-quadratic NURBS surface.
          The control points and weights are given in file \textbf{freeEndsGNLCylinderShellData.m}.
          For symmetry edges, in order to satisfy the
          symmetry condition \ie zero rotation, we constraint two rows of control points in the sense
          that the displacements of the control points right next to the control points locating 
          on the symmetry edges ($\vm{u}_{I1}$) matches $\vm{u}_I$.}
  \label{fig:freeEnd-mesh} 
\end{figure}

We consider the pull out of an open-ended cylindrical shell, see Fig. \ref{fig:freeEndsGeo},
as one of the most common benchmarks for geometrically nonlinear thin shell problems
\cite{Sze20041551,NME:NME4403}.
Due to symmetry only 1/8 of the model is studied. The 1/8 model can be exactly described using one single quadratic-linear NURBS surface as shown in Fig.\ref{fig:freeEnd-mesh}. Refined meshes are then obtained from this initial mesh by using the $k$-refinement. The analysis was performed using a mesh of $8\times8$ 
bi-quartic NURBS elements (144 control points). 
The enforcement of symmetry BCs is achieved by constraining two row of 
control points as shown in the same figure.  These constraints--the so-called multipoint 
linear homogeneous constraints are handled using the master-slave method 
\cite{felippa:note} for the penalty method as described in Section \ref{sec-sym} could endanger the
convergence of the Newton-Raphson solver. Note that for this particular nonlinear 
problem, of which the M file is \textbf{igaGNLFreeEndCylinderShell.m}, the analysis was performed with
a C++ implementation \cite{jemjive}.

\begin{figure}[h!]
  \centering 
  \includegraphics[width=0.25\textwidth]{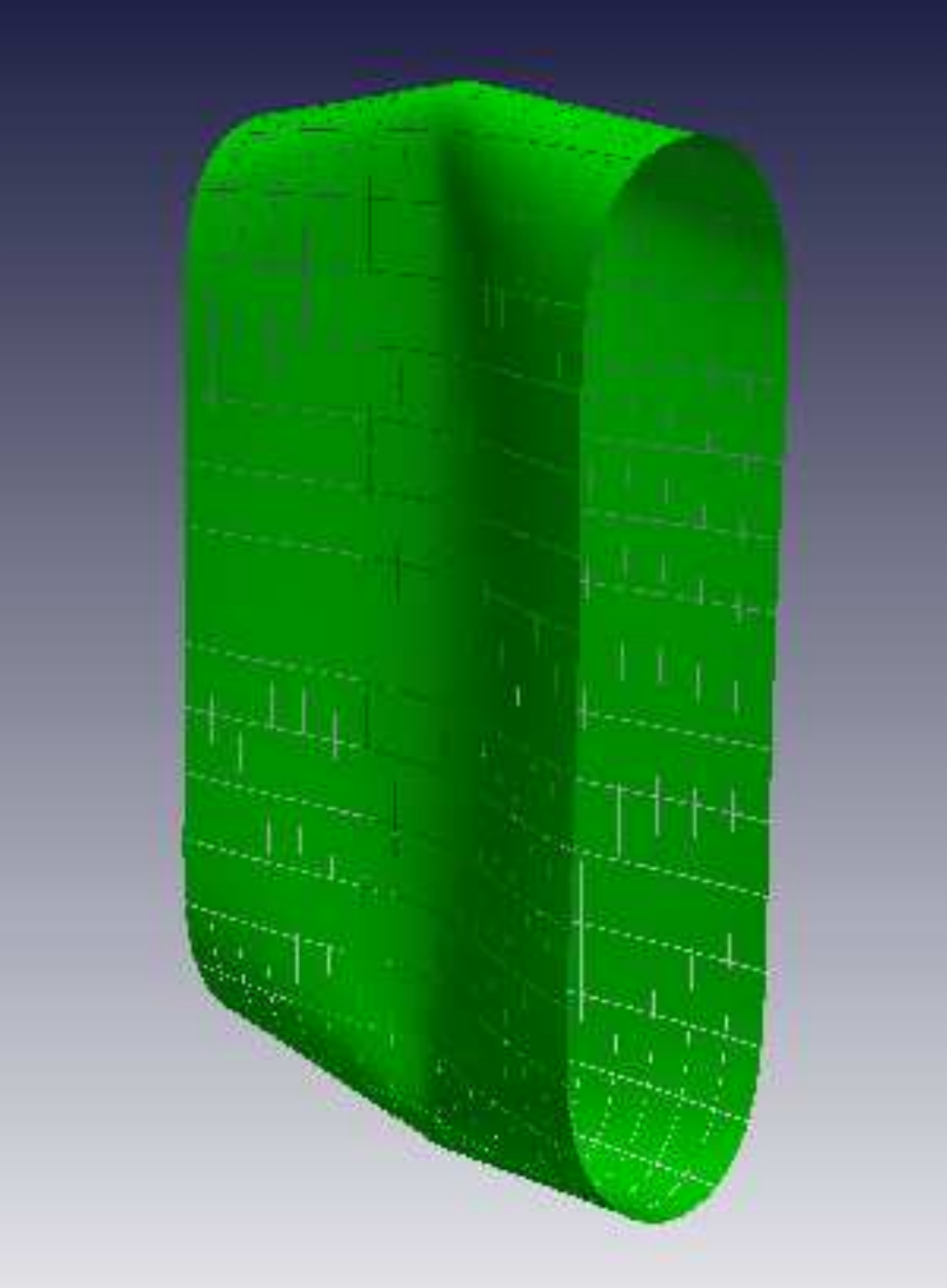}
  \caption{The open-ended cylindrical shell subjected to radial pulling forces: deformed shape without
     scaling.}
  \label{fig:freeEnds-deform} 
\end{figure}

\begin{figure}[h!]
  \centering 
  \includegraphics[width=0.5\textwidth]{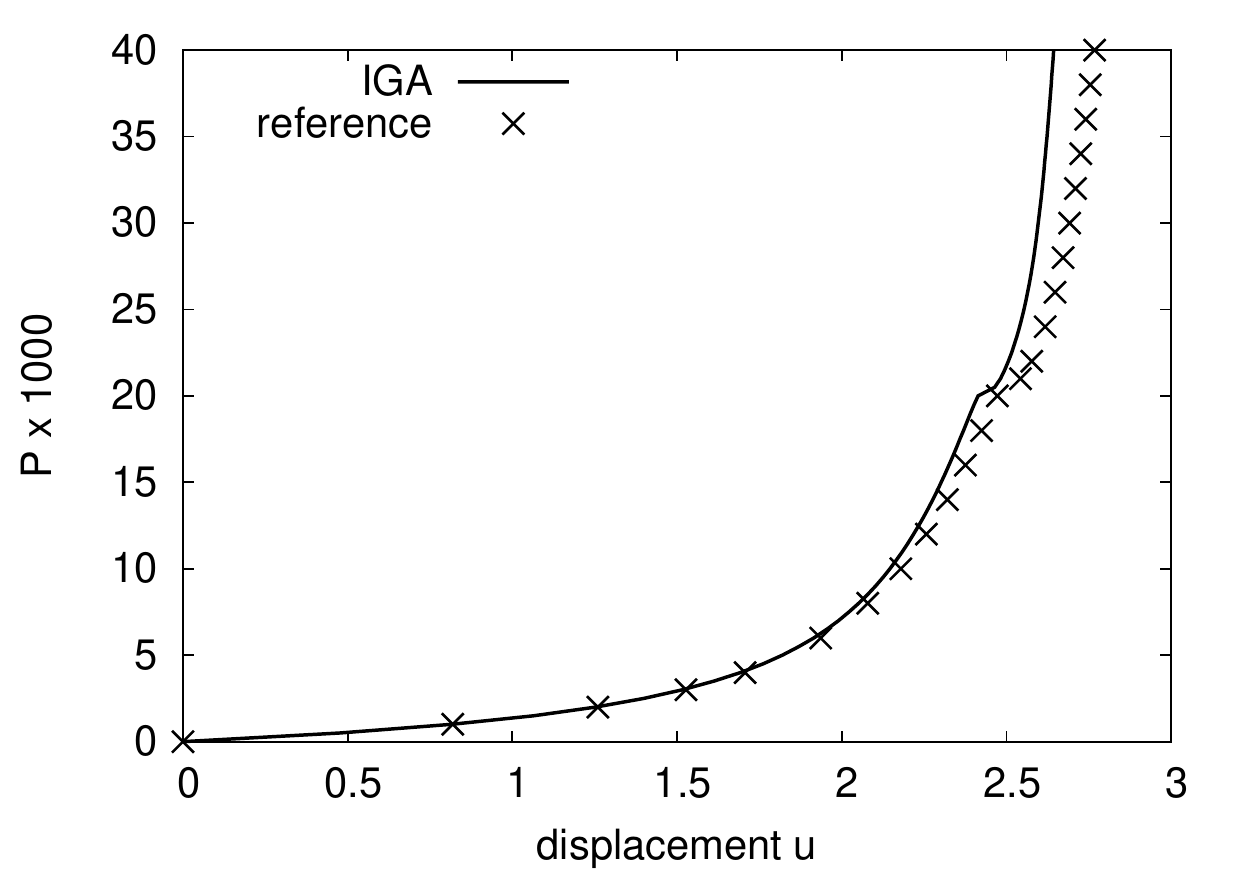}
  \caption{The open-ended cylindrical shell subjected to radial pulling forces: load-displacement responses.}
  \label{fig:freeEnds-lodi} 
\end{figure}

We adopt the Kirchhoff-Love thin shell model, see \eg 
\cite{Dung20082778,NME:NME4403} for details that involves only displacement dofs. The maximum applied load $0.25P_\text{max}$ is applied in 80 equal increments and for each increment the full Newton-Raphson method 
is used to solve for the displacements. 
Fig. \ref{fig:freeEnds-deform} shows the deformed configuration of the shell and 
in Fig.\ref{fig:freeEnds-lodi} the load versus the $z$-displacement at point A is plotted together with
the result reported in \cite{Sze20041551}. 

\section{Conclusion}\label{conclusions}

We presented a Matlab\textsuperscript{\textregistered{}} implementation for one, two and three-dimensional 
isogeometric finite element analysis for solid and structural mechanics. 
This paper is addressed to students or researchers who wish to learn
the concepts of IGA in a clear and concise manner and is especially
suited to those with solid mechanics applications in mind. NURBS are
used throughout, where the underlying construction of the basis
functions is detailed along with associated refinement algorithms
essential for numerical analysis. Differences with conventional FE
implementations are made clear with the use of Matlab\textsuperscript{\textregistered{}} source code to
illustrate isogeometric FE concepts explicitly. In addition, the
implementation of an isogeometric XFEM formulation for both
two-dimensional and three-dimensional problems is described
allowing for simple linear elastic fracture analysis to be performed directly from CAD data. 
The benefits of isogeometric analysis are made evident.
Although not presented, the code supports multi-patch analysis in which  
compatibility between connecting patches is not required \cite{nguyen-nitsche1,nguyen-nitsche2}. 
Geometrical nonlinearities for solid elements under the Total Lagrange framework are provided as well.
In addition, PUM enrichment for holes and inclusions is provided along with implementation of the least squares method for imposing essential boundary conditions. Mass matrices are implemented for almost every elements and popular time integration schemes such as Newmark and central difference explicit are available so that transient 
analysis on CAD objects can be performed.
The code is available for download from \url{https://sourceforge.net/projects/cmcodes/} for Linux and Mac OS
machines. A Windows version can be found at \url{http://www.softpedia.com/get/Science-CAD/igafem.shtml}.

The preliminary concepts and implementation details of isogeometric
analysis have been described, but many challenges remain in the
field. These include the creation of suitable volume discretisations
from given CAD boundary representations, efficient integration
schemes and suitable error estimators.

\section*{Acknowledgements}

The authors would like to acknowledge the partial financial support of the
Framework Programme 7 Initial Training Network Funding under grant number
289361 ``Integrating Numerical Simulation and Geometric Design Technology".
St\'{e}phane Bordas also thanks partial funding for his time provided by
1) the EPSRC under grant EP/G042705/1 Increased Reliability for Industrially
Relevant Automatic Crack Growth Simulation with the eXtended Finite Element
Method and 2) the European Research Council Starting Independent Research
Grant (ERC Stg grant agreement No. 279578) entitled ``Towards real time multiscale simulation of cutting in
non-linear materials with applications to surgical simulation and computer
guided surgery''. The first author would like to express his gratitude towards Professor
L.J. Sluys at Delft University of Technology, The Netherlands for his support to VPN during
the PhD period and the Framework Programme 7 Initial Training Network Funding. The authors acknowledge
the comments of Dr. Robert Simpson at Cardiff University that has improved the paper.

\appendix

\section{Knot vector conventions}
\label{sec:knot-vect-conv}

The definition of an open knot vector as $\Xi = \{ \xi_1, \xi_2,
\ldots,\xi_{n+p+1} \}$ where the first and last knot values are
repeated $p+1$ times is sometimes simplified to remove redundant
information. It is found that if a knot vector is open, the first and
last knot values play no role in the definition of the curve and can
be removed entirely. Knot vectors are then defined as $\Xi = \{ \xi_1, \xi_2,
\ldots,\xi_{n+p-1} \}$ where the first and last knot vectors are
repeated $p$ times. This notation is used frequently through the CAGD
community. 

Furthermore, B-spline algorithms can be written entirely in terms of knot
interval vectors \cite{Se10}, defined as 
\begin{equation}
  \label{eq:knot_inteval_defn}
  \Delta \Xi = \{ \Delta \xi_1, \Delta
  \xi_2,\ldots,\Delta\xi_{n+p-2}\} \quad \textrm{where} \quad
  \Delta\xi_i = \xi_{i+1} - \xi_i
\end{equation}
using the previous simplified knot vector notation. T-splines are
based entirely on knot interval vector notation.

\section{$C^k$ and $G^k$ continuity}
\label{sec:disc-cont}

Two types of continuity are commonplace in the CAGD community,
referred to as $C$ and $G$ continuity. From a numerical analysis
perspective, $C$ continuity refers to the traditional definition of
continuity in which, given a parametric function $f: \mathbb{R} \to \mathbb{R}^2$ which defines a 2D
curve through a parameter $t$, the parameterisation is said
to be $C^k$ continuous at point $a$ if $d^nf/dt^n$ is
continuous (single-valued) for $n=0,1,\ldots,k$ at $t=a$. Intuitively, for
a curve to be $C^0$, it must contain no discontinuities with $t$ 
single-valued throughout its domain. Letting $df/dt$ represent the `speed' 
of the curve, $C^1$ continuity implies that the tangent vector
is equal in magnitude and direction while also assuming $C^0$
continuity. $C^2$ continuity assumes both $C^0$ and $C^1$ continuity
and requires that the `acceleration' of the curve is equal in
magnitude and direction.

$G$ continuity (sometimes referred to as `visual' continuity) is
independent of the parameter $t$. For a curve to be $G^0$ continuous at a
point it must contain no discontinuities, but the parameter value may
have multiple values. For $G^1$ continuity, the tangent of the curve
must be equal in its direction, but its magnitude may differ. And
$G^2$ continuity implies that the direction of acceleration is
continuous, but its magnitude may differ.

If desired, a $G^k$ curve can be made $C^k$-continuous through an
appropriate `reparameterisation'.

\section{Homogenous and non-homogenous coordinates}
\label{sec:homog-non-homog}

It is found that if rational basis functions are used, then it is often more appropriate to work in terms of homogeneous coordinates, since they can be applied straightforwardly to non-rational basis function algorithms. For example, in the case of NURBS, the use of homogenous coordinates allows B-spline algorithms to be used directly.

Given a physical coordinate $\mathbf{P}_A=(x_A, y_A, z_A)$ with weighting $w_A$, the corresponding four-dimensional homogeneous coordinate $\tilde{\mathbf{P}}_A$ is written as
\begin{align}
  \label{eq:homo_coord}
\tilde{\mathbf{P}}_A &= ( x_A w_A, y_A w_A, z_A w_A, w_A )\\
&=  (  \tilde{x}_A, \tilde{y}_A, \tilde{z}_A, w_A ).
\end{align}
To convert back to non-homogenous coordinates, it is a simple case of dividing $\tilde{x}_A, \tilde{y}_A, \tilde{z}_A$ through by $w_A$.

\section{B\'{e}zier extraction}\label{sec:extraction}

B\'{e}zier extraction \cite{borden_isogeometric_2011,scott_isogeometric_2011} provides 
a tool to facilitate the incorporation of NURBS and even T-splines in to any FE codes. 
The idea is based on the fact that any B-splines basis can be written as a linear 
combination of  Bernstein polynomials. Written for element $e$, the shape functions read

\begin{equation}
\vm{N}^e(\bsym{\xi}) = \vm{C}^e \vm{B}(\tilde{\bsym{\xi}})
\label{eq:extraction}
\end{equation}
\noindent where $\vm{C}^e$ denotes the elemental B{\'e}zier extraction operator
and $\vm{B}$ are the Bernstein polynomials which are defined on
the parent element $\tilde{\Omega}$. Index space, knot vectors are all embedded in 
the B{\'e}zier extrators which are computed in a pre-processing step. Therefore, existing
FE solvers can use NURBS/T-splines straightforwardly.
We refer to folder \textbf{bezier-extraction} for examples.

\bibliography{isogeometric}
\bibliographystyle{unsrt}


\end{document}

%% file: pre_commands.tex
\makeatletter
\renewcommand\fs@ruled{\def\@fs@cfont{\bfseries}\let\@fs@capt\floatc@ruled
\def\@fs@pre{\hrule height 1.2pt depth0pt \kern2pt}%
\def\@fs@post{\kern2pt\hrule height 1.2pt depth0pt \kern2pt \relax}%
\def\@fs@mid{\kern2pt\hrule\kern2pt}%
\let\@fs@iftopcapt\iftrue}
\makeatother

\newcommand{\norm}[1]{\left|\left|#1\right|\right|}

\newcommand{\cref}[1]{Chapter~\ref{chapt:#1}}

\newcommand{\D}{\displaystyle}

\newcommand{\di}{\mathrm{d}}

\usepackage{xspace}
\newcommand{\eg}{e.g.,\xspace}

\newcommand{\ie}{i.e.,\xspace}

\newcommand{\etc}{etc.\@\xspace}

\newcommand{\cmatrixb}{\left\{ \begin{matrix}}
\newcommand{\cmatrixe}{\end{matrix} \right\}}

\newcommand{\vm}[1]{\mathbf{#1}}

\newcommand{\bsym}[1]{\boldsymbol{#1}}

\newcommand{\bx}{\boldsymbol{x}}

\newcommand{\trans}{^\mathrm{T}}

\newcommand{\bc}{\begin{center}}
\newcommand{\ec}{\end{center}}
\newcommand{\bitem}{\begin{itemize}}
\newcommand{\eitem}{\end{itemize}}

\newcommand{\ljump}{\lbrack \! \lbrack } 
\newcommand{\rjump}{\rbrack \! \rbrack } 
\newcommand{\jump}[1]{\ljump {#1} \rjump} 

\newcommand{\pderiv}[2]{ \frac{\partial {#1} }{\partial {#2} } }

\newcommand{\beq}{\begin{equation}}
\newcommand{\eeq}{\end{equation}}
\newcommand{\beqa}{\begin{eqnarray}}
\newcommand{\eeqa}{\end{eqnarray}}

\newcommand{\invisible}[1]{}

\newcommand{\bv}{\begin{verbatim}}
\newcommand{\V}{\verb}                  % EX: \V=-d{#@~}= Expr must

\newcommand{\testpix}[1]{\fbox{\begin{minipage}[c]{\textwidth}
                      #1 \end{minipage} }}

\newcommand{\putpstex}[1]{\includegraphics{#1.pstex_t}}